\documentclass{amsart}
\usepackage{amssymb,amsfonts}
\usepackage{mathabx}
\usepackage{enumerate}
\usepackage{mathrsfs}
\usepackage{graphicx}
\usepackage{mathtools}
\usepackage{leftidx}
\usepackage{tikz-cd}
\usepackage{bm}
\usepackage{url}
\newtheorem{thm}{Theorem}[section]
\newtheorem{cor}[thm]{Corollary}
\newtheorem{prop}[thm]{Proposition}
\newtheorem{lem}[thm]{Lemma}

\theoremstyle{definition}
\newtheorem{defn}[thm]{Definition}

\newtheorem{exmp}[thm]{Example}

\theoremstyle{remark}
\newtheorem{rem}[thm]{Remark}

\makeatletter
\let\c@equation\c@thm
\makeatother
\numberwithin{equation}{section}

\def\bthm{\begin{thm}}
\def\ethm{\end{thm}}
\def\blm{\begin{lem}}
\def\elm{\end{lem}}
\def\bdf{\begin{defn}}
\def\edf{\end{defn}}
\def\bpf{\begin{proof}}
\def\epf{\end{proof}}
\def\bpp{\begin{prop}}
\def\epp{\end{prop}}
\def\bcor{\begin{cor}}
\def\ecor{\end{cor}}
\def\brm{\begin{rem}}
\def\erm{\end{rem}}
\def\beg{\begin{exmp}}
\def\eeg{\end{exmp}}
\def\bK{\mathbb{K}}

\def\bN{\mathbb{N}}

\def\bQ{\mathbb{Q}}
\def\bR{\mathbb{R}}
\def\bS{\mathbb{S}}

\def\bZ{\mathbb{Z}}
\def\cA{\mathcal{A}}
\def\cB{\mathcal{B}}
\def\cC{\mathcal{C}}
\def\cD{\mathcal{D}}
\def\cE{\mathcal{E}}
\def\cF{\mathcal{F}}

\def\cM{\mathcal{M}}

\def\cP{\mathcal{P}}
\def\cQ{\mathcal{Q}}
\def\cR{\mathcal{R}}
\def\cS{\mathcal{S}}

\def\cV{\mathcal{V}}

\def\scO{\mathscr{O}}

\def\frF{\mathfrak{F}}
\def\frG{\mathfrak{G}}

\def\frf{\mathfrak{f}}

\newcommand{\raq}{\,\rightarrow \,}

\newcommand{\xraq}[2][]{\, \xrightarrow[#1]{#2} \,}
\newcommand{\xlaq}[2][]{\, \xleftarrow[#1]{#2} \,}

\newcommand{\ra}{\rightarrow}
\newcommand{\la}{\leftarrow}
\newcommand{\rinto}{\hookrightarrow}

\newcommand{\ronto}{\twoheadrightarrow}

\newcommand{\xra}[2][]{\xrightarrow[#1]{#2}}

\makeatletter
\newcommand{\xdashrightarrow}[2][]{\ext@arrow 0359\rightarrowfill@@{#1}{#2}}
\newcommand{\xdashleftarrow}[2][]{\ext@arrow 3095\leftarrowfill@@{#1}{#2}}
\newcommand{\xdashleftrightarrow}[2][]{\ext@arrow 3359\leftrightarrowfill@@{#1}{#2}}
\def\rightarrowfill@@{\arrowfill@@\relax\relbar\rightarrow}
\def\leftarrowfill@@{\arrowfill@@\leftarrow\relbar\relax}
\def\leftrightarrowfill@@{\arrowfill@@\leftarrow\relbar\rightarrow}
\def\arrowfill@@#1#2#3#4{%
	$\m@th\thickmuskip0mu\medmuskip\thickmuskip\thinmuskip\thickmuskip
	\relax#4#1
	\xleaders\hbox{$#4#2$}\hfill
	#3$%
}
\makeatother

\newcommand{\Set}{{\rm Set}}

\newcommand{\Mod}{{\rm Mod}}

\newcommand{\Sym}{{\rm Sym}}

\newcommand{\Ab}{{\rm Ab}}

\newcommand{\Der}{{\rm Der}}

\newcommand{\Ch}{{\rm Ch}}
\newcommand{\Setdel}{{\rm Set}_{\Delta}}

\newcommand{\Tot}{{\rm Tot}}

\newcommand{\RHom}{{\bm R}{\rm Hom}}

\newcommand{\Ob}{\rm{Ob}}
\newcommand{\op}{{\rm op}}

\newcommand{\id}{{\rm id}}

\newcommand{\Hom}{{\rm Hom}}
\newcommand{\Homcom}{\underline{{\rm Hom}}}
\newcommand{\End}{{\rm End}}
\newcommand{\Endcom}{\underline{{\rm End}}}

\newcommand{\Mor}{{\rm Mor}}

\newcommand{\colim}{{\rm colim}}

\newcommand{\cHom}{\mathscr{H}\text{\kern -3pt {\calligra\large om}}\,}

\newcommand{\gr}{{\rm gr}}

\newcommand{\OP}{{\rm Op}}
\newcommand{\COP}{{\rm Op}^c}

\newcommand{\lax}{{\rm lax}}

\newcommand{\Nm}{{\rm Nm}}

\newcommand{\Cinf}{\cC_{{\rm inf}}}

\newcommand{\GrMod}{{\rm GrMod}}

\newcommand{\PCOM}{{\rm Com}^{{\rm pre}}}
\newcommand{\COM}{{\rm Com}}

\newcommand{\MC}{{\rm MC}}

\newcommand{\Tw}{{\rm Tw}}

\newcommand{\Endpr}{\underline{\End}^{{\rm pr}}}
\newcommand{\Endrd}{\underline{\End}^{{\rm rd}}}
\newcommand{\Endmrp}{\underline{\End}^{{\rm mrp}}}

\newcommand{\Cob}{{\rm Cob}}

\newcommand{\Indrm}{{\rm Ind}_{{\rm rd}}^{{\rm mrp}}}

\DeclareMathAlphabet{\mathpzc}{OT1}{pzc}{m}{it}

\title{Ribbon dioperads and modular ribbon properads}
\author{Wai-Kit Yeung}
\address{Kavli IPMU, The University of Tokyo}
\email{wai-kit.yeung@ipmu.jp}
\begin{document}

\begin{abstract}
We define and study the notions of ribbon dioperads and modular ribbon properads. We give a Lie algebra structure on the colimit total object and the limit total object of a ribbon dioperad, and we give a norm map between them. We give a cobar construction for dg ribbon co-dioperads. We also prove similar results for modular ribbon properads. These results are applied to the study of higher Hochschild cochain complexes and pre-Calabi-Yau categories.
\end{abstract}

\maketitle


\tableofcontents

\section{Introduction}

Given an associative algebra $A$ over a field $\bK$, the $p$-th higher (cohomological) Hochschild complex is the complex
\begin{equation}  \label{higher_HH_intro_1}
	\RHom_{(A^{\otimes p})^{e}}(A^{\otimes p}, {}_{\sigma}(A^{\otimes p}))
\end{equation}
where $\Hom_{B^e}(-,-)$ means maps of $B$-bimodules, and the notation ${}_{\sigma}(A^{\otimes p})$ means that the left module structure is twisted by the cyclic rotation map $\sigma : A^{\otimes p} \xra{\cong} A^{\otimes p}$. 

Denote by $\cR(A) \xra{\sim} A$ the standard bar resolution $\cR(A) = \bigoplus_{n \geq 0} A \otimes A[1]^{\otimes n} \otimes A$, then by using the free bimodule resolution $\cR(A)^{\otimes p} \xra{\sim} A^{\otimes p}$, we see that \eqref{higher_HH_intro_1} can be represented by the explicit complex
 \begin{equation}  \label{higher_HH_intro_2}
 	C_H^{(p)}(A) \, = \, \prod_{n_1,\ldots,n_p \geq 0} \, \Hom_{\bK}(A[1]^{\otimes n_1} \otimes \ldots \otimes A[1]^{\otimes n_p}, A^{\otimes p})
 \end{equation}
This complex has a natural action by the cyclic group $C_p := \bZ/(p)$. One of the immediate applications of the theory we develop in this paper is the following result ({\it cf.} \cite{Yeu18, KTV}):
\bthm  \label{Lie_alg_intro}
Both the $C_p$-invariant $\bigoplus_{p \geq 1} C_H^{(p)}(A)^{C_p}[1]$ and the $C_p$-coinvariant $\bigoplus_{p \geq 1} C_H^{(p)}(A)_{C_p}[1]$ admits a weight graded dg Lie algebra structure whose Lie bracket has weight grading $-1$. (Here, we say that $C_H^{(p)}(A)$ has weight grading $p$).
Moreover, the norm map $\Nm : C_H^{(p)}(A)_{C_p} \ra C_H^{(p)}(A)^{C_p}$ obtained by the sum over cyclic rotations, is a map of dg Lie algebras.
\ethm

In fact, the same is true for $A_{\infty}$-algebras (or more generally, for $A_{\infty}$-categories). Moreover, if we forget about the differential, then the Lie bracket is of the form $[x,y] = x \circ y - (-1)^{|x||y|} y \circ x$, where $\circ$ is part of a ``graded pseudo-pre-Lie algebra'' (a generalization of pre-Lie algebras, see Defintion \ref{pseudo_pre_Lie_def} and Lemma \ref{pseudo_pre_Lie_implies_Lie}).
In characteristic zero, one can then define an $n$-pre-Calabi-Yau structure on an $A_{\infty}$-algebra to be a Maurer-Cartan element in a certain $(2-n)$-shifted version of $\prod_{p \geq 2} C_H^{(p)}(A)_{C_p}[1] \cong \prod_{p \geq 2} C_H^{(p)}(A)^{C_p}[1]$ (see Section \ref{PCY_sec} for details).

In this paper, we define and study a notion of ribbon dioperads, which is a natural conceptual framework to understand and prove Theorem \ref{Lie_alg_intro}. 
This will also allow us to replace the $C_p$-(co)invariants by its homotopy version. We expect this homotopy version to be useful in positive characteristics, as well as in applications to Fukaya categories (see Remark \ref{Fuk_remark} below).
Our conceptual framework also clarifies the origin of the shifts in \eqref{higher_HH_intro_2} as a Hadamard product, and hence gives a recipe to handle Koszul signs.
To define ribbon dioperads, we recall a reformulation of the notion of operads in \cite{Get09} and \cite{KW17}. 

Recall that a (non-unital) operad in a symmetric monoidal category $\cC$ consists of a sequence of objects $\cP(n) \in \cC$ with $S_n$-action, together with maps $\circ_i : \cP(m) \otimes \cP(n) \ra \cP(m+n-1)$ for $1 \leq i \leq m$, satisfying certain relations. One may regard $\cP(n)$ to be the object associated to a rooted corolla $C(n,1)$ of $n$ leaves (and $1$ root); the $S_n$-action to be induced by automorphisms of the rooted corolla $C(n,1)$; and the maps $\circ_i$ to be induced by the operation that grafts the root of $C(n,1)$ into the $i$-th leaf of $C(m,1)$ to become $C(m+n-1,1)$. In fact, the relations that the maps $\circ_i$ are required to satisfy also follows naturally from the combinatorics of the grafting operations of rooted corolla. 

These grafting operations of rooted corolla can be organized into a category, which then gives a succinct definition of a (non-unital) operad. Namely, let $\frG$ be the category whose objects are disjoint unions of rooted corolla, and whose morphisms are ``merging operations'', as in the following diagram (see Section \ref{ribbon_diop_sec} for precise definition):
\begin{equation*}
	\includegraphics[scale=0.2]{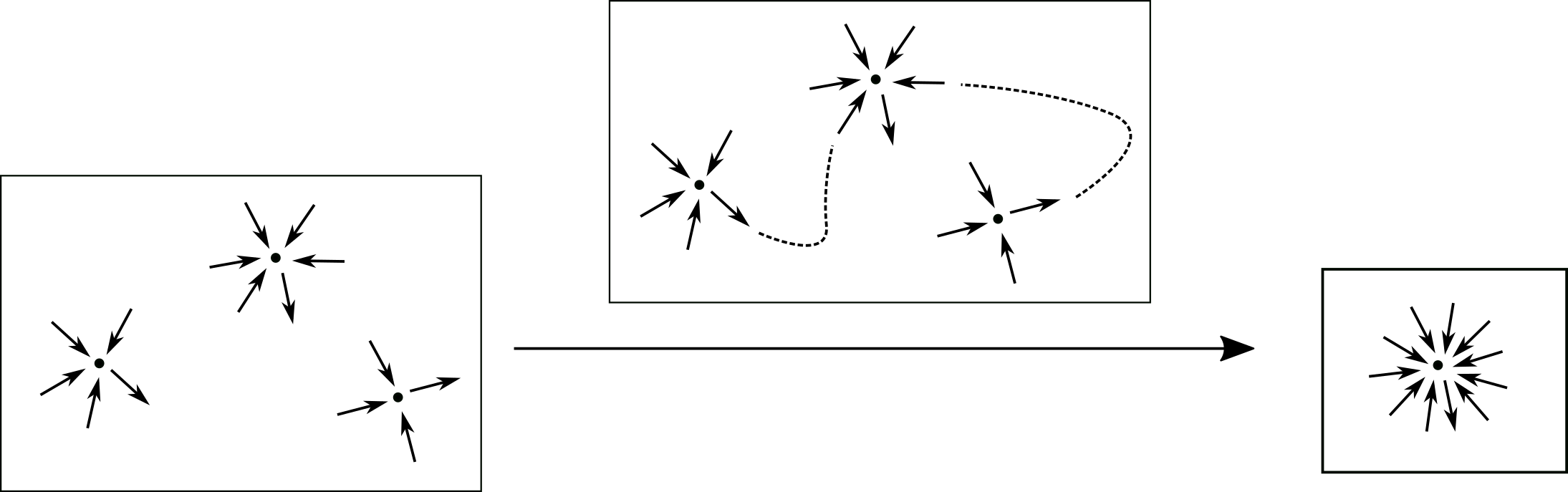}
\end{equation*}

This category $\frG$ has a symmetric monoidal structure given by taking disjoint union. Then a (non-unital) operad is simply a symmetric monoidal functor $\frG \ra \cC$. 

A large part of the theory of operads relies on the basic fact that free operads exist, and has a simple description.
In \cite{Get09} and \cite{KW17}, this basic fact is regarded as a consequence of a crucial property of $\frG$, namely that it is a regular pattern (we adopt the terminology of \cite{Get09} here%
\footnote{In \cite{KW17}, an essentially equivalent notion is called a Feynman category. While the main examples we consider indeed has a flavor of Feynman diagrams to it, the notion itself is of a general categorical nature, and is not necessarily related to Feynman diagrams. Hence, we prefer the broader name ``regular pattern'' (although perhaps a name more descriptive of the specific property may be more preferable).}).
One may regard the theory of regular patterns to be a further ``operization'' of the theory of operads, so that we have a sequence of ``de-operizations'' ({\it cf.} the introduction to \cite{GK94}):
\begin{equation}  \label{deoperize}
	\text{Modules} \quad \leftsquigarrow \quad \text{Algebras} \quad \leftsquigarrow \quad \text{Operads} \quad \leftsquigarrow \quad \text{Regular patterns}
\end{equation}

One advantage of this reformulation of the theory of operads is that we can then consider other regular pattern $\frG'$, and define a $\frG'$-operad to be a symmetric monoidal functor $\frG' \ra \cC$. Some examples are given in the following table:
\begin{center}
	\begin{tabular}{l | l | l}
		$\Ob(\frG)$ & $\Mor(\frG)$ & Symmetric monoidal $\frG \ra \cC$ \\
		\hline 
		Rooted corolla & Merging operation (automatically no cycles) & operads \\
		\hline
		Oriented corolla & Merging operations with no cycles & dioperads \\
		\hline
		Oriented corolla & Merging operations with no oriented cycles & properads
	\end{tabular}
\end{center}

In this paper, we consider the analogue of the three rows in this table, where the corolla are given a ribbon structure ({\it i.e.,} a cyclic ordering of the flags of the corolla). The first row then becomes the theory of non-symmetric operads, which is well studied. The second row becomes the theory of ribbon dioperads, which is the main focus of this paper. For the ribbon analogue of the third row, it is natural to also remember certain ``modular data'', and we consider accordingly a theory of modular ribbon properads (see Section \ref{mrp_sec}).

We now focus on the notion of ribbon dioperads. 
We mentioned above that the underlying $\Sigma$-module ({\it i.e.,} the sequence $\cP(n)$ with $S_n$-action) of an operad can be regarded as the assignments on rooted corolla and isomorphisms between them, and the composition morphisms as the assignments on merging operations.
Replacing the combinatorics of rooted corolla with that of ribbon oriented corolla (and merging operations with no cycles), one can characterize a ribbon dioperad in $\cC$ as the data that associates an object $\cP(C) \in \cC$ to each oriented corolla $C$, and an isomorphism $\cP(C) \cong \cP(C')$ to each isomorphism $C \cong C'$ (the part of this data is called the underlying $\frf$-module of $\cP$), together with composition maps $\cP(C_1) \otimes \cP(C_2) \ra \cP(C_0)$ for any $1$-edge contraction from $C_1 \amalg C_2$ to $C_0$.

We now focus on $\cC = \Ch(\bK)$. Recall that given any chain complex $A \in \cC$, one can form a dg operad $\Endcom^{{\rm operad}}(A)$ given by $\Endcom^{{\rm operad}}(A)(n) = \Homcom(A^{\otimes n},A)$. Similarly, in the theory of dioperads or properads, since one allow more than one output, there is also a corresponding dioperad/properad given by $\cP(n,m) = \Homcom(A^{\otimes n},A^{\otimes m})$. 
In the theory of ribbon dioperads, one also remembers the cyclic ordering of the flags of a corolla. Accordingly, the corresponding endomorphism ribbon dioperad $\cP = \Endrd(A)$ organizes the data \eqref{higher_HH_intro_2}. Then, taking the $C_p$-invariant and $C_p$-coinvariant as in Theorem \ref{Lie_alg_intro} corresponds respectively to taking limits and colimits of the underlying $\frf$-module of $\cP$. One can show that (see Section \ref{PCY_sec}) Theorem \ref{Lie_alg_intro} follows from a general statement that says that the total limit object and the total colimit object of a ribbon dioperad have canonical Lie algebra structures (Theorem \ref{Lie_colim_lim_diop}), and there is a norm map between them (Theorem \ref{norm_map_Lie}).

This hints at a Koszul duality theory for ribbon dioperads. We do not develop such a theory in full here. Instead, we will be content with only giving a cobar construction (see Section \ref{cobar_sec}). This will be applied in Section \ref{PCY_sec} to construct a dg ribbon dioperad $\cP_{{\rm PCY}}$ whose algebras are precisely pre-Calabi-Yau algebras. To each dg ribbon dioperad $\cP$, there is an associated PROP $\Cob( \Indrm(\cP)^+ )$ that controls the same class of algebras. 
We expect (a suitable modification of) the associated PROP of $\cP = \cP_{{\rm PCY}}$ to play the same role as the combinatorial model for the (colored) PROP in \cite{Cos07} that controls open topological conformal field theories. This question is also studied in \cite{KTV}.

\vspace{0.2cm}

\textbf{Acknowledgement.} 
The author would like to thank Sheel Ganatra, Ezra Getzler, and Mikhail Kapranov for helpful discussions.

\vspace{0.2cm}

\textbf{Conventions and notations.} 
We work with homological grading, so that $\deg(d) = -1$, and $V[1]_n = V_{n-1}$. The notation $\bK$ will always refer to a fixed commutative ring with $1$. Further assumptions on $\bK$ will be specified if needed.
The notation $\Ab, \Set, \Setdel, \Mod(\bK), \GrMod(\bK), \Ch(\bK)$ are respectively the category of abelian groups, sets, simplicial sets, $\bK$-modules, graded $\bK$-modules, and chain complexes of $\bK$-modules. All these are endowed with the usual symmetric monoidal structures, where the ones on $\GrMod(\bK), \Ch(\bK)$ have Koszul signs.

\section{Regular patterns}

Fix a complete and cocomplete closed symmetric monoidal category $(\cV,\boxtimes)$. We will mostly be interested in the cases $(\Set,\times)$, $(\Setdel,\times)$, $(\Ch(k),\otimes)$. 
For a $\cV$-category $\cC$, denote the underlying category by $|\cC|$. {\it i.e.,} we have $|\cC|(x,y) = \cV(\mathbf{1}_{\cV},\cC(x,y))$.
For $\cV$-categories $\cC,\cD$, denote by $\cC \boxtimes \cD$ the $\cV$-category with object set $\Ob(\cC) \times \Ob(\cD)$, and with Hom-objects given by $(\cC \boxtimes \cD)((c,d),(c',d')) = \cC(c,c') \boxtimes \cD(d,d')$.
A \emph{symmetric monoidal $\cV$-category} is a $\cV$-category $\cC$ equipped with a $\cV$-functor $\otimes : \cC \boxtimes \cC \ra \cC$, together with an extension of the underlying functor $\otimes: |\cC| \times |\cC| \ra |\cC|$ to a symmetric monoidal category structure. 
Given symmetric monoidal $\cV$-categories $\cC,\cD$, then a \emph{(strong/lax/oplax) symmetric monoidal $\cV$-functor} is a $\cV$-functor $F : \cC \ra \cD$ together with an extension of $F : |\cC| \ra |\cD|$ to a (strong/lax/oplax) symmetric monoidal functor. We will often omit the word ``strong''.
A \emph{(symmetric) monoidal natural transformation} is simply a monoidal natural transformation between the (strong/lax/oplax) symmetric monoidal functors on the underlying categories.

For example, for $\cV = \Ch(\bK)$, this gives a notion of symmetric monoidal dg category $\cC$. By our definition, the tensor bifunctor is a dg bifunctor, while the associators, unitors, and symmetrizers are required to be morphisms in $|\cC| = Z^0(\cC)$. Similarly, all the structure maps for symmetric monoidal functors from $\cC$ to $\cD$, and natural transformations between them, are required to be in $Z^0(\cD)$.

\blm  \label{strictify_lemma}
\begin{enumerate}
	\item If $F : \cC \ra \cD$ is a symmetric monoidal $\cV$-functor that is an equivalence of $\cV$-categories, then there exists a symmetric monoidal $\cV$-functor $G : \cC \ra \cD$, with monoidal natural isomorphisms $F \circ G \cong \id$ and $G \circ F \cong \id$.
	\item Every symmetric monoidal $\cV$-functor $F : \cC \ra \cD$ admits a canonical factorization $\cC \xra{F'} \cE' \xra{u} \cE'' \xra{i} \cD$, where $\Ob(\cC) = \Ob(\cE')$, $F'$ is identity on object set and is strictly symmetric monoidal, $u$ is a symmetric monoidal equivalence ({\it i.e.,} it satisfies the conditions of (1)), and $i : \cE'' \subset \cD$ is a full subcategory closed under $\otimes$. In particular, if $F$ is fully faithful ({\it i.e.,} the maps $F : \cC(c,c') \ra \cD(F(c),F(c'))$ in $\cV$ are isomorphisms), then $F'$ is an isomorphism of $\cV$-categories; if $F$ is essentially surjective ({\it i.e.,} the functor $F : |\cC| \ra |\cD|$ is essentially surjective), then $i$ is an equality of $\cV$-categories.
\end{enumerate}
\elm

\bpf
For (1), choose a quasi-inverse $\cV$-functor $G$. Since $|G|$ is a right adjoint of $|F|$, it has a canonical lax monoidal structure by doctrinal adjunction, which is clearly strong. One can verify that the natural isomorphisms $|F| \circ |G| \cong \id_{\cD}$ and $|G| \circ |F| \cong \id_{\cC}$ are monoidal.

For (2), let $\cE''$ be the essential image of $|F|$, which is clearly closed under $\otimes$. Let $\cE'$ be the $\cV$-category whose object set is $\Ob(\cC)$, and whose morphism objects are $\cE'(c,c') := \cD(F(c),F(c'))$. It is clear how to define the symmetric monoidal structure on $\cE'$, as well as the symmetric monoidal functors $F'$ and $u$. One then checks the compatibilities of these structures.
\epf

Given $\cV$-categories $\cA$ and $\cC$ where $\cA$ is small, denote by $[\cA,\cC]$ the $\cV$-category of $\cV$-functors between them (see, {\it e.g.}, \cite{Kel05}). If $\cA$ and $\cC$ are symmetric monoidal $\cV$-categories (with $\cA$ small), then denote by $[\![\cA,\cC]\!]$ the ordinary category of (strong) symmetric monoidal $\cV$-functor and monoidal natural transformations. 

For any small $\cV$-category $\cA$, denote by $\Mod(\cA) := [\cA^{\op},\cV]$, the $\cV$-category of right modules. One can tensor a left and a right module over $\cA$. Namely, for $M \in \Mod(\cA)$ and $N \in \Mod(\cA^{\op})$, define
$M \boxtimes_{\cA} N \in \Ob(\cV)$ as the $\cV$-enriched coend
\begin{equation*}
	M \boxtimes_{\cA} N \, := \, \int^{a \in \cA} M(a) \boxtimes N(a)
\end{equation*}
Given a $\cV$-functor $\tau : \cA \ra \cB$, then for any $M \in \Mod(\cA)$, one can use a similar coend to define $M \boxtimes_{\cA} \cB \in \Mod(\cB)$.
This fits into an extension-restriction adjunction
\begin{equation*}
	\begin{tikzcd}
		\tau_! = - \boxtimes_{\cA} \cB \, : \, \Mod(\cA) \ar[r, shift left] & \Mod(\cB) \, : \, \tau^* = (-)|_{\cA} \ar[l, shift left]
	\end{tikzcd}
\end{equation*}
For example, dg algebras are $\cV$-categories with $1$ object for $\cV = \Ch(\bK)$. For dg algebras and for maps between them, these give the usual notions of left/right dg modules, tensor products, extension/restrictions.

Now if $\cA$ is a small symmetric monoidal $\cV$-category, then there is an induced symmetric monoidal $\cV$-category structure on $\Mod(\cA)$, given by the Day convolution. Namely, for $M,N \in \Mod(\cA)$, it is clear how to define $M \boxtimes N \in \Mod(\cA \boxtimes \cA)$. The Day convolution $M * N \in \Mod(\cA)$ is then defined as the extension of this along the $\cV$-functor $\otimes : \cA \boxtimes \cA \ra \cA$. {\it i.e.,} we have
\begin{equation*}
	M * N \, := \, (M \boxtimes N)  \boxtimes_{\cA \boxtimes \cA} \cA
\end{equation*}
This gives a symmetric monoidal structure on $\Mod(\cA)$, which is cocontinuous in each variable, and the Yoneda functor $\cA \ra \Mod(\cA)$ is symmetric monoidal. 

If $\tau : \cA \ra \cB$ is a symmetric monoidal $\cV$-functor between small symmetric monoidal $\cV$-categories, then for any $M,N \in \Mod(\cB)$, there is a canonical map of $\cA$-modules
\begin{equation}  \label{Day_restr_lax}
	(\tau^*M) * (\tau^*N) = ((\tau \boxtimes \tau)^*(M \boxtimes N))  \boxtimes_{\cA \boxtimes \cA} \cA 
	\raq (M \boxtimes N) \boxtimes_{\cB \boxtimes \cB} \cB|_{\cA} = \tau^*(M * N)
\end{equation}
which gives a lax monoidal functor structure on $\tau^* : \Mod(\cB) \ra \Mod(\cA)$.

\bdf[\cite{Get09}]
A symmetric monoidal $\cV$-functor $\tau : \cA \ra \cB$ between small symmetric monoidal $\cV$-categories is said to be a \emph{pattern} if the map \eqref{Day_restr_lax} is an isomorphism for any $M, N \in \Mod(\cA)$.

By cocontinuity, it suffices to check it for representable modules $M,N$.
\edf

We unravel this definition for $\cV = \Set$. Namely, consider representable $M,N$, say $M = \cB(-,b_1)$ and $N = \cB(-,b_2)$ for $b_1,b_2 \in \Ob(\cB)$, then the map \eqref{Day_restr_lax} of functors $\cA^{\op} \ra \Set$ is given on $a \in \Ob(\cA)$ by
\begin{equation}  \label{Day_restr_lax_set}
	\{ a \xra{\cA} a_1 \otimes a_2  \xra{\cB \times \cB} b_1 \otimes b_2\}/_{\sim}
	\raq 
	\{ a \xra{\cB} b_1 \otimes b_2\}
\end{equation}
On the left hand side, we consider two objects $a_1, a_2 \in \Ob(\cA)$, a morphism $f : a \ra a_1 \otimes a_2$ in $\cA$, and morphisms $g_i : \tau(a_i) \ra b$ in $\cB$ for $i = 1,2$. 
We say that $(a_1,a_2,f,g_1,g_2)$ is equivalent to $(a'_1,a'_2,f',g'_1,g'_2)$ if there is are maps $\theta_i : a_i \ra a_i'$ in $\cA$ intertwining in the obvious sense with $f,f',g_i,g_i'$. Take the equivalence relation $\sim$ generated by this.
The left hand side of \eqref{Day_restr_lax_set} then says we quotient out by this equivalence relation.
The right hand side of \eqref{Day_restr_lax_set} is the set $\cB(\tau(a),b_1 \otimes b_2)$. 
The map \eqref{Day_restr_lax_set} is the obvious composition map.

A composition of a pattern is a pattern. This follows from the fact that the lax monoidal structure \eqref{Day_restr_lax} composes:
\blm  \label{Day_lax_compose}
Given symmetric monoidal $\cV$-functors $\cA_1 \xra{\tau_1} \cA_2 \xra{\tau_2} \cA_3$, denote by $\varphi_1$, $\varphi_2$ and $\varphi$ the map \eqref{Day_restr_lax} for $\tau_1$, $\tau_2$ and $\tau := \tau_2 \circ \tau_1$ respectively. Then for any $M,N \in \Mod(\cA_3)$, the composition
\begin{equation*}
	(\tau_1^* \tau_2^*M) * (\tau_1^* \tau_2^* N)
	\xraq{\varphi_1} \tau_1^*( (\tau_2^*M ) * ( \tau_2^* N ) )
	\xraq{\tau_1^*( \varphi_2 )} \tau_1^* \tau_2^* (M*N)
\end{equation*}
is equal to $\varphi$.
\elm

\bpf
This is clear from the tensor notation for the map \eqref{Day_restr_lax}.
\epf

Given a $\cV$-functor $\tau : \cA \ra \cB$ between small $\cV$-categories,
and given any $\cV$-cocomplete $\cV$-category $\cC$ there is a $\cV$-adjunction
\begin{equation}  \label{restr_ext_adj_2}
	\begin{tikzcd}
	\tau_! \, : \, {[}\cA,\cC{]} \ar[r, shift left] & {[}\cB,\cC{]} \, : \, \tau^* \ar[l, shift left]
	\end{tikzcd}
\end{equation}
where $\tau^*$ is the restriction functor, and $\tau_!$ is the left Kan extension. Explicitly, for $F \in [\cA,\cC]$, it is given by
\begin{equation}  \label{left_Kan_tensor}
	(\tau_! F)(b) \, = \, \tau^*(h_b) \boxtimes_{\cA} F
\end{equation}
Here, $\tau^*(h_b) \in [\cA^{\op},\cV]$ is the pullback along $\tau$ of the functor $h_b := \cB(-,b)$ representable by $b \in \Ob(\cB)$. Since $\cC$ is $\cV$-cocomplete, it is in particular $\cV$-copowered, so that we may tensor $\tau^*(h_b) \in [\cA^{\op},\cV]$ with $F \in [\cA,\cC]$ over $\cA$. The result is $\tau^*(h_b) \boxtimes_{\cA} F \in \cC$, which is the definition of the right hand side of \eqref{left_Kan_tensor}.

\bpp[\cite{Get09}]  \label{pattern_adj}
If $\tau : \cA \ra \cB$ is a pattern, and if $\cC$ is a $\cV$-cocomplete symmetric monoidal $\cV$-category such that $\otimes$ is cocontinuous in each variable, then the adjunction \eqref{restr_ext_adj_2} canonically induces an adjunction
	\begin{equation}  \label{restr_ext_adj_3}
		\begin{tikzcd}
			\tau_! \, : \, {[\![}\cA,\cC{]\!]} \ar[r, shift left] & {[\![}\cB,\cC{]\!]} \, : \, \tau^* \ar[l, shift left]
		\end{tikzcd}
	\end{equation}
\epp

\bpf
Clearly, if $G : \cB \ra \cC$ is symmetric monoidal, then $\tau^*(G)$ has a canonical symmetric monoidal structure. Thus, the non-trivial content of the statement is that the left Kan extension of a symmetric monoidal functors is a symmetric monoidal functor. Thus, given $F \in [\![\cA,\cC]\!]$, denote by $F^{(2)}$ the composition $\cA \boxtimes \cA \xra{F \boxtimes F} \cC \boxtimes \cC \xra{\otimes} \cC$. Then we use
\eqref{left_Kan_tensor} to compute
\begin{equation*}
	\begin{split}
		& (\tau_! F)(b_1) \otimes (\tau_! F)(b_2) \, = \, 
		( \tau^*(h_{b_1}) \boxtimes_{\cA} F ) \otimes ( \tau^*(h_{b_1}) \boxtimes_{\cA} F ) 
		\, \stackrel{(1)}{=} \, ( \tau^*(h_{b_1}) \boxtimes \tau^*(h_{b_1}) ) \boxtimes_{\cA \boxtimes \cA} F^{(2)}
		\\
		& \stackrel{(2)}{=} \, (\tau^*(h_{b_1}) * \tau^*(h_{b_2})) \boxtimes_{\cA} F \xraq{\cong} \tau^*(h_{b_1} * h_{b_2}) \boxtimes_{\cA} F = (\tau_! F)(b_1 \otimes b_2) 
	\end{split}
\end{equation*}
where (1) uses the fact that $\otimes$ is cocontinuous in each variable in $\cC$; while (2) uses the alternative expression for $F^{(2)}$ as the composition $\cA \boxtimes \cA \xra{\otimes} \cA \xra{F} \cC$.
To show that \eqref{restr_ext_adj_3} forms an adjunction, verify that the unit and counit maps of \eqref{restr_ext_adj_2} are monoidal natural transformations.
\epf

For any $\cV$-category $\cA$, denote by $\bS^n(\cA)$ the $\cV$-category whose objects are $n$-tuples $(a_1,\ldots,a_n)$ of objects in $\cA$, and whose morphisms are given by 
\begin{equation*}
	\bS^n(\cA)((a_1,\ldots,a_n),(a_1',\ldots,a_n')) \, = \, \coprod_{\sigma \in S_n} \, \cC(a_1, a'_{\sigma(1)}) \boxtimes \ldots \boxtimes  \cC(a_n, a'_{\sigma(n)})
\end{equation*}
(where coproduct is taken in $\cV$) with the obvious composition. 
Denote by $\bS(\cA)$ the disjoint union $\bS(\cA) := \coprod_{n \geq 0} \bS^n(\cA)$, which clearly has a canonical symmetric monoidal structure.

Choose, for each $n \geq 3$, a bracketing of the sequence of $n$ dots. For example, we may fix the choice $(\ldots ((\bullet \, \bullet) \bullet ) \ldots \bullet)$.
Then for any symmetric monoidal $\cV$-category $\cC$, and any $\cV$ functor $F : \cA \ra \cC$, there is a symmetric monoidal $\cV$-functor $F^{\otimes} : \bS(\cA) \ra \cC$ whose effect on objects is specified by the chosen bracketing. Indeed, the fact that the obvious specification on morphisms is well-defined follows from MacLane's coherence theorem.
If $\cA$ is small, this induces an equivalence
\begin{equation*}
	|[\cA ,\cC]| \, \simeq \, [ \! [  \bS(\cA) , \cC ] \! ] 
\end{equation*}
between the category $[ \! [  \bS(\cA) , \cC ] \! ] $ of symmetric monoidal $\cV$-functors and the underlying category $|[\cA ,\cC]|$ of the $\cV$-category $[\cA,\cC]$ of $\cV$-functors.

\bdf[\cite{Get09}]
A \emph{regular pattern} in $\cV$ consists of the data $\frf \rinto \frF \ra \frG$ where
\begin{enumerate}
	\item $\frF \ra \frG$ is a $\cV$-pattern that is essentially surjective ({\it i.e.,} $|\frF| \ra |\frG|$ is essentially surjective).
	\item $\frf \ra \frF$ is a $\cV$-functor (necessarily $\cV$-fully faithful) where the induced symmetric monoidal functor $\bS(\frf) \ra \frF$ is an equivalence.
\end{enumerate}

A \emph{map} between regular patterns is a diagram
\begin{equation}  \label{map_of_reg_patterns}
	\begin{tikzcd}[row sep = 10]
		\frf_1 \ar[r] \ar[d, "f"]& \frF_1 \ar[r] \ar[d, "F"] & \frG_1  \ar[d, "G"] \\
		 \frf_2 \ar[r] & \frF_2 \ar[r] & \frG_2
	\end{tikzcd}
\end{equation}
where $f$ is a $\cV$-functor, $F,G$ are symmetric monoidal $\cV$-functors. The first square is commutative up to isomorphism as $\cV$-functors, and the second square is commutative up to isomorphism as symmetric monoidal $\cV$-functors.

By a slight abuse of notation, we will often specify only $\frG$ when we refer to regular patterns and maps between them.
\edf

\brm  \label{strictify_reg_pattern_rem}
We will be interested in categories $[ \! [ \frG , \cC ] \! ]$ of symmetric monoidal $\cV$-functors out of $\frG$. By Lemma \ref{strictify_lemma}(1), the category $[ \! [ \frG , \cC ] \! ]$ remains equivalent if $\frG$ is replaced by $\frG'$ that is symmetric monoidal equivalent to it.
Thus, in view of Lemma \ref{strictify_lemma}(2), we may assume that $\Ob(\cF) = \Ob(\frG)$ and $\frF \ra \frG$ is the identity on object sets and is strictly symmetric monoidal. We say that $\frF \ra \frG$ is \emph{strictified} in this case. We may also assume that $\frF = \bS(\frf)$.
\erm

\bdf
Let $\cC$ be a symmetric monoidal $\cV$-category.
Given a regular $\cV$-pattern $\frf \rinto \frF \ra \frG$, define the (ordinary) category of $\frG$-operads (resp. $\frG$-cooperads) in $\cC$ to be
\begin{equation*}
	\begin{split}
	\OP(\frG,\cC) \, &:= \, [\![ \frG , \cC ]\!]  \\
	\COP(\frG,\cC) \, &:= \, [\![ \frG^{\op} , \cC ]\!]
	\end{split}
\end{equation*}
The restriction of a (co)operad to $\frf$ is called the \emph{underlying $\frf$-module} (which is a left $\frf$-module in $\cC$ for $\frG$-operads, and a right $\frf$-module in $\cC$ for $\frG$-cooperads).
\edf

By Proposition \ref{pattern_adj}, if $\cC$ is a $\cV$-cocomplete symmetric monoidal $\cV$-category such that $\otimes$ is cocontinuous in each variable, then there is an adjunction
\begin{equation}  \label{operad_free_forget_adj}
\frG_{\frf}(-) \, : \, \frf \text{-} \Mod_{\cC} \, \simeq \, [ \! [ \frF , \cC ] \! ] \qquad \substack{\longrightarrow\\  \longleftarrow}
		\qquad  [ \! [  \frG , \cC ] \! ] \, = \,  \OP(\frG,\cC) \, : \, (-)|_{\frf}
\end{equation}
where the left adjoint $\frG_{\frf}(-)$ is given by left Kan extension from $\frF$ to $\frG$. For an $\frf$-module $V$, the $\frG$-operad $\frG_{\frf}(V)$ is called the \emph{free $\frG$-operad} on $V$. 

\blm  \label{SF_is_pattern}
For any $\cV$-functor $\tau : \cA \ra \cB$, the induced symmetric monoidal $\cV$-functor $\bS(\tau) : \bS(\cA) \ra \bS(\cB)$ is a pattern.
\elm

\bpf
Given $a \in \bS(\cA)$ and $b_1,b_2 \in \bS(\cB)$, we verify that the $\cV$-enriched version of \eqref{Day_restr_lax_set} is an isomorphism.
Say $b_1 = (b_1^{(1)},\ldots,b_1^{(p)})$, $b_2 = (b_2^{(1)},\ldots,b_2^{(q)})$, and $a = (a^{(1)},\ldots,a^{(n)})$. Then the $\cV$-enriched version of both sides of \eqref{Day_restr_lax_set} are given by
\begin{equation*}
	\coprod_{\varphi} \, \cB(\tau(a^{(1)}) , b^{(\varphi(1))}) \otimes \ldots \otimes \cB(\tau(a^{(n)}) , b^{(\varphi(n))})
\end{equation*}
where the coproduct in $\cV$ is taken over the set of bijections $\varphi : \{1,\ldots,n\} \ra \{1,\ldots,p\} \amalg \{ 1,\ldots q\}$, and we write $b^{(\varphi(i))}$ for either $b_1^{(\varphi(i))}$ or $b_2^{(\varphi(i))}$ depending on whether $\varphi(i)$ is in $ \{1,\ldots,p\}$ or $\{ 1,\ldots q\}$.
\epf

\bpp
Given a map \eqref{map_of_reg_patterns} of regular patterns, then both $F$ and $G$ are patterns.
\epp

\bpf
We have seen from Lemma \ref{SF_is_pattern} that $F$ is a pattern. Since patterns compose (see Lemma \ref{Day_lax_compose}), we see that the composition $\frF_1 \xra{F} \frF_2 \ra \frG_2$ is a pattern. By commutativity, this means that the composition $\frF_1 \ra \frG_1 \xra{G} \frG_2$ is a pattern. Since $\frF_1 \ra \frG_1$ is essentially surjective, the pullback functor detects isomorphism. Thus, applying Lemma \ref{Day_lax_compose} to the composition $\frF_1 \ra \frG_1 \xra{G} \frG_2$, we see that $G$ is a pattern.
\epf

\bcor
Given a map \eqref{map_of_reg_patterns} of regular patterns, then for any $\cV$-cocomplete symmetric monoidal $\cV$-category $\cC$ such that $\otimes$ is cocontinuous in each variable, there is an adjunction
\begin{equation*}
	\begin{tikzcd}
		G_! \, : \, \OP(\frG_1,\cC) \ar[r, shift left]
		& \OP(\frG_2,\cC) \, : \, G^* \ar[l, shift left]
	\end{tikzcd}
\end{equation*}
where $G_!$ is given by left Kan extension, and $G^*$ is given by restriction.
\ecor

Define the category of lax $\frG$-operads in $\cC$ to be the category of lax symmetric monoidal $\cV$-functors from $\frG$ to $\cC$:
\begin{equation*}
		\OP^{\lax}(\frG,\cC) \, := \, [\![ \frG , \cC ]\!]^{\lax}
\end{equation*}

The following proposition shows that there is a $\frG$-operad $\cP^a$ associated to every lax $\frG$-operad $\cP$, whose underlying $\frf$-module is unchanged. In practice, this last property $\cP^a|_{\frf} \cong \cP|_{\frf}$ often makes it obvious how to define $\cP^a$ (for example, a $\frG$-operad can often be described as an $\frf$-module with extra structures, see Proposition \ref{frG_gen_rel} and Remark \ref{frG_gen_rel_rem} below, and it will be clear in those cases how to read off these extra structures on $\cP^a$ from the lax monoidal structure on $\cP$). 
\bpp  \label{lax_assoc_operad}
The inclusion $\OP^{\lax}(\frG,\cC)\rinto \OP(\frG,\cC)$ has a right adjoint $\cP \mapsto \cP^a$. Moreover, the adjunction counit $\cP^a \ra \cP$ restricts to an isomorphism $\cP^a|_{\frf} \xra{\cong} \cP|_{\frf}$ on $\frf$.
\epp

\bpf
Assume that $\Ob(\frF)=\Ob(\frG)$, and $\iota : \frF \ra \frG$ is the identity on object set and is strictly symmetric monoidal (see Remark \ref{strictify_reg_pattern_rem}). Let $\cP^a_0 : \frF \ra \cC$ be the symmetric monoidal $\cV$-functor that extends $\cP|_{\frf}$. The lax monoidal structure on $\cP$ then induces a natural transformation $\alpha : \cP^a_0 \ra \iota^* \cP$. It suffices to define $\cP^a$ on morphisms, because its effect on objects and its symmetric monoidal structure maps are already specified by $\cP^a_0$.

Consider first the unenriched case ({\it i.e.,} the case $\cV = \Set$) for simplicity of notation. Given $a,b \in \frF$, assume without loss of generality that $b = b_1 \otimes \ldots \otimes b_r$ for $b_i \in \frf$. Then by \eqref{Day_restr_lax_set}, any morphism from $a$ to $b$ can be written in the form
\begin{equation*}
	a \xraq{f \in \frF} a_1 \otimes \ldots \otimes a_r \xraq{(g_1,\ldots,g_r) \in \frG \times \ldots \times \frG} b_1 \otimes \ldots \otimes b_r
\end{equation*}
for some $a_i \in \frF$ (not necessarily in $\frf$).
The image of this morphism under $\cP^a$ is defined to be the composition
\begin{equation*}
	\cP_0^a(a) \xra{\cP_0^a(f)} \cP_0^a(a_1) \otimes \ldots \otimes \cP_0^a(a_r) \xra{\alpha} \cP( a_1) \otimes \ldots \otimes \cP(a_r) \xra{\cP(g_i)} \cP(b_1) \otimes \ldots \otimes \cP(b_r) = \cP^a_0(b)
\end{equation*}
It is clear that this descends through the equivalence relation in \eqref{Day_restr_lax_set}. One can also check that it is functorial, and that the symmetric monoidal structure maps specified by $\cP^a_0$ are compatible with the functor $\cP^a$, and hence endow $\cP^a$ with a structure of a symmetric monoidal functor. Also, the natural transformation $\alpha : \cP^a_0 \ra \iota^* \cP$ is also natural with respect to morphisms in $\frG$, so that we have $\alpha : \cP^a \ra \cP$, which is a lax monoidal natural transformation.

It is clear that $\cP^a$ is strong symmetric monoidal, and that $\alpha : \cP^a \ra \cP$ is an isomorphism if $\cP$ is strong symmetric monoidal. From this, we see that $\cP \mapsto \cP^a$ is the left adjoint to inclusion, and $\alpha$ is the adjunction counit.
This completes the proof for the unenriched case. The enriched case is completely parallel: one simply replace \eqref{Day_restr_lax_set} with its enriched version, in terms of a coend, and recast the argument in categorical terms.
\epf

Throughout the rest of this section, we consider the unenriched case $\cV = \Set$. Thus, $\frG$ is set-theoretical, but $\frG$-operads are allowed to take value in any symmetric monoidal category $\cC$.

\brm
To dispel a possible confusion about the role of enrichment, notice that enrichment can be imposed at each stage of ``de-operization'' of \eqref{deoperize}. 
More precisely, if $\cV_0$ is a symmetric monoidal category, and for $i = 0,1,2$, $\cV_{i+1}$ is a symmetric monoidal $\cV_i$-category, then (modulo technical details) one can expect to define a regular pattern $\frG$ in $\cV_0$, then a $\frG$-operad $\cP$ in $\cV_1$, then a $\cP$-algebra $A$ in $\cV_2$, then an $A$-module $M$ in $\cV_3$.
%
%
Thus, even though we mostly work with set-theoretic regular patterns, operads over them are often still enriched, and often have a homotopy theory.
\erm

If $\cP_1$, $\cP_2$ are $\frG$-operads in $\cC$, then the functor
\begin{equation}  \label{Hadamard_1}
	\cP_1 \otimes \cP_2 \, : \, \frG \xraq{\cP_1 \times \cP_2} \cC \times \cC \xraq{\otimes} \cC
\end{equation}
is clearly symmetric monoidal, and hence defines a $\frG$-operads in $\cC$, called their \emph{Hadamard product} $\cP_1 \otimes \cP_2$.

Assume that $\cC$ is closed, then one can also define a Hom object from a $\frG$-cooperad $\cQ$ to a $\frG$-operad $\cP$. Namely, first consider the functor
\begin{equation}  \label{Hadamard_2}
	[\cQ , \cP]^{\lax} \, : \, \frG \xraq{\cQ \times \cP} \cC^{\op} \times \cC \xraq{[-,-]} \cC
\end{equation}
which has an obvious lax monoidal structure given by the canonical maps
\begin{equation*}
	[\cQ(a_1), \cP(a_1)] \otimes [\cQ(a_2), \cP(a_2)] \raq [\cQ(a_1) \otimes \cQ(a_2), \cP(a_1) \otimes \cP(a_2)]
\end{equation*}
The \emph{Hadamard Hom} from $\cQ$ to $\cP$ is the $\frG$-operad $[\cQ , \cP]$ associated to the lax $\frG$-operad $[\cQ , \cP]^{\lax}$ as in Proposition \ref{lax_assoc_operad}.

\brm
The Hadamard product and the Hadamard Hom is defined for regular patterns in Cartesian monoidal category $(\cV,\times)$ ({\it i.e.,} the monoidal product is the Cartesian product), such as $(\Set,\times)$ or $(\Setdel,\times)$, but not for $(\Ab,\otimes)$. For example, if $\frG$ is $\Ab$-enriched, then the obvious guess of a Hadamard product would be the functor $(\cP_1 \otimes \cP_2)(a) := \cP_1(a) \otimes \cP_2(a)$. But this functor is non-linear on Hom-sets, so it is not an $\Ab$-enriched functor.

Conceptually, this is because the category $\cC \times \cC$ in \eqref{Hadamard_1} appears in a dual role. On the one hand, it is a product in the category of categories, so that the first functor is well-defined. On the other hand, it is also $\cC \boxtimes \cC$ for $(\cV,\boxtimes) = (\Set,\times)$, so that the second functor is well-defined.
\erm

Now we consider modules and derivations for $\frG$-operads. We will formulate the statements for set-theoretic $\frG$ ({\it i.e.,} for $\cV=\Set$). Since we will be working with $\frG$-operads over symmetric monoidal $\Ab$-categories $\cC$, it is in fact more natural to work with regular patterns over $\cV = \Ab$, which would then recover the statements for set-theoretic $\frG$ by applying them to the regular pattern $\bZ[\frG]$ over $\Ab$ obtained by applying $\bZ[-] : (\Set,\times)\ra(\Ab,\otimes)$ to the Hom sets. However, we stick with the set-theoretic case for simplicity of terminology.

\bdf  \label{primitive_def}
Let $\frf \subset \frF \ra \frG$ be a set-theoretic regular pattern. A morphism $\varphi \in \frG(a,b)$ is said to be \emph{primitive} if its target $b$ is in $\frf$. By \eqref{Day_restr_lax_set}, every morphism in $\frG$ is a composition of a morphism in $\frF$, followed by a $\otimes$-product of primitive morphisms. 
We will abbreviate this by saying that $\frG$ is $\otimes$-generated by primitive morphisms over $\frF$.
\edf

Let $\cC$ be a symmetric monoidal additive category ({\it i.e.,} symmetric monoidal $\cV$-category for $\cV = \Ab$) such that $\otimes$ commutes with finite coproduct in each variable.
Let $\Cinf$ be the category $\cC \times \cC$, with symmetric monoidal structure $(V_0,V_1) \otimes (V_0',V_1') := (V_0 \otimes V_0', (V_0 \otimes V_1') \oplus (V_1 \otimes V_0') )$. 
The functor $\pi_0 : \Cinf \ra \cC$ given by $(V_0,V_1) \mapsto V_0$ is clearly symmetric monoidal.
The functor $\Tot : \Cinf \ra \cC$ given by $(V_0 , V_1) \mapsto V_0 \oplus V_1$ has a lax symmetric monoidal structure given by the map $\Tot(V_0,V_1) \otimes \Tot(V_0',V_1') \ra \Tot( (V_0,V_1) \otimes (V_0',V_1'))$ that sends $V_1 \otimes V_1'$ to zero.
The natural transformation $\pi : \Tot \Rightarrow \pi_0$ given by $V_0 \oplus V_1 \ra V_0$ that sends $V_1$ to zero is clearly a monoidal natural transformation (between lax monoidal functors).

\bdf  \label{P_module_def}
Denote by $\pi_0 : \OP(\frG,\Cinf) \ra \OP(\frG,\cC)$ the post-composition with $\pi_0 : \Cinf \ra \cC$. Then for any $\cP \in \OP(\frG,\cC)$, the category of \emph{$\cP$-modules} is the category $\pi_0^{-1}(\cP)$. 
Thus, a module is a symmetric monoidal functor of the form $\frG \xra{(\cP,M_{\cP})} \cC \times \cC = \Cinf$. We denote by $M$ the underlying $\frf$-module $M = (M_{\cP})|_{\frf}$. 
\edf

A $\cP$-module may be described as the underlying $\frf$-module $M : \frf \ra \cC$ with extra structures. Assume that $\frF \ra \frG$ is strictified as in Remark \ref{strictify_reg_pattern_rem}. Given an $\frf$-module $M$, then the functor $(\cP|_{\frf},M) : \frf \ra \Cinf$ determines a symmetric monoidal functor $\frF \simeq \bS(\frf) \xra{(\cP,M_{\cP})} \cC \times \cC = \Cinf$. 
To extend this from $\frF$ to $\frG$, one has to specify the effect of $(\cP,M_{\cP})$ on morphisms $\varphi \in \frG(a,b)$. Since $\frG$ is $\otimes$-generated by primitive morphisms over $\frF$, it suffices to consider the case when $\varphi$ is primitive ({\it i.e.,} when $b \in \frf$). Write $a = a_1 \otimes \ldots \otimes a_r$ for $a_i \in \frf$. Thus, to each such $\varphi \in \frG(a_1 \otimes \ldots \otimes a_r, b)$, the functoriality of $M_{\cP} : \frG \ra \cC$ gives a map
\begin{equation}  \label{module_structure_map}
	\alpha_{M}(\varphi) \, : \, \bigoplus_{i=1}^r \cP(a_1) \otimes \ldots \otimes \cP(a_{i-1}) \otimes M(a_i) \otimes \cP(a_{i+1}) \otimes \ldots \otimes \cP(a_r) \raq M(b)
\end{equation}
satisfying suitable compatibility conditions with the original $\frF$-module structure, and with compositions in $\frG$ (which we will not explicate). 

In practice, $\frG$ is often specified by a set of primitive generating morphisms $\varphi$ and relations $\psi$ over $\frF \simeq \bS(\frf)$ (see Proposition \ref{frG_gen_rel} and Remark \ref{frG_gen_rel_rem} below). In this case, it suffices to specify \eqref{module_structure_map} for these generating morphisms $\varphi$, and verify the corresponding relations for $\psi$. 
For example, for the regular patterns that controls (colored) ribbon dioperads, one can use Proposition \ref{frG_gen_rel} to specify the precise structure of a $\cP$-module from Definition \ref{P_module_def}.

\bdf  \label{P_derivation_def}
Given a module $(\cP,M_{\cP}) : \frG \ra \Cinf$, denote by $\Tot^{\lax}(\cP,M_{\cP})$ the post-composition with the lax symmetric monoidal functor $\Tot : \Cinf \ra \cC$.
The monoidal natural transformation $\pi : \Tot \Rightarrow \pi_0$ then induces a map $\pi : \Tot^{\lax}(\cP,M_{\cP}) \ra \cP$ of lax $\frG$-operads.
A \emph{derivation} from $\cP$ to $M$ is a section of $\pi : \Tot^{\lax}(\cP,M_{\cP}) \ra \cP$ in the category of lax $\frG$-operads%
\footnote{If we denote by $\Tot(\cP,M_{\cP})$ the $\frG$-operad associated to the lax $\frG$-operad $\Tot^{\lax}(\cP,M_{\cP})$ as in Proposition \ref{lax_assoc_operad}, then a derivation from $\cP$ to $M$ may alternatively be defined as a section of $\pi : \Tot(\cP,M_{\cP}) \ra \cP$ in the category of $\frG$-operads.}.
Denote the set of derivations from $\cP$ to $M$ by $\Der(\cP,M)$.
\edf

Notice that the underlying $\frf$-module of $\Tot^{\lax}(\cP,M_{\cP})$ is given by $\cP|_{\frf} \oplus M$. 
Clearly, any derivation is determined by the underlying map of $\frf$-modules $D : \cP|_{\frf} \ra M$. This map is required to satisfy a certain derivation property with respect to the $\cP$-module structure of $M$. Namely, for each primitive $\varphi \in \frG(a_1 \otimes \ldots \otimes a_r, b)$, we require the following diagram to commute (see \eqref{module_structure_map}):
\begin{equation}  \label{derivation_comm_diag}
	\begin{tikzcd}
		\cP(a_1) \otimes \ldots  \otimes \cP(a_r) \ar[r, "\cP(\varphi)"] \ar[d, "\id \otimes D(a_i) \otimes \id"']
		& \cP(b) \ar[d, "D(b)"] \\
		\bigoplus_{i=1}^r \, \cP(a_1) \otimes \ldots \otimes \cP(a_{i-1}) \otimes M(a_i) \otimes \cP(a_{i+1}) \otimes \ldots \otimes \cP(a_r) \ar[r, "\alpha_M(\varphi)"]
		& M(b)
	\end{tikzcd}
\end{equation}
We will simply say that $D : \cP|_{\frf} \ra M$ is a derivation if it satisfies this property.
If we specify $\frG$ by generators $\varphi$ and relations $\psi$ over $\frF$ (see Proposition \ref{frG_gen_rel} and Remark \ref{frG_gen_rel_rem} below), then it suffices to verify the derivation property for each such generating primitive morphism $\varphi$.

\blm  \label{der_free_operad}
If $\cP = \frG_{\frf}(V)$ is a free $\frG$-operad on an $\frf$-module $V$, then the restriction along the canonical map $V \ra \cP|_{\frf}$ (the adjunction unit of \eqref{operad_free_forget_adj}) induces a bijection
\begin{equation*}
	\Der(\cP,M) \, \cong \, \Hom_{\frf\text{-}\Mod}(V,M)
\end{equation*}
\elm

\bpf
As we noted above (see footnotes in Definition \ref{P_derivation_def}), a derivation can be defined as a section of $\pi : \Tot(\cP,M_{\cP}) \ra \cP$ in the category of $\frG$-operads. The statement then follows from freeness of $\cP$.
\epf

Our next goal is to show that a dg $\frG$-operad is the same as a graded $\frG$-operad together with a square-zero derivation of degree $-1$ (we work with homological grading). We will formulate this in a more general categorical setting, so that Koszul signs are handled once and for all.

Let $\cC$ be as above ({\it i.e.,}  $\cC$ is a symmetric monoidal additive category such that $\otimes$ commutes with finite coproduct in each variable). Let $\Sigma : \cC \ra \cC$ be an additive endofunctor, together with isomorphisms
\begin{equation}  \label{Sigma_phi_maps}
	(\Sigma V) \otimes W  \xlaq[\cong]{\phi^L_{V,W}} \Sigma (V \otimes W) \xraq[\cong]{\phi^R_{V,W}}  V \otimes (\Sigma W) 
\end{equation}
natural in $V,W \in \cC$.

For each $V_1,\ldots,V_n \in \cC$, there are many different ways to use these maps $\phi^L$, $\phi^R$ to obtain an isomorphism $\Sigma(V_1 \otimes \ldots \otimes V_n) \cong V_1 \otimes \ldots \otimes V_{i-1} \otimes \Sigma V_i \otimes V_{i+1} \otimes \ldots \otimes V_n$, one for each bracketing of $n$ dots. We require that all these maps are the same. It is clear that one suffices to impose this relation for $n=3$. {\it i.e.,} we require 
\begin{equation}  \label{Sigma_assoc}
	\begin{split}
	\phi^L_{V_1, V_2 \otimes V_3} = (\phi^L_{V_1, V_2} \otimes \id_{V_3}) \circ \phi^L_{V_1 \otimes V_2, V_3}  \, &: \, \Sigma(V_1 \otimes V_2 \otimes V_3) \raq \Sigma(V_1) \otimes V_2 \otimes V_3 \\
	(\id_{V_1} \otimes \phi^L_{V_2, V_3}) \circ \phi^R_{V_1, V_2 \otimes V_3} = (\phi^R_{V_1, V_2} \otimes \id_{V_3}) \circ \phi^L_{V_1 \otimes V_2, V_3} \, &: \, \Sigma(V_1 \otimes V_2 \otimes V_3) \raq V_1 \otimes \Sigma(V_2) \otimes V_3 \\
	(\id_{V_1} \otimes \phi^R_{V_2, V_3}) \circ \phi^R_{V_1, V_2 \otimes V_3} = \phi^R_{V_1 \otimes V_2, V_3} \, &: \, \Sigma(V_1 \otimes V_2 \otimes V_3) \raq V_1 \otimes V_2 \otimes \Sigma(V_3) 
	\end{split}
\end{equation}
Since all these isomorphisms are the same, we will often denote them simply by $\phi$, without fear of confusion.
We will also require the following compatibility with the unitality structure:
\begin{equation}  \label{Sigma_unital}
	\begin{split}
\Sigma(u^L_V) = u^L_{\Sigma(V)} \circ \phi^L_{V,\mathbf{1}} \, &: \, \Sigma(V \otimes \mathbf{1}) \raq \Sigma(V) \\
\Sigma(u^R_V) = u^R_{\Sigma(V)} \circ \phi^R_{\mathbf{1},V} \, &: \, \Sigma( \mathbf{1} \otimes V) \raq \Sigma(V)
\end{split}
\end{equation}
where $u^L_{W} : W \otimes \mathbf{1} \xra{\cong} W$ and $u^R_{W} : \mathbf{1} \otimes W \xra{\cong} W$ are the structure isomorphisms for the unit objects.
We also require that the following diagram commutes:
\begin{equation}  \label{Sigma_comm}
	\begin{tikzcd}[column sep = 40]
		\Sigma(V \otimes W)  \ar[r, "\Sigma(\sigma_{V,W})"] \ar[d, "\phi^L_{V,W}"'] &  \Sigma(W \otimes V) \ar[d, "\phi^R_{W,V}"] \\
		\Sigma (V) \otimes W \ar[r, "\sigma_{\Sigma(V),W}" ] & W \otimes \Sigma(V)
	\end{tikzcd}
\end{equation}

\bdf  \label{inf_sym_mon_def}
Let $\cC$ be a symmetric monoidal additive category such that $\otimes$ commutes with finite coproduct in each variable, then an \emph{infinitesimally symmetric monoidal endofunctor} is an additive endofunctor $\Sigma : \cC \ra \cC$ together with isomorphisms \eqref{Sigma_phi_maps} natural in $V,W$, satisfying \eqref{Sigma_assoc}, \eqref{Sigma_unital} and \eqref{Sigma_comm}.
\edf

\beg  \label{graded_module_shift_example}
Take $\cC = \GrMod(k)$, the category of graded modules, where the symmetry operations have Koszul signs. Take $\Sigma(V) = V[1]$, and define the maps $\phi^L$, $\phi^R$ by
\begin{equation*}
	\phi^L(s(x \otimes y)) = (sx) \otimes y
	\qquad \text{and} \qquad 
	\phi^R(s(x \otimes y)) = (-1)^{|x|} x \otimes (sy)
\end{equation*}
One can verify easily that all the above conditions \eqref{Sigma_assoc}, \eqref{Sigma_unital} and \eqref{Sigma_comm} are satisfied. It also satisfies \eqref{Sigma_odd} below.
\eeg

Let $\PCOM(\cC,\Sigma)$ be the category whose objects consist of a pair $(V,d)$, where $V \in \cC$ and $d : V \ra \Sigma(V)$. Morphisms in $\PCOM(\cC,\Sigma)$ are morphisms $V \ra W$ in $\cC$ that intertwine with $d$. This category is clearly additive, and has a symmetric monoidal structure given by $(V,d_V) \otimes (W,d_W) := (V \otimes W, d_{V \otimes W})$, where 
\begin{equation}  \label{d_on_tensor}
	d_{V \otimes W} := (\phi^L_{V,W})^{-1} \circ (d_V \otimes \id_W) + (\phi^R_{V,W})^{-1} \circ (\id_V \otimes d_W) \, : \, V \otimes W \raq \Sigma(V \otimes W)
\end{equation}
Indeed, the fact that the associativity, unitality, and commutativity structure maps on $(\cC,\otimes)$ are morphisms in $\PCOM(\cC,\Sigma)$ follow from the requirements \eqref{Sigma_assoc}, \eqref{Sigma_unital} and \eqref{Sigma_comm} respectively.

Consider the functor $\Delta_{\Sigma} : \cC \ra \Cinf$ given by $V \mapsto (V_0,V_1) := (V,\Sigma(V))$. It has a lax symmetric monoidal structure given by $\Delta_{\Sigma}(V) \otimes \Delta_{\Sigma}(W) \ra \Delta_{\Sigma}(V \otimes W)$ that sends $V \otimes \Sigma(W)$ and $\Sigma(V) \otimes W$ to $\Sigma(V\otimes W)$ via the maps $(\phi^R)^{-1}$ and $(\phi^L)^{-1}$ respectively. 
Indeed, the requirements \eqref{Sigma_assoc}, \eqref{Sigma_unital} and \eqref{Sigma_comm} are precisely the condition that this is a lax symmetric monoidal functor.

For any $\frG$-operad $\cP : \frG \ra \cC$, its post-composition $\Delta_{\Sigma} \circ \cP$ is a lax $\frG$-operad $\frG \ra \Cinf$, whose associated operad $(\cP,(\Sigma \cP)_{\cP}) : \frG \ra \Cinf$ is therefore a $\cP$-module, known as the \emph{$\Sigma$-shifted diagonal module}. Notice that its underlying $\frf$-module is simply $\Sigma(\cP|_{\frf})$.

Explicitly, the functor $(\Sigma \cP)_{\cP} : \frG \ra \cC$ is given by
\begin{equation*}
	(\Sigma \cP)_{\cP}(a_1 \otimes \ldots \otimes a_r) \, = \, \bigoplus_{i=1}^r \, \cP(a_1) \otimes \ldots \otimes \cP(a_{i-1}) \otimes \Sigma(\cP(a_i)) \otimes \cP(a_{i+1}) \otimes \ldots \otimes \cP(a_r)
\end{equation*}
for $a_1,\ldots, a_r \in \frf$. 
For a primitive morphism $\varphi \in \frG(a_1 \otimes \ldots \otimes a_r,b)$ where $b \in \frf$, the $\cP$-module structure map \eqref{module_structure_map}
is given by
\begin{equation}  \label{module_map_Sigma_diag}
	 \alpha_i =  \Sigma(\cP(\varphi)) \circ \phi^{-1} \, : \,  \cP(a_1) \otimes \ldots \otimes \cP(a_{i-1}) \otimes \Sigma(\cP(a_i)) \otimes \cP(a_{i+1}) \otimes \ldots \otimes \cP(a_r) \raq \Sigma(\cP(b))
\end{equation}

\bpp   \label{oeprad_in_precomp_der}
Giving a $\frG$-operad in $\PCOM(\cC,\Sigma)$ is equivalent to giving a $\frG$-operad $\cP : \frG \ra \cC$ in $\cC$, together with a derivation $D : \cP|_{\frf} \ra \Sigma(\cP|_{\frf})$ from $\cP$ to the $\Sigma$-shifted diagonal module.
\epp

\bpf
Assume that $\frF \ra \frG$ is strictified as in Remark \ref{strictify_reg_pattern_rem}.
Any map $D : \cP|_{\frf} \ra \Sigma(\cP|_{\frf})$ of $\frf$-modules determines an $\frf$-module $(\cP|_{\frf},D) : \frf \ra \PCOM(\cC,\Sigma)$, which then determines a symmetric monoidal functor $(\cP|_{\frf},D)^{\otimes} : \frF \ra \PCOM(\cC,\Sigma)$. If we forget about the differential, then $(\cP|_{\frf},D)^{\otimes}$, as a functor to $\cC$, is given by $(\cP|_{\frf})^{\otimes} \cong \cP|_{\frF}$.
In other words, we have specified a lifting of $\cP|_{\frF} : \frF \ra \cC$ to $(\cP|_{\frF},\widetilde{D}): \frF \ra \PCOM(\cC,\Sigma)$.
Since $\cP_{\frF}$ has a given extension to $\cP : \frG \ra \cC$, we would obtain an extension $(\cP,\widetilde{D}): \frG \ra \PCOM(\cC,\Sigma)$ if for all $\varphi \in \frG(a,b)$, the induced map $\cP(\varphi) : \cP(a) \ra \cP(b)$ intertwine with the induced differential $\widetilde{D}$. 
Again, it suffices to check it for primitive morphisms $\varphi$, which asserts the commutativity of
\begin{equation}  \label{derivation_D_tilde}
	\begin{tikzcd}
		\cP(a_1) \otimes \ldots  \otimes \cP(a_r) \ar[r, "\cP(\varphi)"] \ar[d, "\widetilde{D}(a_1{,}\ldots{,}a_r)"']
		& \cP(b) \ar[d, "D(b)"] \\
		\Sigma(\cP(a_1) \otimes \ldots  \otimes \cP(a_r)) \ar[r, "\Sigma(\cP(\varphi))"]
		& \Sigma(\cP(b))
	\end{tikzcd}
\end{equation}
where $\widetilde{D}(a_1,\ldots,a_r) := \sum_{i=1}^r \id \otimes D(a_i) \otimes \id$. 
The condition is precisely the commutativity of \eqref{derivation_comm_diag} for $M = \Sigma(\cP|_{\frf})$, where the bottom horizontal map is given in \eqref{module_map_Sigma_diag}. In other words,  the induced map $\cP(\varphi) : \cP(a) \ra \cP(b)$ intertwine with the induced differential $\widetilde{D}$ if and only if $D$ is a derivation.
\epf

Let $\COM(\cC,\Sigma) \subset \PCOM(\cC,\Sigma)$ be the full subcategory consisting of objects $(V,d)$ such that
\begin{equation*}
	0 = \Sigma(d) \circ d  \, : \, V \raq \Sigma(\Sigma V) 
\end{equation*}
In order for $\COM(\cC,\Sigma) \subset \PCOM(\cC,\Sigma)$ to be closed under $\otimes$, we will impose one more condition:
\begin{equation}  \label{Sigma_odd}
	\phi^L_{V,\Sigma W} \circ \Sigma(\phi^R_{V,W}) = - \phi^R_{\Sigma V,W} \circ \Sigma(\phi^L_{V,W})  \, : \, \Sigma(\Sigma(V \otimes W)) \raq \Sigma(V) \otimes \Sigma(W)
\end{equation}
If this holds for all $V,W$, we say that $\Sigma$ is \emph{odd}. If the analogue of \eqref{Sigma_odd} without the minus sign holds for all $V,W$, we say that $\Sigma$ is \emph{even}.
Clearly, if $\Sigma$ is odd, then $\COM(\cC,\Sigma) \subset \PCOM(\cC,\Sigma)$ is closed under $\otimes$.

\beg \label{chain_complex_example}
For $(\cC,\Sigma)$ in Example \ref{graded_module_shift_example}, the category $\COM(\cC,\Sigma)$ is precisely the category of chain complexes.
To reconcile with the signs of differentials on tensor product, for a map $d : V \ra \Sigma(V)$ in $\cC$, write $d = sd'$, where $d' : V \ra V$ is a map of degree $-1$. Then \eqref{d_on_tensor} reads $d'(x \otimes y) = d'(x) \otimes y + (-1)^{|x|} x \otimes d'(y)$, which is of course the usual sign convention.
\eeg

\beg \label{chain_complex_example_2}
Given $(\cC,\Sigma)$ with $\Sigma$ odd, then $\Sigma$ induces an endofunctor, still denoted as $\Sigma$, on $\COM(\cC,\Sigma)$. Namely, $\Sigma(V,d) := (\Sigma(V), -\Sigma(d))$. Notice that this last minus sign ensures that, for $V,W \in \COM(\cC,\Sigma)$, the morphisms \eqref{Sigma_phi_maps} are morphisms in $\COM(\cC,\Sigma)$. Hence, $\Sigma$ is still an odd infinitesimally symmetric monoidal endofunctor on $\COM(\cC,\Sigma)$.
Objects in $\COM(\COM(\cC,\Sigma),\Sigma)$ then consists of an object $V \in \cC$, together with two square zero maps $d_1, d_2 : V \ra \Sigma(V)$ that anticommute with each other: $\Sigma(d_1) \circ d_2 + \Sigma(d_2) \circ d_1 = 0$.
\eeg

\brm  \label{coherence_for_odd}
Recall that the axiom \eqref{Sigma_assoc} implies that there is a well-defined isomorphism
$\Sigma(V_1 \otimes \ldots \otimes V_n) \cong V_1 \otimes \ldots \otimes V_{i-1} \otimes \Sigma V_i \otimes V_{i+1} \otimes \ldots \otimes V_n$. 
A repeated application of this shows that for each $\sigma \in S_n$, there is a isomorphism
\begin{equation*}
	\psi_{\sigma} \,  : \, \Sigma^n(V_1 \otimes \ldots \otimes V_n) \raq \Sigma(V_1) \otimes \ldots \otimes \Sigma(V_n)
\end{equation*}
One can show by induction that if $\Sigma$ is odd, then we have $(-1)^{{\rm sgn}(\sigma)} \psi_{\sigma} = (-1)^{{\rm sgn}(\sigma')} \psi_{\sigma'}$.
\erm

\brm
The formalization of an infinitesimally symmetric monoidal endofunctor is in fact rather redundant. If we fix any $V_0 \in \cC$, then (we believe that) the endofunctor $\Sigma_{V_0}$ defined by $\Sigma_{V_0}(V) := V_0 \otimes V$ is infinitesimally symmetric monoidal. Conversely, given $\Sigma$, then (we believe that) one can use $\phi^L$ to identify $\Sigma \cong \Sigma_{V_0}$ for $V_0 = \Sigma(\mathbf{1})$. The even/odd condition then translates to $\sigma_{\Sigma(\mathbf{1}),\Sigma(\mathbf{1})} = \pm \id$.
However, we still find it benefitial to formulate the notion of an (odd) infinitesimally symmetric monoidal endofunctor because the axioms are exactly what one needs in practice, and also because we find the signs less confusing in this formulation.
\erm

Notice that if $\Sigma : \cC \ra \cC$ is an infinitesimally symmetric monoidal endofunctor, then so is $\Sigma^2 : \cC \ra \cC$, with the obvious induced structure maps $\phi^L$, $\phi^R$ (the verification of \eqref{Sigma_assoc} for $\Sigma^2$ is slightly more tedious, but in any case follows from applying the axiom \eqref{Sigma_assoc} for $\Sigma$ two times, and using the naturality of $\phi^L$ and $\phi^R$). 
Thus, one can define the $\Sigma^2$-shifted diagonal module $\Sigma^2(\cP|_{\frf})$.

\blm  \label{square_zero_deriv_lemma}
If $\Sigma$ is odd, and if $D : \cP|_{\frf} \ra \Sigma(\cP|_{\frf})$ is a derivation to the $\Sigma$-shifted diagonal module, then the composition
$\cP|_{\frf} \xra{D} \Sigma(\cP|_{\frf}) \xra{\Sigma(D)} \Sigma^2(\cP|_{\frf})$ is a derivation to the $\Sigma^2$-shifted diagonal module.
Moreover, this composition $\Sigma(D) \circ D$ is zero if and only if the functor $(\cP,\widetilde{D}) : \frG \ra \PCOM(\cC,\Sigma)$ induced by $D$ as in Proposition \ref{oeprad_in_precomp_der} lands in $\COM(\cC,\Sigma) \subset \PCOM(\cC,\Sigma)$.
\elm

\bpf
Recall that the derivation property is equivalent to the commutativity of \eqref{derivation_D_tilde}. Apply it two times, and we have the commutativity of 
\begin{equation}  \label{D_square_der_diag}
	\begin{tikzcd}
		\cP(a_1) \otimes \ldots  \otimes \cP(a_r) \ar[r, "\cP(\varphi)"] \ar[d, "\widetilde{D}"']
		& \cP(b) \ar[d, "D(b)"] \\
		\Sigma(\cP(a_1) \otimes \ldots  \otimes \cP(a_r)) \ar[r, "\Sigma(\cP(\varphi))"] \ar[d, "\Sigma(\widetilde{D})"']
		& \Sigma(\cP(b))  \ar[d, "\Sigma(D(b))"]\\
		\Sigma^2(\cP(a_1) \otimes \ldots  \otimes \cP(a_r)) \ar[r, "\Sigma^2(\cP(\varphi))"]
		& \Sigma^2(\cP(b))
	\end{tikzcd}
\end{equation}
The composition $\Sigma(\widetilde{D}) \circ \widetilde{D}$ of the two vertical maps on the left can be written as a sum $\Sigma(\widetilde{D}) \circ \widetilde{D} = \sum_{1 \leq i,j \leq r} \widetilde{D}^{(2)}_{i,j}$ in the obvious way. By assumption, $\Sigma$ is odd, so that $\widetilde{D}^{(2)}_{i,j} + \widetilde{D}^{(2)}_{j,i} = 0$ if $i \neq j$. Thus, we are left with a sum over $i=j$, so that the commutativity of \eqref{D_square_der_diag} asserts that $\Sigma(D) \circ D$ is a derivation.
This computation of $\Sigma(\widetilde{D}) \circ \widetilde{D}$ also shows the second statement.
\epf

\section{Ribbon dioperads}  \label{ribbon_diop_sec}

We focus on the unenriched case $\cV = \Set$ from now on. We first define several regular patterns in terms of the combinatorics of graphs, following \cite{Get09} and \cite{KW17}.

\bdf  \label{graph_def_1}
A \emph{graph} consists of the data $\Gamma = (V,F,p,\sigma)$ where $V,F$ are finite sets, $p : F \ra V$ is a map of sets, and $\sigma : F \ra F$ is an involution.
We call $V$ the set of vertices; $F$ the set of flags; and $p^{-1}(v)$ the set of flags meeting at the vertex $v \in V$. A fixed point of $\sigma$ is called a leg, the set of which is denoted by $L(\Gamma)$. A $2$-element orbit of $\sigma$ is called an edge, the set of which is denoted by $E(\Gamma)$.

A \emph{multi-corolla} is a triple $(V,F,p)$ where $V,F$ are finite sets, $p : F \ra V$ is a map of sets. A \emph{corolla} is a multi-corolla where $V$ is a singleton.
A (multi-)corolla is often regarded as a graph by setting $\sigma = \id$.

An \emph{orientation} on a graph $\Gamma$ is a map of sets $\theta : F \ra \{+,-\}$ that is injective on each $\sigma$-orbit. An \emph{orientation} on a (multi-)corolla is a map of sets $\theta : F \ra \{+,-\}$.

A \emph{ribbon structure} on a graph or (multi-)corolla is a permutation $\mu : F \ra F$ such that for each $v \in V$, the subset $p^{-1}(v)$ is either empty or consists of a single $\mu$-orbit.

Given a ribbon graph $\Gamma = (V,F,p,\sigma,\mu)$, take the disjoint union $F' \amalg F''$ of two copies of $F$, and consider the equivalence relation $\approx$ generated by the relations $f' \approx \mu(f)''$ for all $f \in F$, and $f' \approx \sigma(f)''$ for all $f \in F$ such that $\sigma(f) \neq f$. Then the  \emph{set of boundary components} of $\Gamma$, denoted as $\pi_0(\partial \Sigma(\Gamma) )$, is the disjoint union of $(F' \amalg F'')/_{\approx}$ and $(V \setminus p(F))$.
Given a set $\scO$, then an $\scO$-coloring of $\Gamma$ is a map of sets  $\pi_0(\partial \Sigma(\Gamma) ) \ra \scO$. 
Consider the equivalence relation $\sim$ generated by $\approx$ and $f' \sim f''$ for all leg $f \in L(\Gamma)$. Then the \emph{set of compactified boundary components} of $\Gamma$, denoted as $\pi_0(\partial \overline{\Sigma}(\Gamma) )$, is the disjoint union of $(F' \amalg F'')/_{\sim}$ and $(V \setminus p(F))$. Notice that $(F' \amalg F'')/_{\sim}$ can also be identified with the set of orbits of $F$ under $\sigma \circ \mu$. Namely, the identification $F = F'$ induces a bijection $F/(\sigma \circ \mu) \xra{\cong} (F' \amalg F'')/_{\sim}$.
\edf

One can define the $1$-dimensional CW complex associated to a graph, as well as the standard notions of connected components, (oriented) cycles, trees, forests, {\it etc}. Details will be omitted. 
A \emph{rooted tree} is an oriented graph that is a tree, such that there is a unique leg $f$ with $\theta(f) = -$. A \emph{rooted forest} is an oriented graph each of its component is a rooted tree. One can also define \emph{rooted (multi-)corolla} in the same way.

A ribbon structure on a graph allows one to thicken the associated $1$-dimensional CW complex to form an oriented surface $\Sigma(\Gamma)$ with boundary, where the original CW complex embeds as a deformation retract. The set of boundary components $\pi_0(\partial \Sigma(\Gamma))$, as defined combinatorially above, is canonically identified with the set of boundary components of $\Sigma(\Gamma)$. The surface $\Sigma(\Gamma)$ can be compactified by adding one point for each leg of $\Gamma$. We denote the compactified surface with boundary by $\overline{\Sigma}(\Gamma)$. Again, the set of compactified boundary components $\pi_0(\partial \overline{\Sigma}(\Gamma) )$, as defined combinatorially above, is canonically identified with the set of boundary components of $\overline{\Sigma}(\Gamma)$. Our conventions is fixed by the following example:

\[
\includegraphics[scale=0.5]{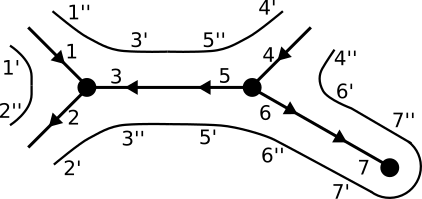}
\]
The graph here is defined by $V = \{a,b,c\}$, $F = \{1,2,3,4,5,6,7\}$, $p^{-1}(a,b,c) = (\{1,2,3\} ,\{4,5,6\} ,\{7\})$, and orbits of $\sigma$ are $\{3,5\}, \{6,7\}, \{1\}, \{2\}, \{4\}$. The orientation is specified by $\theta(\{1,3,4,7\}) = \{+\}$ and $\theta(\{2,5,6\}) = \{-\}$. The ribbon structure $\mu$ is defined by the counterclockwise order at each vertex, {\it i.e.}, $1 \mapsto 2 \mapsto 3 \mapsto 1$ and $5 \mapsto 6 \mapsto 4 \mapsto 5$ and $7 \mapsto 7$.
To identify the combinatorially defined set $\pi_0(\partial \Sigma(\Gamma))$ with the set of connected components of the boundary of the surface $\Sigma(\Gamma)$, for each given flag $f \in F$, stand at the vertex of that flag, facing towards the direction of the flag. Hold both hands forward and slightly open, then the left hand points to the connected component containing $f'$, while the right hand points to the connected component containing $f''$. The meaning of the equivalence relation $\approx$ is then clear.

\bdf
Given graphs $\Gamma$ and $\Gamma'$, then an \emph{edge contraction map} from $\Gamma$ to $\Gamma'$ consists of an injective map $\alpha : F' \rinto F$ intertwining $\sigma$ and $\sigma'$, which restricts to a bijection $\alpha : L(F') \xra{\cong} L(F)$, and a surjective map $\beta : V \ronto V'$ respecting $p$ and $p'$ in the sense that $\beta(p(\alpha(f'))) = p'(f')$ for all $f' \in F'$.

Given an edge contraction map, let $F_0 = F \setminus \alpha(F')$, and let $V_0 = p(F_0)$, then $\Gamma_0 = (V_0,F_0,p,\sigma)$ is called the \emph{contraction subgraph}.
\edf

\blm  \label{contr_ori_ribbon}
Given an edge contraction map from $\Gamma$ to $\Gamma'$.
\begin{enumerate}
	\item If $\Gamma$ has an orientation $\theta$, then $\Gamma'$ has an induced orientation $\theta' = \theta \circ \alpha$.
	\item If $\Gamma$ has a ribbon structure $\mu$, and if the contracted subgraph has no cycles, then $\Gamma'$ has an induced ribbon structure $\mu' : F' \ra F'$ defined by 
	$\alpha(\mu'(f)) = (\mu \circ \sigma)^m(\mu(\alpha(f)))$ where $m \geq 0$ is the smallest integer such that $(\mu \circ \sigma)^m(\mu(\alpha(f))) \in \alpha(F')$.
	\item If $\Gamma' \ra \Gamma''$ is another edge contraction, then the contraction subgraph of $\Gamma \ra \Gamma''$ has no cycles if and only if both $\Gamma \ra \Gamma'$ and $\Gamma' \ra \Gamma''$ have this property. In this case, the procedure in (2) is transitive. {\it i.e.}, starting with a ribbon structure $\mu$ on $\Gamma$, then the ribbon structure on $\Gamma''$ induced by the two-step contraction $\Gamma \ra \Gamma' \ra \Gamma''$ coincides with the one induced the composition $\Gamma \ra \Gamma''$.
\end{enumerate}
\elm

\bpf
(1) is clear. (2) and (3) can be proved by the usual fact that each ribbon forest can be embedded into $\bR^2$, uniquely up to isotopy, in a way that the cyclic ordering on vertices become the counterclockwise ordering of the embedded graph. In particular, for any edge contraction map from $\Gamma$ to $\Gamma'$ whose contraction subgraph has no cycles, let $\Lambda$ be the ribbon forest given by the contraction subgraph together with the flags attached to it. Embed $\Lambda$ into $\bR^2$ as above. Then the ribbon structure $\mu'$ in (2) is the induced counterclockwise order. 
For example, the following diagram is an example of $\Lambda$, where the contracted subgraph is drawn in bold lines:
\begin{equation*}
	\includegraphics[scale=0.4]{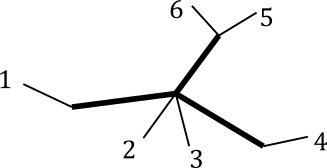}
\end{equation*}
One can verify that in this example, the definition of $\mu'$ in (2) coincides with the counterclockwise order $1 \mapsto 2 \mapsto 3 \mapsto 4 \mapsto 5 \mapsto 6 \mapsto 1$.
\epf

\bdf
Given multi-corolla $C = (V,F,p)$ and $C' = (V',F',p')$, then a \emph{merging operation} from $C$ to $C'$ is an extension of $C$ into a graph $\Gamma = (V,F,p,\sigma)$, together with an edge contraction map from $\Gamma$ to $\Gamma' = (V',F',p',\sigma' = \id)$. We may express this schematically as
\begin{equation}  \label{merge_diag_1}
	C = (V,F,p) \xdashleftarrow{\text{cut all edges}} \Gamma = (V,F,p,\sigma)
	\xra{\text{contract all edges}} C' = (V',F',p')
\end{equation}

If $C$ and $C'$ are both oriented, then we require that the corresponding structure on $C'$ to be the one induced from the edge contraction $\Gamma$ to $\Gamma'$ as in Lemma \ref{contr_ori_ribbon}. We express this in the following schematic diagrams:
\begin{equation}  \label{merge_diag_2}
 C = (V,F,p,\theta) \xdashleftarrow{\text{cut all edges}} \Gamma = (V,F,p,\sigma,\theta)
	\xra[\text{Lemma \ref{contr_ori_ribbon}}]{\text{contract all edges}} C' = (V',F',p',\theta')
\end{equation}

If $C$ and $C'$ are both oriented ribbon multi-corolla, then we furthermore require that $\Gamma$ have no cycles, so that we may apply Lemma \ref{contr_ori_ribbon}:
\begin{equation}  \label{merge_diag_3}
C = (V,F,p,\mu,\theta) \xdashleftarrow{\text{cut all edges}} \parbox{9em}{$\Gamma = (V,F,p,\sigma,\mu,\theta)$\\with no cycles}
	\xra[\text{Lemma \ref{contr_ori_ribbon}}]{\text{contract all edges}} C' = (V',F',p',\mu',\theta')
\end{equation}
In this case, any $\scO$-coloring on $\Gamma$ induces $\scO$-colorings on $C$ and $C'$. If $C$ and $C'$ are given $\scO$-colorings, then a merging operation as $\scO$-colored oriented ribbon multi-corolla is a diagram \eqref{merge_diag_3} where $\Gamma$ is endowed with a given $\scO$-coloring that induces the given ones on $C$ and $C'$.
\edf

Notice that we may compose merging operations. For example, if we are given 
\begin{equation}  \label{merge_compose}
	C \dashleftarrow \Gamma \ra C' \dashleftarrow \Gamma' \ra C''
\end{equation}
in the context of \eqref{merge_diag_1}, then the data $C' \dashleftarrow \Gamma'$ induces a lift to $\Gamma$, so that we have $C \dashleftarrow \Gamma \dashleftarrow \widetilde{\Gamma'}$. The composition of \eqref{merge_compose} is then the merging operation $C \dashleftarrow \Gamma \dashleftarrow \widetilde{\Gamma'} \ra \Gamma' \ra C''$.
This same discussion clearly carries over to the context of \eqref{merge_diag_2}. It also carries over to  the context of \eqref{merge_diag_3} because of Lemma \ref{contr_ori_ribbon}(3).

\bdf  \label{graphical_reg_pat_def}
Let $\frG_1$ be the category whose objects are rooted multi-corolla, and whose morphisms are merging operations between them (as oriented multi-corolla). Notice that in this case, the oriented graph $\Gamma$ in \eqref{merge_diag_2} is automatically a rooted forest.

Let $\frG_2$ be the category whose objects are rooted ribbon multi-corolla, and whose morphisms are merging operations between them (as oriented ribbon multi-corolla).

Let $\frG_3$ be the category whose objects are oriented multi-corolla, and whose morphisms are merging operations between them (as oriented multi-corolla) such that the oriented graph $\Gamma$ in \eqref{merge_diag_2} is a forest ({\it i.e.,} has no cycles).

Let $\frG_4$ be the category whose objects are oriented multi-corolla, and whose morphisms are merging operations between them (as oriented multi-corolla) such that the oriented graph $\Gamma$ in \eqref{merge_diag_2} has no oriented cycles.

Let $\frG_5$ be the category whose objects are multi-corolla, and whose morphisms are merging operations between them such that the graph $\Gamma$ in \eqref{merge_diag_1} is a forest ({\it i.e.,} has no cycles).

Let $\frG_6$ be the category whose objects are oriented ribbon multi-corolla, and whose morphisms are merging operations between them (as oriented ribbon multi-corolla, see \eqref{merge_diag_3}). Recall that we require $\Gamma$ to have no cycles.

Fix a set $\scO$. Let $\frG_7$ be the category whose objects are $\scO$-colored oriented ribbon multi-corolla, and whose morphisms are merging operations between them (as $\scO$-colored oriented ribbon multi-corolla).
\edf

Each of these categories has a symmetric monoidal structure obtained by disjoint union of graphs. For $i = 1,\ldots,7$, consider the subcategory $\frF_i \subset \frG_i$ with all objects and with all the invertible morphisms. Let $\frf_i \subset \frF_i$ be the full subcategory consisting of corolla. Then we have

\bpp
For each $i = 1,\ldots,7$, the data $\frf_i \subset \frF_i \ra \frG_i$ is a regular pattern. 
\epp

\bpf
It suffices to check that the map of sets \eqref{Day_restr_lax_set} for $\cA = \frF_i$ and $\cB = \frG_i$ is a bijection for each $a_1, a_2, b \in \Ob(\frG_i)$, which is obvious in each case.
\epf

\bdf
An operad over $\frG_1$ is called an \emph{operad}.

An operad over $\frG_2$ is called a \emph{non-symmetric operad}.

An operad over $\frG_3$ is called a \emph{dioperad}.

An operad over $\frG_4$ is called a \emph{properad}.

An operad over $\frG_5$ is called a \emph{cyclic operad}.

An operad over $\frG_6$ is called a \emph{ribbon dioperad}.

An operad over $\frG_7$ is called an \emph{$\scO$-colored ribbon dioperad}.
\edf

\brm
Our definition of operad/dioperad/properad/cyclic operad corresponds to the non-unital versions of these notions in the literature. 
\erm

\brm
By comparing the regular patterns $\frG_1$, $\frG_2$, $\frG_3$ and $\frG_6$, one may consider the relation between ribbon dioperads and dioperads to be similar to the relation between non-symmetric operads and operads. This explains the importance of ribbon dioperads in noncommutative algebraic geometry.
\erm

In this paper, we will focus on the cases $\frG_6$ and $\frG_7$. The categories $\frf_i \subset \frF_i \ra \frG_i$ will be renamed as $\frf \subset \frF \ra \frG$ and $\frf^{\scO} \subset \frF^{\scO} \ra \frG^{\scO}$ respectively.
Of course, $\frG \simeq \frG^{\scO}$ for $\scO = \{*\}$.

We now present $\frG^{\scO}$ by generators and relations over $\frF^{\scO}$.
Intuitively, the following should be clear:
\begin{equation}  \label{gen_rel_diop}
	\parbox{40em}{$\frG^{\scO}$ is $\otimes$-generated over $\frF^{\scO}$ by the merging operations that consist of $1$-edge contractions. Moreover, for any $2$-edge contraction whose contraction graph is of the form $[\bullet \la \bullet \la \bullet]$ or $[\bullet \ra \bullet \la \bullet]$ or $[\bullet \la \bullet \ra \bullet]$, there is a relation that says that the two ways of writing it as a composition of two $1$-edge contractions yield the same map. Thus, $\frG^{\scO}$ is $\otimes$-generated over $\frF^{\scO}$ by these generating morphisms, modulo the $\otimes$-ideal generated by these relations.}
\end{equation}

The same statement holds verbatim for the non-colored non-ribbon case, {\it i.e.,} for the regular pattern $\frG_3$ in Definition \ref{graphical_reg_pat_def} that controls dioperads. Below, we formulate \eqref{gen_rel_diop} in a more precise statement in Proposition \ref{frG_gen_rel}.

For $C_1,C_2,C_0 \in \frf^{\scO}$, let 
\begin{equation}  \label{M12_subset}
	M_{1 \la 2}(C_1,C_2,C_0) \, \subset \, \frG^{\scO}(C_1 \otimes C_2,C_0)
\end{equation}
be the subset consisting of merging operations whose contraction graph is of the shape $[1 \la 2]$ where the vertex $i$ correspond to the corolla $C_i$.
This gives a functor $M_{1 \la 2} : (\frf^{\scO})^{\op}  \times (\frf^{\scO})^{\op} \times \frf^{\scO} \ra \Set$.
Similarly, define the subsets $M_{1 \ra 2 \ra 3}(C_1,C_2,C_3,C_0)$, $M_{2 \ra 1 \la 3}(C_1,C_2,C_3,C_0)$ and $M_{1 \la 3 \ra 2}(C_1,C_2,C_3,C_0)$ of $\frG^{\scO}(C_1 \otimes C_2 \otimes C_3,C_0)$ to be those merging operations whose contraction graph is of the shape specified in the subscript. Let $M_{\text{tree}}(C_1,C_2,C_3,C_0)$ be the union of these three subsets.

\bpp  \label{frG_gen_rel}
Let $\cC$ be a symmetric monoidal category, and let $\cP_0 : \frf^{\scO} \ra \cC$ be a functor. Suppose we are given maps of sets 
\begin{equation}  \label{frG_gen_rel_alpha}
	\alpha : M_{1 \la 2}(C_1,C_2,C_0) \ra \cC(\cP_0(C_1) \otimes \cP_0(C_2) , \cP_0(C_0))
\end{equation}
functorial in $(C_1,C_2,C_0) \in (\frf^{\scO})^{\op}  \times (\frf^{\scO})^{\op} \times \frf^{\scO}$, then each element $\varphi \in M_{\text{tree}}(C_1,C_2,C_3,C_0)$ determine two well-defined%
\footnote{To define these map, one writes $\varphi$ as a composition of two morphisms in $\frG^{\op}$, and apply $\alpha$ to these morphisms. By the $\frf^{\scO}$-functoriality of $\alpha$, this is independent of the choice of the corolla that we use to represent the intermediary contraction, so that the final result is well-defined.} 
maps in $\cC$
\begin{equation*}
	\alpha'(\varphi) , \alpha''(\varphi) : \cP_0(C_1) \otimes \cP_0(C_2) \otimes \cP_0(C_3) \ra \cP_0(C_0)
\end{equation*}
defined by breaking up the $\varphi$-contraction into a $2$-step contraction by the $2$ edges. Suppose that the relation 
\begin{equation*}
	\alpha'(\varphi) \, = \, \alpha''(\varphi)
\end{equation*}
holds for each $\varphi \in M_{\text{tree}}(C_1,C_2,C_3,C_0)$, and for each $C_1,C_2,C_3,C_0 \in \frf^{\scO}$, then there exists a symmetric monoidal functor $\cP : \frG^{\scO} \ra \cC$, unique up to canonical isomorphism, whose restriction to $\frf^{\scO}$ is equal to $\cP_0$, and whose restriction to the subsets $M_{1 \la 2}(C_1,C_2,C_0) \subset \frG^{\scO}(C_1 \otimes C_2,C_0)$ is given by $\alpha$.
\epp

\bpf
For any $C',C'' \in \frF^{\scO}$, define the set $\widetilde{M}_{1\text{-edge}}(C',C'')$ by 
\begin{equation*}
	\widetilde{M}_{1\text{-edge}}(C',C'') \, = \,
	\begin{cases*}
		\frG^{\scO}(C',C'') & if $\chi(C'') = \chi(C') - 1$ \\
		\emptyset & otherwise
	\end{cases*}
\end{equation*} 

This is an $\frF^{\scO}$-bimodule $\widetilde{M}_{1\text{-edge}} : (\frF^{\scO})^{\op} \times \frF^{\scO} \ra \Set$. We may therefore define the category $T_{\frF^{\scO}}( \widetilde{M}_{1\text{-edge}} )$ with the same objects of $\frF^{\scO}$, and whose morphisms are given by
\begin{equation*}
	T_{\frF^{\scO}}( \widetilde{M}_{1\text{-edge}} ) \, := \, \frF^{\scO} \amalg (\, \widetilde{M}_{1\text{-edge}} \,) \, \amalg \, ( \, \widetilde{M}_{1\text{-edge}} \times_{\frF^{\scO}} \widetilde{M}_{1\text{-edge}} )  \, \amalg \,  \ldots
\end{equation*}
The map of $\frF^{\scO}$-bimodule $\widetilde{M}_{1\text{-edge}} \rinto \frG^{\scO}$ induces a functor 
\begin{equation*}
	\pi \, : \, T_{\frF^{\scO}}( \widetilde{M}_{1\text{-edge}} ) \raq \frG^{\scO}
\end{equation*}
which is clearly surjective on Hom sets.

If $C',C'' \in \frF^{\scO}$ are such that $\chi(C'') = \chi(C') - 2$, then the map of Hom-sets
\begin{equation}  \label{pi_map_2_edge}
	\pi  \, :  \, ( \, \widetilde{M}_{1\text{-edge}} \times_{\frF^{\scO}} \widetilde{M}_{1\text{-edge}} )(C',C'') \raq \frG^{\scO}(C',C'')
\end{equation} 
is $2$-to-$1$. In this case, the merging operations $\varphi \in \frG^{\scO}(C',C'')$ have contraction subgraph of the four possible forms:
\begin{equation}  \label{2_edge_4_types}
	[\bullet \la \bullet \la \bullet] \quad \text{or} \quad [\bullet \ra \bullet \la \bullet] \quad \text{or} \quad [\bullet \la \bullet \ra \bullet] \quad \text{or} \quad [\bullet \la \bullet \quad \bullet \la \bullet]
\end{equation}
For each of these $\varphi$, let $\widetilde{\varphi}'$ and $\widetilde{\varphi}''$ be the two preimages of $\varphi$ under \eqref{pi_map_2_edge}. Let $\sim$ be the smallest equivalence relations on the Hom sets of  $T_{\frF^{\scO}}( \widetilde{M}_{1\text{-edge}})$ containing $\widetilde{\varphi}' \sim \widetilde{\varphi}''$ and closed under composition with any other morphisms.
Then the induced functor 
\begin{equation*}
	\pi \, : \, T_{\frF^{\scO}}( \widetilde{M}_{1\text{-edge}} )/_{\sim} \raq \frG^{\scO}
\end{equation*}
is an isomorphism of categories.

Now suppose we are given the data $\cP_0$ and $\alpha$ as in the Proposition. First, $\cP_0$ extends to a symmetric monoidal functor $\cP_0^{\otimes} : \frF^{\scO} \ra \cC$, which is unique up to canonical isomorphism. The functor $\cP_0^{\otimes}$ allows one to realize $\cC$ as an $\frF^{\scO}$-bimodule. Then maps \eqref{frG_gen_rel_alpha} then canonically induces a map of $\frF^{\scO}$-bimodules $\widetilde{\alpha} : \widetilde{M}_{1\text{-edge}} \ra \cC$, which therefore induces a functor $\tau : T_{\frF^{\scO}}( \widetilde{M}_{1\text{-edge}} ) \ra \cC$. 

If $\varphi$ is of the last type of \eqref{2_edge_4_types}, then this functor $\tau$ automatically satisfies $\tau(\widetilde{\varphi}') = \tau(\widetilde{\varphi}'')$, so that it suffices to require the relation for the other three types of \eqref{2_edge_4_types}. This is precisely the relation $\alpha'(\varphi) = \alpha''(\varphi)$ in the Proposition. 
\epf

\brm  \label{frG_gen_rel_rem}
This is an example of presenting a regular pattern $\frG^{\scO}$ by generating morphisms and relations over $\frF^{\scO} \simeq \bS(\frf^{\scO})$, although we will not formulate this notion rigorously here. In general, one can take the generating morphisms to be primitive morphisms (see Definition \ref{primitive_def}). Notice that the generating morphisms in Proposition \ref{frG_gen_rel} are of the form $\alpha(\varphi) : \cP_0 \otimes \cP_0 \ra \cP_0$, while the relations are of the form $\cP_0 \otimes \cP_0 \otimes \cP_0 \ra \cP_0$. Thus, we may regard $\frG^{\scO}$ as a ``quadratic regular pattern''.
\erm

\beg  \label{ab_twist}
Let $\cC = \Ch(\bK)$, and let  $\Sigma : \cC \ra \cC$ be given by $\Sigma V = V[1]$. Fix integers $a,b \in \bZ$. Let $\frG_3$ be regular pattern that controls dioperads, as in Definition \ref{graphical_reg_pat_def}. We now define a dioperad $\cS^{a,b} : \frG_3 \ra \cC$. On $\frF_3 \subset \frG_3$, it is given by
\begin{equation*}
\cS^{a,b}(C) \, := \, \bigotimes_{f \in F^+} \Sigma^{-a} \bK \, \, \otimes \, \, \bigotimes_{f \in F^-} \Sigma^{-b} \bK \, \, \otimes \, \bigotimes_{v \in V} \Sigma^{a+b} \bK
\end{equation*}
for any oriented multi-corolla $C$, where we write $F^{\pm} = \theta^{-1}(\pm)$. Consider the isomorphism
\begin{equation}  \label{Sab_fixed_isom}
	\Sigma^{a+b} \bK \otimes \Sigma^{-a} \bK \otimes \Sigma^{-b} \bK \otimes \Sigma^{a+b} \bK \xraq{\cong} \Sigma^{a+b} \bK
\end{equation}
given by $s^{a+b}1 \otimes s^{-a}1 \otimes s^{-b}1 \otimes s^{a+b}1 \mapsto s^{a+b}1$.

Use the isomorphism \eqref{Sab_fixed_isom} to define the morphisms that are associated to $1$-edge contractions by the functor $\cS^{a,b}$. This assignment can be shown to satisfy the three types of relations in \eqref{gen_rel_diop} (for the two types $[\bullet \ra \bullet \la \bullet]$ and $[\bullet \la \bullet \ra \bullet]$, there is a non-trivial but simple check of consistency of signs).
Thus, by Proposition \ref{frG_gen_rel} (which is a precise formulation of \eqref{gen_rel_diop}), we have a dioperad $\cS^{a,b} : \frG_3 \ra \cC$.

One can then define an $\scO$-colored ribbon dioperad, still denoted as $\cS^{a,b}$, by precomposing with the obvious symmetric monoidal functor $\frG^{\scO} \ra \frG_3$.
For any $\scO$-colored dg ribbon dioperad $\cP$, its \emph{$(a,b)$-twist} is defined as the Hadamard product $\cP \otimes \cS^{a,b}$.

Notice that in the dioperad $\cS^{a,b} : \frG_3 \ra \cC$, all the edge contraction maps are invertible morphisms in $\cC$. Hence, $\cS^{a,b}$ lands in the subcategory $\cC^{\circ} \subset \cC$ of invertible morphisms. Taking inverse everywhere then gives a co-dioperad $\cS^{a,b} : \frG_3^{\op} \ra \cC$. Thus, one can also use this to define the $(a,b)$-twist of a $\scO$-colored dg ribbon co-dioperad $\cQ$.
\eeg

\beg  \label{end_diop}
We define the endomorphism ribbon dioperad $\Endrd(\cA)$ here. First, we consider a colored version of the regular pattern $\frG_4$ that controls (non-ribbon) properads. An $\scO^2$-colored graph is an graph $\Gamma = (V,F,p,\sigma)$ together with a map $\nu : F \ra \scO^2$ such that $\nu(f) = \nu(\sigma(f))$ for all $f \in F$. 
Merging operations between $\scO^2$-colored multi-corolla are required to preserve the $\scO^2$-colorings in the obvious sense. In this way, one can define the $\scO^2$-colored analogue $\frG_4^{\scO^2}$ of $\frG_4$ in Definition \ref{graphical_reg_pat_def}.

Given an oriented ribbon graph $\Gamma = (V,F,p,\sigma,\mu,\theta)$. Then any $\scO$-coloring $o : \pi_0(\partial \Sigma(\Gamma)) \ra \scO$ of $\Gamma$ induces an $\scO^2$-coloring of the underlying non-ribbon graph $ (V,F,p,\sigma)$ by the assignment $\nu(f) = (o(f''),o(f'))$ if $\theta(f) = +$ and $\nu(f) = (o(f'),o(f''))$ if $\theta(f) = -$ (see the definition of $\pi_0(\partial \Sigma(\Gamma))$ in Definition \ref{graph_def_1} for the terminology $f'$ and $f''$).
This assignment preserves merging operations, and hence gives a symmetric monoidal functor $\frG^{\scO} \ra \frG_4^{\scO^2}$.

Suppose $(\cC_1,\otimes)$ is a symmetric monoidal category, and $(\cC_2,\otimes)$ is a symmetric monoidal $\cC_1$-category. A typical example is when $(\cC_1,\otimes) = (\cC_2,\otimes)$ is closed symmetric monoidal, and the self-enrichment is given by the internal Hom.

Suppose we are given a collection $\cA$ that associates an object $\cA(x,y) \in \cC_2$ to each pair $(x,y) \in \scO^2$. Then we define a symmetric monoidal functor $\Endpr(\cA) : \frG_4^{\scO^2} \ra \cC_1$ by sending each $\scO^2$-colored oriented corolla $C$ to
\begin{equation*}
	\Endpr(\cA)(C) \, := \, \Hom \Bigl( \, \bigotimes_{f \in F^+} \, \cA(\nu(f))  \, , \, \bigotimes_{f \in F^-} \, \cA(\nu(f))  \, \Bigr)
\end{equation*}
where we write $F^{\pm} = \theta^{-1}(\pm)$,
and $\Hom(-,-)$ is the enriched Hom, so that it is an object in $\cC_1$. If $C = \coprod_{\alpha} C_{\alpha}$ is an $\scO^2$-colored oriented multi-corolla, then we take $\Endpr(\cA)(C) = \bigotimes_{\alpha} \Endpr(\cA)(C_{\alpha})$. Since the merging operations in $\frG_4^{\scO^2}$ have no oriented cycles, it is clear how to assign their image under the functor $\Endpr(\cA) : \frG_4^{\scO^2} \ra \cC_1$.

The ribbon dioperad $\Endrd(\cA) : \frG^{\scO} \ra \cC_1$ is the precomposition of $\Endpr(\cA)$ with the above defined forgetful functor $\frG^{\scO} \ra \frG_4^{\scO^2}$.
\eeg

\bdf
Let $\cP : \frG^{\scO} \ra \cC_1$ be a ribbon dioperad in a symmetric monoidal category $(\cC_1,\otimes)$. Suppose that $(\cC_2,\otimes)$ is a symmetric monoidal $\cC_1$-category, then a $\cP$-algebra in $\cC_2$ is a collection  $\cA$ that associates an object $\cA(x,y) \in \cC_2$ to each pair $(x,y) \in \scO^2$, together with a map $\cP \ra \Endrd(\cA)$ of ribbon dioperads in $\cC_1$, where $\Endrd(\cA)$ is as defined in Example \ref{end_diop}.
\edf

\bdf  \label{unital_def_1}
For any $x,y \in \scO$, let $C_{x,y}$ be the $\scO$-colored oriented ribbon corolla with two flags $F = \{ f^+,f^-\}$ with $\theta(f^{\pm}) = \pm$, and colored by $o((f^-)') = o((f^+)'') = x$ and  $o((f^-)'') = o((f^+)') = y$.
An $\scO$-colored ribbon dioperad $\cP : \frG^{\scO} \ra \cC$ is said to be \emph{unital} if for any $x,y \in \scO$, there exists a map $\id_{x,y} : \mathbf{1} \ra \cP(C_{x,y})$ in $\cC$ that acts as an identity element under $1$-edge contractions with $C_{x,y}$. 

Maps between unital $\scO$-colored ribbon dioperads are required to preserve the maps $\id_{x,y}$. The category of unital $\scO$-colored ribbon dioperads is denoted as $\OP^{u}(\frG^{\scO},\cC)$.
\edf

\blm  \label{unital_ext_1}
Let $(\cC,\otimes)$ be symmetric monoidal category with finite coproducts such that $\otimes$ preserves finite coproducts in each variable. 
The forgetful functor $\OP^{u}(\frG^{\scO},\cC) \ra \OP(\frG^{\scO},\cC)$ has a left adjoint $\cP \mapsto \cP^+$. Moreover, the underlying $\frf^{\scO}$-module of $\cP^+$ is given by
\begin{equation}  \label{P_plus_def}
	\cP^+(C) \, = \, \begin{cases*}
		 \cP(C) \amalg \mathbf{1} &  if  $C \cong C_{x,y}$ for some $x,y \in \scO$\\
		 \cP(C) & otherwise
	\end{cases*}
\end{equation}
for any $\scO$-colored ribbon corolla $C$.
\elm

\bpf
Define $\cP^+|_{\frf^{\scO}}$ by \eqref{P_plus_def}. To extend it to an $\scO$-colored ribbon dioperad $\cP^+ : \frG^{\scO} \ra \cC$, it suffices, by Proposition \ref{frG_gen_rel}, to specify the effects of $1$-edge contractions. We take the obvious choice. It is clear that for any $\cQ \in \OP^{u}(\frG^{\scO},\cC)$, we have $\Hom_{\OP^{u}(\frG^{\scO},\cC)}(\cP^+, \cQ) \cong \Hom_{\OP(\frG^{\scO},\cC)}(\cP, \cQ)$ (simply check that the obvious maps in both directions are well-defined).
\epf

\brm  \label{comb_descr_rmk}
By using Proposition \ref{frG_gen_rel}, one can give a combinatorial definition of a ribbon dioperad (in fact, there are two natural formulation of such). 
For simplicity, we will work with the uncolored case. 
The first step is to describe an $\frf$-module. Given a tuple $(n_1,\ldots,n_p) \in \bN^p$ of non-negative integers, for $p \geq 1$, one can form an oriented ribbon corolla with $|F^-| = p$ (where $F^{\pm} = \theta^{-1}(\pm)$) such that, if we name them $F^- = \{1,\ldots,p\}$ in anti-clockwise order, then there are $n_i$ flags in $F^+$ between $i$ and $i+1$. Every oriented ribbon corolla with $|F^-| \geq 1$ can be constructed in this way. The remaining case of $|F^-| = 0$ is specified by the integer $|F^+| = n$.
Under this description, an $\frf$-module is specified by the following data:
\begin{enumerate}
	\item An object $\cP(n_1,\ldots,n_p) \in \cC$ for each $p \geq 1$ and each tuple $(n_1,\ldots,n_p) \in \bN^p$ of non-negative integers;
	\item isomorphisms $\sigma_{(n_1,\ldots,n_p)} : \cP(n_1,\ldots,n_p) \xra{\cong} \cP(n_2,\ldots,n_p,n_1)$, such that 
	$\sigma_{(n_p,n_1,\ldots,n_{p-1})} \circ \ldots \circ \sigma_{(n_2,\ldots,n_p,n_1)} \circ \sigma_{(n_1,\ldots,n_p)} = \id$;
	\item an object $\cP_0(n)$ for each $n \in \bN$; and
	\item an automorphism $\sigma : \cP_0(n) \ra \cP_0(n)$ satisfying $\sigma^n = \id$.
\end{enumerate} 

Given tuples $\vec{m} = (m_1,\ldots,m_q) \in \bN^q$ and $\vec{n} = (n_1,\ldots,n_p) \in \bN^p$ for $p,q \geq 1$, then for each $1 \leq i \leq p$, $1 \leq j \leq n_i$, $1 \leq k \leq q$, define $\vec{m} *_{i,j,k} \vec{n} \in \bN^{p+q-1}$ by
\begin{equation*}
	\vec{n} *_{i,j,k} \vec{m} \, := \, (m_1,\ldots,m_{k-1}, m_k + n_i - j , n_{i+1},\ldots,n_p,n_1,\ldots,n_{i-1},j-1+m_{k+1},m_{k+2},\ldots,m_q)
\end{equation*}
Then part of the data of a ribbon dioperad consists of a map 
\begin{equation*}
	\alpha_{i,j,k} \, : \, \cP(\vec{n}) \otimes \cP(\vec{m}) \raq \cP(\vec{n} *_{i,j,k} \vec{m})
\end{equation*}
as well as certain maps on $\cP_0(n) \otimes \cP(\vec{m})$, and they are required to satisfy certain relations as specified by Proposition \ref{frG_gen_rel}. 

In another formulation, we may describe an oriented ribbon corolla by taking a tuple $\theta = (\theta_1,\ldots,\theta_n) \in \{+,-\}^n$, $n \geq 0$, and placing them in counter-clockwise order. The tuple $\theta^{\sigma} \in \{+,-\}^n$ defined by $\theta^{\sigma} := (\theta_2, \ldots,\theta_n, \theta_1)$ clearly gives an isomorphic oriented ribbon corolla. Accordingly, an $\frf$-module consists of an object $\cP(\theta) \in \cC$ for each $\theta \in \{+,-\}^n$, $n \geq 0$, together with isomorphisms $\sigma : \cP(\theta) \xra{\cong} \cP(\theta^{\sigma})$, such that 
$\sigma \circ \stackrel{(n)}{\ldots} \circ \sigma = \id$ on $\cP(\theta)$.

Given tuples $\theta \in \{+,-\}^n$ and $\theta' \in \{+,-\}^m$ with $m,n \geq 1$,
and $1 \leq i \leq n$ and $1 \leq j \leq m$ such that $\theta_i = +$ and $\theta'_j = i$, define
\begin{equation*}
	\theta *_{i,j} \theta' = (\theta_1,\ldots,\theta_{i-1},\theta'_{j+1},\ldots,\theta'_m,\theta'_1,\ldots,\theta'_{j-1},\theta_{i+1},\ldots,\theta_n)
\end{equation*}
then the data of a ribbon dioperad consists of a map
\begin{equation*}
	\alpha_{i,j} \, : \, \cP(\theta) \otimes \cP(\theta') \raq \cP(\theta *_{i,j} \theta')
\end{equation*}
which is required to satisfy certain relations as specified by Proposition \ref{frG_gen_rel} (there are relations that specify the effects of the maps $\sigma : \cP(\theta) \xra{\cong} \cP(\theta^{\sigma})$ under the composition maps $\alpha_{i,j}$, as well as $3$ kinds of associativity relations for $\alpha_{i,j}$).

In both of these formulations, we have found it rather tedious to write down the precise relations that the composition maps are required to satisfy, which makes these combinatorial formulations difficult to work with. In this paper, we will always work with the general (abstract) definition of a ribbon dioperad instead of these combinatorial formulation. Any result stated in terms of the general definition can be easily translated into these combinatorial formulations.
\erm

\section{Lie algebra structure}  \label{Lie_alg_sec}

Let $(\cC,\otimes)$ be a cocomplete symmetric monoidal additive category ({\it i.e.,} a cocomplete symmetric monoidal $\cV$-category for $(\cV,\boxtimes) = (\Ab,\otimes)$) such that $\otimes$ is cocontinuous in each variable. We also assume that finite limits in $\cC$ exist.
We are mostly interested in the case $\cC = \Ch(\bK)$, but it is convenient to work in the more general setting so that one does not have to check Koszul signs.
We will often argue as if $\cC = \Ab$, so that for $V \in \cC$ we may write $x,y \in V$ and $x-y \in V$, {\it etc.} All the statements and proofs can be rewritten categorically, so that they hold for general such $(\cC,\otimes)$.

\bdf  \label{pseudo_pre_Lie_def}
A \emph{pseudo-pre-Lie algebra} in $\cC$ is an object $V \in \cC$, together with three maps
\begin{equation*}
	\begin{split}
		m_{1 \la 2} &: V \otimes V \ra V \\
		m_{2 \ra 1 \la 3} &: V \otimes V \otimes V \ra V \\
		m_{1 \la 3 \ra 2} &: V \otimes V \otimes V \ra V
	\end{split}
\end{equation*}
We will write $x \circ y := m_{1 \la 2}(x,y)$. Then these relations are required to satisfy (see the paragraph above)
\begin{equation*}
	\begin{split}
		m_{2 \ra 1 \la 3}(x,y,z) &= m_{2 \ra 1 \la 3}(x,z,y) \\
		m_{1 \la 3 \ra 2}(x,y,z) &= m_{1 \la 3 \ra 2}(y,x,z) \\
		(x \circ y) \circ z - x \circ (y \circ z) &= m_{2 \ra 1 \la 3}(x,y,z) - m_{1 \la 3 \ra 2}(x,y,z)
	\end{split}
\end{equation*}
\edf

Notice that a pre-Lie algebra is a pseudo-pre-Lie algebra such that $m_{1 \la 3 \ra 2} = 0$, while an associative algebra is a pseudo-pre-Lie algebra such that $m_{1 \la 3 \ra 2} = m_{2 \ra 1 \la 3} = 0$.

\blm  \label{pseudo_pre_Lie_implies_Lie}
For any pseudo-pre-Lie algebra $V$, the map $[x,y] := x \circ y - y \circ x$ defines a Lie algebra structure on $V$.
\elm

\bpf
Write $f(x_1,x_2,x_3)+ \text{cyc.} := \sum_{\sigma \in C_3} f(x_{\sigma(1)}, x_{\sigma(2)}, x_{\sigma(3)})$. Then compute
\begin{equation*}
	\begin{split}
	 [[x,y],z] + \text{cyc.} 
	 &= (x \circ y) \circ z - z \circ (x \circ y) - (y \circ x) \circ z + z \circ (y \circ x)  + \text{cyc.} \\
	 &= (x \circ y) \circ z - x \circ (y \circ z) - (y \circ x) \circ z + y \circ (x \circ z)  + \text{cyc.} \\
	 &= m_{2 \ra 1 \la 3}(x,y,z) - m_{1 \la 3 \ra 2}(x,y,z) - m_{2 \ra 1 \la 3}(y,x,z) + m_{1 \la 3 \ra 2}(y,x,z) + \text{cyc.} \\
	\end{split}
\end{equation*}
which is zero because we have $m_{2 \ra 1 \la 3}(x,y,z) + \text{cyc.} = m_{2 \ra 1 \la 3}(y,x,z) + \text{cyc.}$ and $m_{1 \la 3 \ra 2}(x,y,z) = m_{1 \la 3 \ra 2}(y,x,z)$.
\epf

\brm
When $\cC = \Mod(\bK)$, then a Lie algebra object in $\cC$ is a Lie algebra in the usual sense, where the anti-symmetry axiom is defined to be $[x,y] = -[y,x]$. The other often used convention $[x,x]=0$ (equivalent if $2$ is invertible in $\bK$) is not well-defined for general additive $(\cC,\otimes)$.
\erm

For $a,b \in \bZ_{\geq -1}$, define $\frf^{\scO}_{a,b} \subset \frf^{\scO}$ to be the subcategory consisting of $\scO$-colored oriented ribbon corolla with $|F^+|-1 = a$ and $|F^-|-1 = b$. 

Let $\cP : \frG^{\scO} \ra \cC$ be a $\scO$-colored ribbon dioperad in $\cC$. Let $\lim_{\frf^{\scO}}^{\gr}(\cP)$ and $\colim_{\frf^{\scO}}^{\gr}(\cP)$ be the $(\bZ_{\geq -1})^2$-graded object in $\cC$ defined by 
\begin{equation*}
	\begin{split}
		{\rm lim}_{\frf^{\scO}}^{\gr}(\cP)_{a,b} \, &:= \, \lim \, [ \, \cP  :  \frf^{\scO}_{a,b} \ra \cC \, ]  \\
		{\rm colim}_{\frf^{\scO}}^{\gr}(\cP)_{a,b} \, &:= \, \colim \, [ \, \cP  :  \frf^{\scO}_{a,b} \ra \cC \, ] 
	\end{split}
\end{equation*}
We now construct an $(\bZ_{\geq -1})^2$-graded pseudo-pre-Lie algebra structure on ${\rm lim}_{\frf^{\scO}}^{\gr}(\cP)$ and ${\rm colim}_{\frf^{\scO}}^{\gr}(\cP)$.
The reader is advised to first refer to  \eqref{m_12_lim_diag} and  \eqref{m_12_colim_diag}, as well as the two paragraphs that follows, for a quick overview of the definition of $m_{1 \la 2}$. We now give a precise formulation of \eqref{m_12_lim_diag} and  \eqref{m_12_colim_diag}, which ensures that they are well-defined (and also clarify any possible confusion).

Let $\cM_{1 \la 2}^{\scO}$ be the category whose objects are $(C_1,C_2,C_0,\varphi)$ where $C_1,C_2,C_0 \in \frf^{\scO}$, and $\varphi$ is a morphism in $\frG^{\scO}$ from $C_1 \otimes C_2$ to $C_0$ whose contraction graph is of the shape $[1 \la 2]$ where the vertex $i$ correspond to the corolla $C_i$. A morphism in $\cM_{1 \la 2}^{\scO}$ consist of three maps $C_i \ra C_i'$ in $\frf^{\scO}$ that intertwine with $\varphi$ and $\varphi'$.
Let $\pi_i : \cM_{1 \la 2}^{\scO} \ra \frf^{\scO}$ be the functor that sends $(C_1,C_2,C_0,\varphi)$ to $C_i$.

\blm  \label{M_12_lemma}
\begin{enumerate}
	\item The category $\cM_{1 \la 2}^{\scO}$ is an equivalence relation. {\it i.e.,} it is a groupoid in which there is at most $1$ morphism between any two objects.
	\item The functor $\pi = (\pi_1,\pi_2,\pi_0) : \cM_{1 \la 2}^{\scO} \ra \frf^{\scO} \times \frf^{\scO} \times \frf^{\scO}$ is an isofibration. {\it i.e.,} for every $(C_1,C_2,C_0,\varphi) \in \cM_{1 \la 2}^{\scO}$, and every triple of isomorphisms $\alpha_i : C_i \cong C_i'$ in $\frf^{\scO}$, there exists $\varphi' \in \frG^{\scO}(C_1' \otimes C_2',C_0')$ such that $(C_1',C_2',C_0',\varphi') \in \cM_{1 \la 2}^{\scO}$ and $(\alpha_1,\alpha_2,\alpha_0)$ is a morphism in $\cM_{1 \la 2}^{\scO}$.
\end{enumerate}
\elm

\bpf
It is clear that $\cM_{1 \la 2}^{\scO}$ is a groupoid. For any object $(C_1,C_2,C_0,\varphi) \in \cM_{1 \la 2}^{\scO}$, the map $\varphi$ determines a distinguished flag in $C_1$ and $C_2$. A morphism $\alpha_i : C_i \cong C_i'$ in $\cM_{1 \la 2}^{\scO}$ must preserve this distinguish flag, so that $\alpha_1$ and $\alpha_2$ is determined by this because they must also respect the cyclic order. This determines the map on the set of all flags, which therefore also determines $\alpha_3$. This proves (1). Statement (2) is obvious: $\varphi'$ is in fact unique.
\epf

Denote by $\cP^{(2)}$ the composition $\frG^{\scO} \times \frG^{\scO} \xra{\cP \times \cP} \cC \times \cC \xra{\otimes} \cC$. By monoidicity of $\cP$, it is also isomorphic to the composition $\frG^{\scO} \times \frG^{\scO} \xra{\otimes} \frG^{\scO} \xra{\cP} \cC$. We first construct maps
\begin{equation}  \label{m12_lim_colim}
	\begin{split}
		m'_{1 \la 2} \, &: \, {\rm lim}_{\frf^{\scO}_{a_1,b_1} \times \frf^{\scO}_{a_2,b_2}}(\cP^{(2)}) \raq   {\rm lim}_{\frf^{\scO}_{a_1+a_2,b_1+b_2}}(\cP) \\
		m''_{1 \la 2} \, &: \, {\rm colim}_{\frf^{\scO}_{a_1,b_1} \times \frf^{\scO}_{a_2,b_2}}(\cP^{(2)}) \raq   {\rm colim}_{\frf^{\scO}_{a_1+a_2,b_1+b_2}}(\cP)
	\end{split}
\end{equation}

Let $(\cM_{1 \la 2}^{\scO})_{(a_1,b_1),(a_2,b_2)}$ be preimage of $\frf^{\scO}_{(a_1,b_1)} \times \frf^{\scO}_{(a_2,b_2)}$ under $\pi_1 \times \pi_2$. In particular, $\pi_0$ sends $(\cM_{1 \la 2}^{\scO})_{(a_1,b_1),(a_2,b_2)}$ to $\frf^{\scO}_{a_1+a_2,b_1+b_2}$.

For each $C_0 \in \frf^{\scO}_{a_1+a_2,b_1+b_2}$, consider the preimage subcategory 
$\pi_0^{-1}(C_0)_{(a_1,b_1),(a_2,b_2)} \subset (\cM_{1 \la 2}^{\scO})_{(a_1,b_1),(a_2,b_2)}$ (recall that morphisms are required to be sent to the identity on $C_0$ under $\pi_0$). 
Each $(C_1,C_2,C_0,\varphi) \in \pi_0^{-1}(C_0)_{(a_1,b_1),(a_2,b_2)}$ determines a map 
\begin{equation}  \label{m12_lim_1}
	(m'_{1 \la 2})_{(C_1,C_2,C_0,\varphi)} \, : \, {\rm lim}_{\frf^{\scO}_{a_1,b_1} \times \frf^{\scO}_{a_2,b_2}}(\cP^{(2)}) \raq \cP^{(2)}(C_1 \times C_2) \xraq{\cP(\varphi)} \cP(C_0)
\end{equation}
Clearly, the map depends only the connected component of $(C_1,C_2,C_0,\varphi)$ inside $\pi_0^{-1}(C_0)_{(a_1,b_1),(a_2,b_2)}$.
Hence we take
\begin{equation}  \label{m12_lim_2}
	(m'_{1 \la 2})_{C_0} \, : \, {\rm lim}_{\frf^{\scO}_{a_1,b_1} \times \frf^{\scO}_{a_2,b_2}}(\cP^{(2)}) \raq \cP(C_0)
\end{equation}
given by the sum of the maps \eqref{m12_lim_1}, one for each connected component of $\pi_0^{-1}(C_0)_{(a_1,b_1),(a_2,b_2)}$.

By Lemma \ref{M_12_lemma}(2), $\pi_0$ is an isofibration, so that each isomorphism $\psi : C_0 \xra{\cong} C_0'$ in $\frf^{\scO}$ induces a bijection between the connected components of $\pi_0^{-1}(C_0)_{(a_1,b_1),(a_2,b_2)}$ and $\pi_0^{-1}(C_0')_{(a_1,b_1),(a_2,b_2)}$. From this, we see that 
$(m'_{1 \la 2})_{C_0'} = \cP(\psi) \circ (m'_{1 \la 2})_{C_0}$, so that the map \eqref{m12_lim_2} passes to a map $m'_{1 \la 2}$ as in \eqref{m12_lim_colim}.

The construction of the map $m''_{1 \la 2}$ in \eqref{m12_lim_colim} is similar. For each $(C_1,C_2) \in \frf^{\scO}_{(a_1,b_1)} \times \frf^{\scO}_{(a_2,b_2)}$, consider the subcategory $(\pi_1 \times \pi_2)^{-1}(C_1,C_2) \subset (\cM_{1 \la 2}^{\scO})_{(a_1,b_1),(a_2,b_2)}$. Each object $(C_1,C_2,C_0,\varphi) \in (\pi_1 \times \pi_2)^{-1}(C_1,C_2)$ determines a map
\begin{equation}  \label{m12_colim_1}
	(m''_{1 \la 2})_{(C_1,C_2,C_0,\varphi)} \, : \, \cP(C_1 \otimes C_2) \xraq{\cP(\varphi)} \cP(C_0) \raq   {\rm colim}_{\frf^{\scO}_{a_1+a_2,b_1+b_2}}(\cP)
\end{equation}
Again, this map depends only on the connected component of $(C_1,C_2,C_0,\varphi)$ inside $(\pi_1 \times \pi_2)^{-1}(C_1,C_2)$, and we take
\begin{equation}  \label{m12_colim_2}
	(m''_{1 \la 2})_{C_1,C_2} \, : \, \cP(C_1 \otimes C_2) \raq   {\rm colim}_{\frf^{\scO}_{a_1+a_2,b_1+b_2}}(\cP)
\end{equation}
given by the sum of the maps \eqref{m12_colim_1}, one for each connected component of $(\pi_1 \times \pi_2)^{-1}(C_1,C_2)$.

By Lemma \ref{M_12_lemma}(2), $\pi_1 \times \pi_2$ is an isofibration, so that each $(\psi_1,\psi_2) : (C_1,C_2) \xra{\cong} (C_1',C_2')$ in $\frf^{\scO} \times \frf^{\scO}$ induces a bijection between the connected components of $(\pi_1 \times \pi_2)^{-1}(C_1,C_2)$ and $(\pi_1 \times \pi_2)^{-1}(C'_1,C'_2)$. From this, we see that $(m''_{1 \la 2})_{C_1,C_2} = (m''_{1 \la 2})_{C'_1,C'_2} \circ \cP(\psi_1 \otimes \psi_2)$, so that \eqref{m12_colim_2} pass to a map $m''_{1 \la 2}$ as in \eqref{m12_lim_colim}.

To construct the maps $m_{1 \la 2}$ on $\lim_{\frf^{\scO}}^{\gr}(\cP)$ and $\colim_{\frf^{\scO}}^{\gr}(\cP)$, notice that by the universal properties of limits and colimits, there are canonical map
\begin{equation*}
	\begin{split}
	{\rm lim}_{\frf^{\scO}_{a_1,b_1}} (\cP)  \otimes {\rm lim}_{\frf^{\scO}_{a_2,b_2}}(\cP) &\raq 
	{\rm lim}_{\frf^{\scO}_{a_1,b_1} \times \frf^{\scO}_{a_2,b_2}}(\cP^{(2)}) \\
	{\rm colim}_{\frf^{\scO}_{a_1,b_1}} (\cP)  \otimes {\rm colim}_{\frf^{\scO}_{a_2,b_2}}(\cP) &\xlaq{\cong} 
	{\rm colim}_{\frf^{\scO}_{a_1,b_1} \times \frf^{\scO}_{a_2,b_2}}(\cP^{(2)})
	\end{split}
\end{equation*}
where the map on colimits is an isomorphism by the assumption that $\otimes$ is cocontinuous in each variable.
Combining this with \eqref{m12_lim_colim}, we obtain maps
\begin{equation}  \label{m12_lim_colim_2}
	\begin{split}
		m_{1 \la 2} \, &: \, {\rm lim}_{\frf^{\scO}}^{\gr}(\cP) \otimes {\rm lim}_{\frf^{\scO}}^{\gr}(\cP) \raq  {\rm lim}_{\frf^{\scO}}^{\gr}(\cP) \\
		m_{1 \la 2} \, &: \, {\rm colim}_{\frf^{\scO}}^{\gr}(\cP) \otimes {\rm colim}_{\frf^{\scO}}^{\gr}(\cP) \raq  {\rm colim}_{\frf^{\scO}}^{\gr}(\cP)
	\end{split}
\end{equation}

This construction can be summarized by the following diagrams:
\begin{equation}  \label{m_12_lim_diag}
	\includegraphics[scale=0.4]{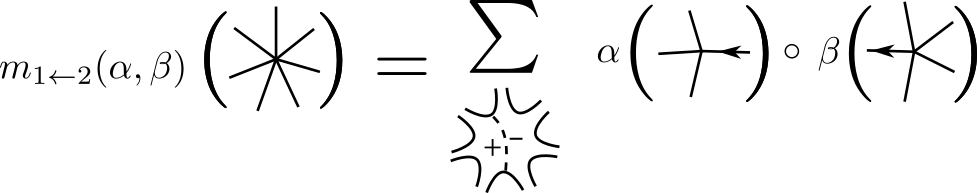}
\end{equation}
\begin{equation}  \label{m_12_colim_diag}
	\includegraphics[scale=0.4]{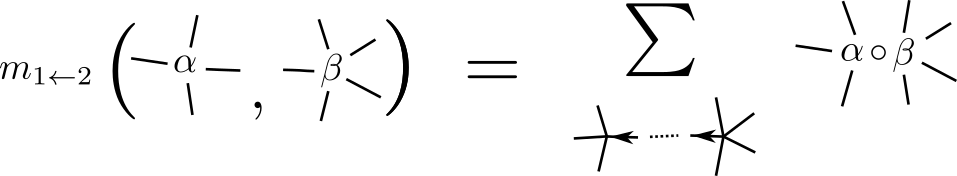}
\end{equation}

Here, \eqref{m_12_lim_diag} defines the map $m_{1 \la 2}$ for ${\rm lim}_{\frf^{\scO}}(\cP)$, while \eqref{m_12_colim_diag} defines the map $m_{1 \la 2}$ for ${\rm colim}_{\frf^{\scO}}(\cP)$. 
In these diagrams, we argue as if objects of $\cC$ had underlying sets, so we may talk about their elements. The reader may assume $\cC = \Ch(\bK)$ in the present discussion. In particular, an element $\alpha$ of ${\rm lim}_{\frf^{\scO}}(\cP)$ assigns an element $\alpha(C) \in \cP(C)$ for each $C \in \frf^{\scO}$ in a functorial way. Given $\alpha , \beta \in {\rm lim}_{\frf^{\scO}}(\cP)$, then \eqref{m_12_lim_diag} asserts that the element $m_{1 \la 2}(\alpha , \beta) \in {\rm lim}_{\frf^{\scO}}(\cP)$ assigns to $C_0 \in \frf^{\scO}$ the sum on the right hand side. Here, the sum is taken over all possible ways of writing $C_0$ as a merging of two oriented ribbon corolla. As shown in \eqref{m_12_lim_diag}, such a sum is indexed by the set $(\partial \Sigma(C_0))^2$. Namely, given any pair $(b',b'')$ of connected components of the boundary of the surface $\Sigma(C_0)$ obtained by thickening $C_0$, one can draw a line from $b'$ to $b''$, and put ``$-$'' on the right hand side, and ``$+$'' on the left hand side. This cuts $\Sigma(C_0)$ into two halves, each of which is naturally the thickening of an $\scO$-colored oriented ribbon corolla.
We apply $\alpha$ and $\beta$ to these two $\scO$-colored oriented ribbon corolla and compose along the merging edge (this composition map is part of the structure of $\cP$).

On the other hand, an element of ${\rm colim}_{\frf^{\scO}}(\cP)$ is a finite sum of pairs $(C_1,\alpha)$ where $C_1 \in \frf^{\scO}$ and $\alpha \in \cP(C_1)$. In \eqref{m_12_colim_diag}, the element $m_{1 \la 2}((C_1,\alpha) , (C_2,\beta))$ is defined to be a sum over all possible ways of pairing a ``$-$'' leg of $C_2$ with a ``$+$'' leg of $C_2$ that is compatible with the $\scO$-colorings. Any such pairing defines a merging operation $\varphi \in \frG^{\scO}(C_1 \otimes C_2, C_0)$, and we may compose $\alpha$ and $\beta$ along $\varphi$ to obtain an element $\alpha \circ \beta \in \cP(C_0)$, which is the term on the right hand side of \eqref{m_12_colim_diag}.

The construction of $m_{1 \la 2}$ can be generalized to other labelled oriented trees.
Thus, let $\Gamma$ be an oriented tree with $n$ vertices, $n \geq 2$, and with no legs, and let $\phi : \{1,\ldots,n\} \cong V(\Gamma)$ be a numbering of the vertices. One can then define the category $\cM_{(\Gamma,\phi)}^{\scO}$ in the same way as above. The analogue of Lemma \ref{M_12_lemma} can be proved in the same way.
Following the same procedure, one can construct the maps
\begin{equation}  \label{m_Gamma_maps}
	\begin{split}
		m_{(\Gamma,\phi)} \, &: \, {\rm lim}_{\frf^{\scO}}^{\gr}(\cP)^{\otimes n} \raq  {\rm lim}_{\frf^{\scO}}^{\gr}(\cP) \\
		m_{(\Gamma,\phi)} \, &: \, {\rm colim}_{\frf^{\scO}}^{\gr}(\cP)^{\otimes n} \raq  {\rm colim}_{\frf^{\scO}}^{\gr}(\cP)
	\end{split}
\end{equation}
These maps only depend on the isomorphism type of $(\Gamma,\phi)$.

If $\sigma \in S_n$ is a permutation, then changing the numbering $\phi$ by $\sigma$ changes the map $m_{(\Gamma,\phi)}$ by a precomposition by $\sigma$. {\it i.e.,} we have $m_{(\Gamma,\phi)} \circ \sigma_* = m_{(\Gamma,\phi \circ \sigma)}$, where $\sigma_*$ is the automorphism of ${\rm lim}_{\frf^{\scO}}^{\gr}(\cP)^{\otimes n}$ (resp. ${\rm colim}_{\frf^{\scO}}^{\gr}(\cP)^{\otimes n}$) given by the symmetric monoidal structure on $\cC$. In particular, $m_{2 \ra 1 \la 3}$ is symmetric in the second and third variables, while $m_{1 \la 3 \ra 2}$ is symmetric in the first and second variables.

\bthm  \label{Lie_colim_lim_diop}
The following holds for both ${\rm lim}_{\frf^{\scO}}^{\gr}(\cP)$ and ${\rm colim}_{\frf^{\scO}}^{\gr}(\cP)$ (see the paragraph preceding Definition \ref{pseudo_pre_Lie_def} for our notation):
 \begin{equation*}
 	\begin{split}
 		m_{2 \ra 1 \la 3}(x,y,z) \, &= \, m_{2 \ra 1 \la 3}(x,z,y) \\
 		m_{1 \la 3 \ra 2}(x,y,z)  \, &= \, m_{1 \la 3 \ra 2}(y,x,z) \\
 		m_{1 \la 2}( m_{1 \la 2}(x,y),z ) \, &= \, m_{1 \la 2 \la 3}(x,y,z) + m_{2 \ra 1 \la 3}(x,y,z) \\
 		m_{1 \la 2}( x, m_{1 \la 2}(y,z)) \, &= \, m_{1 \la 2 \la 3}(x,y,z) + m_{1 \la 3 \ra 2}(x,y,z)
 	\end{split}
 \end{equation*}
Hence, the maps $(m_{1 \la 2},m_{2 \ra 1 \la 3},m_{1 \la 3 \ra 2})$ give pseudo-pre-Lie algebra structures on both ${\rm lim}_{\frf^{\scO}}^{\gr}(\cP)$ and ${\rm colim}_{\frf^{\scO}}^{\gr}(\cP)$. In particular, $[x,y] := m_{1 \la 2}(x,y) - m_{1 \la 2}(y,x)$ is a Lie algebra structure on both.
\ethm

\bpf
We have already seen the validity of the first two equations in the preceding paragraph. The remaining two equations can be seen easily from the diagrams  \eqref{m_12_lim_diag} and  \eqref{m_12_colim_diag}.
\epf

Notice that the Hom-sets of $\frf^{\scO}$ are finite. Therefore, for each pair $C_1, C_2 \in \frf^{\scO}$, we may take the sum
\begin{equation*}
N_{C_1,C_2} = \sum_{\phi \in \frf^{\scO}(C_1,C_2)} \, \phi_*  \, : \, \cP(C_1) \raq \cP(C_2) 
\end{equation*}
This map is invariant under pre-compositions and post-compositions. Namely, if we are given $\phi_1 : C_1' \ra C_1$ and $\phi_2 : C_2 \ra C_2'$, then we have 
$(\phi_2)_* \circ N_{C_1,C_2} \circ (\phi_1)_* = N_{C'_1,C'_2}$. Therefore, $N_{C_1,C_2}$ descends to a map
\begin{equation}  \label{norm_map_def}
	\Nm \, : \, {\rm colim}_{\frf^{\scO}}^{\gr}(\cP) \raq {\rm lim}_{\frf^{\scO}}^{\gr}(\cP)
\end{equation}
called the \emph{norm map}.

\bthm  \label{norm_map_Lie}
The norm map \eqref{norm_map_def} intertwines with the operations $m_{1 \la 2}$, $m_{2 \ra 1 \la 3}$, $m_{1 \la 3 \ra 2}$, and is therefore a map of graded pseudo-pre-Lie algebra objects.
\ethm

\bpf
We will only give the proof for $m_{1 \la 2}$. The proofs for the other operations are completely parallel.

Both the maps $F' = \Nm \circ m_{1 \la 2}$ and $F'' = m_{1 \la 2} \circ (\Nm \otimes \Nm)$ are descended from maps of the form
\begin{equation*}
F'_{(C_1,C_2,C_0)} , \, F''_{(C_1,C_2,C_0)}  \,: \, \cP(C_1) \otimes \cP(C_2) \raq \cP(C_0)
\end{equation*}
through the universal properties of colimits and limits.
Both of these are in terms of a finite sum, described in \eqref{m_12_lim_diag} and \eqref{m_12_colim_diag}. We will verify that the sets that index these respective sums are in a canonical bijection, and the corresponding maps are the same under this bijection.

The sum $F'_{(C_1,C_2,C_0)}$ is indexed over the set $S'$ of choices, defined as follows. For each composable pair $(f^+,f^-)$ where $f^+$ is a positive flag from $C_1$ and $f^-$ is a negative flag from $C_2$, fix a choice $C_0'$ of contraction of $(C_1,C_2)$ along $(f^+,f^-)$. Then we take $S' := \coprod_{(f^+,f^-)} \frf^{\scO}(C_0',C_0)$.
There is an obvious bijection $S' \cong M_{1\la 2}(C_1,C_2,C_0)$ with the set defined in \eqref{M12_subset}.

The sum $F''_{(C_1,C_2,C_0)}$ is indexed over the set $S''$ of choices, defined as follows. For each choice of a pair of boundary $(b',b'')$ of $C_0$, fix a choice $(C'_1,C'_2,C_0,\varphi)$ that represents the corresponding contraction. Then we take $S'' := \coprod_{(b',b'')} \frf^{\scO}(C_1,C'_1) \times \frf^{\scO}(C_2,C'_2)$.
There is again an obvious bijection $S'' \cong M_{1\la 2}(C_1,C_2,C_0)$.

It is also clear that under both of these bijections, each map in the sum is simply $\phi_*$ for $\phi \in M_{1\la 2}(C_1,C_2,C_0) \subset \frG^{\scO}(C_1 \otimes C_2,C_0)$. In other words, both $F'_{(C_1,C_2,C_0)}$ and $F''_{(C_1,C_2,C_0)}$ are equal to the map
\begin{equation*}
\sum_{\phi \, \in \, M_{1\la 2}(C_1,C_2,C_0)} \, \phi_* \, : \, \cP(C_1) \otimes \cP(C_2) \raq \cP(C_0)
\end{equation*}
\epf

Now we assume that $(\cC,\otimes)$ is closed.
Let $\cP \in \OP(\frG^{\scO},\cC)$ be an $\scO$-colored ribbon dioperad and $\cQ \in \COP(\frG^{\scO},\cC)$ be an $\scO$-colored ribbon co-dioperad. 
Let $[\cQ , \cP] \in \OP(\frG^{\scO},\cC)$ be the Hadamard Hom (see the discussion at \eqref{Hadamard_2}). Then we define
\begin{equation}  \label{conv_Lie}
	[\cQ , \cP]_{\frf^{\scO}} \, := \, {\rm lim}_{\frf^{\scO}}([\cQ , \cP])
	\qquad \quad  \text{and} \qquad \quad  
	[\cQ , \cP]^{\gr}_{\frf^{\scO}} \, := \,{\rm lim}_{\frf^{\scO}}^{\gr}([\cQ , \cP])
\end{equation}
which are (graded) Lie algebra objects in $\cC$, which we call the \emph{convolution Lie algebra}.

\section{The cobar construction}  \label{cobar_sec}

Let $(\cC,\otimes)$ be a cocomplete symmetric monoidal additive category such that $\otimes$ is cocontinuous in each variable, and let $\Sigma : \cC \ra \cC$ be an odd infinitesimally symmetric monoidal endofunctor (see Definition \ref{inf_sym_mon_def} and \eqref{Sigma_odd}). 
We assume that $\Sigma$ is an autoequivalence, with an inverse $\Sigma^{-1}$, which is then automatically odd infinitesimally symmetric monoidal.
We are mostly interested in the example $\cC = \GrMod(k)$ and $\cC = \Ch(k)$ (see Example \ref{graded_module_shift_example} and \ref{chain_complex_example_2}), but we work in the more general categorical setup so that Koszul signs are automatically handled.

Let $\cQ : (\frG^{\scO})^{\op} \ra \cC$ be an $\scO$-colored ribbon co-dioperad in $\cC$. Since $\frf^{\scO}$ is a groupoid, the right $\frf^{\scO}$-module $\cQ|_{(\frf^{\scO})^{\op}}$ can be regarded as a left $\frf^{\scO}$-module, which will be denoted as $\cQ|_{\frf^{\scO}}$.
Consider the free $\scO$-colored ribbon dioperad (we omit the superscripts $\scO$ and denote the free functor simply by $\frG_{\frf} : \frf \text{-} \Mod_{\cC} \ra \OP(\frG^{\scO},\cC)$)
\begin{equation*}
	\Omega( \cQ ) \, := \, \frG_{\frf}( \Sigma^{-1} (\cQ|_{\frf^{\scO}}))
\end{equation*}

By definition as a left Kan extension, the underlying $\frf^{\scO}$-module of $\Omega( \cQ )$ may be described as the following coend in $\cC$:
\begin{equation}  \label{cobar_f_Mod_1}
	\Omega( \cQ )(C_0) \, := \, \bigoplus_{r \geq 0} \, \int^{(C_1 \otimes \ldots \otimes C_r) \in \frF^{\scO}} \frG(C_1 \otimes \ldots \otimes C_r, C_0) \boxtimes ( \Sigma^{-1} (\cQ(C_1)) \otimes \ldots \otimes \Sigma^{-1}(\cQ(C_r)))
\end{equation}
where we denote by $S \boxtimes V := V^{\oplus S}$ for $S \in \Set$ and $V \in \cC$. Also, as usual, we denote by $C_i$ oriented ribbon corolla $C_i \in \frf^{\scO}$.

Let $\cM^{(r)}(C_0)$ be the category whose objects are pairs $(C_1 \otimes \ldots \otimes C_r,\varphi)$ where $C_1 \otimes \ldots \otimes C_r$ is an object in $\bS^r(\frf^{\scO}) \subset \frF^{\scO}$, and $\varphi \in \frG(C_1 \otimes \ldots \otimes C_r , C_0)$. Morphisms are those in $\frF^{\scO}$ that intertwine with the maps $\varphi$. Then \eqref{cobar_f_Mod_1} may be rewritten as
\begin{equation}  \label{cobar_f_Mod_2}
	\Omega( \cQ )(C_0) \, = \, \bigoplus_{r \geq 0} \, \Omega^{(r)}(\cQ)(C_0)  \, := \, \bigoplus_{r \geq 0} \, \, 
	\underset{(C_1 \otimes \ldots \otimes C_r,\varphi) \in \cM^{(r)}(C_0)}{\colim} \, \,  \Sigma^{-1} (\cQ(C_1)) \otimes \ldots \otimes \Sigma^{-1}(\cQ(C_r))
\end{equation}

Recall the category $\cM_{1 \la 2}^{\scO}$ defined in the paragraph preceding Lemma \ref{M_12_lemma}. The subcategory $\pi_0^{-1}(C_0) \subset \cM_{1 \la 2}^{\scO}$ differ from $\cM^{(2)}(C_0)$ in two (related) ways. Firstly, the objects are different: the objects in $\pi_0^{-1}(C_0)$ have contraction graph of the form $[1 \la 2]$, while there is no such restriction in $\cM^{(2)}(C_0)$. Secondly, the morphisms in $\pi_0^{-1}(C_0)$ are morphisms in $\frf^{\scO} \times \frf^{\scO}$ intertwining with $\varphi$, while the morphisms in $\cM^{(2)}(C_0)$ are morphisms in $\bS^2(\frf^{\scO})$ intertwining with $\varphi$.
In any case, there is a well-defined functor $\pi_0^{-1}(C_0) \ra \cM^{(2)}(C_0)$ sending $(C_1,C_2,C_0,\varphi)$ to $(C_1 \otimes C_2, \varphi)$.

Our next goal is to use an analogue of \eqref{m_12_lim_diag} to define a differential on $\Omega( \cQ )$. An often confusing aspect of this is the sign, which we specify now. Choose a convention for defining the isomorphism
\begin{equation*}
	\psi_{V_1,V_2} \, : \, \Sigma^{-1}(V_1 \otimes V_2) \xraq{\cong} \Sigma( \Sigma^{-1}(V_1) \otimes \Sigma^{-1}(V_2) )
\end{equation*}
For example, we may require that in the target, $\Sigma$ annihilate with $\Sigma^{-1}$ that attaches to $V_1$ (rigorously, write $\psi_{V_1,V_2}$ in the form $\Sigma^{-2}(V_1 \otimes V_2) \ra \Sigma^{-1}(V_1) \otimes \Sigma^{-1}(V_2)$, and use Remark \ref{coherence_for_odd} to define the map). We call this map $\psi_{V_1,V_2}^L$. Notice that the other convention differs by a minus sign because $\Sigma$ is odd. They also differ by a conjugation:
\begin{equation}  \label{psi_L_sigma}
	\begin{tikzcd} [column sep = 70]
		\Sigma^{-1}(V_1 \otimes V_2) \ar[r, "\psi^L_{V_1,V_2}"] \ar[d, "\Sigma^{-1}(\sigma_{V_1,V_2})"'] & \Sigma(\Sigma^{-1}(V_1) \otimes \Sigma^{-1}(V_2)) \ar[d, "\Sigma(\sigma_{\Sigma^{-1}(V_1), \Sigma^{-1}(V_2)})"] \\
		\Sigma^{-1}(V_2 \otimes V_1) \ar[r, "\psi^R_{V_2,V_1} = -\psi^L_{V_2,V_1}"] & \Sigma(\Sigma^{-1}(V_2) \otimes \Sigma^{-1}(V_1)) 
	\end{tikzcd}
\end{equation}

\blm  \label{psi_V123_lemma}
For any $V_1,V_2,V_3 \in \cC$, consider the compositions
\begin{equation*}
	\begin{split}
		\psi^L_{((V_1,V_2),V_3)} \, : \, &\Sigma^{-1}(V_1 \otimes V_2 \otimes V_3) 
		\xra{\psi^L_{V_1 \otimes V_2,V_3}} \Sigma( \Sigma^{-1}(V_1 \otimes V_2) \otimes \Sigma^{-1}(V_3) ) \xra{\Sigma( \psi^L_{V_1,V_2} \otimes \id )}  \\
		&\Sigma( \Sigma ( \Sigma^{-1}(V_1) \otimes \Sigma^{-1}(V_2)) \otimes \Sigma^{-1}(V_3) ) \xra{\Sigma(\phi^{-1})} \Sigma^{2}( \Sigma^{-1}(V_1) \otimes \Sigma^{-1}(V_2) \otimes \Sigma^{-1}(V_3))  \\
		\psi^L_{(V_1,(V_2,V_3))} \, : \, &\Sigma^{-1}(V_1 \otimes V_2 \otimes V_3) 
		\xra{\psi^L_{V_1 , V_2 \otimes V_3}} \Sigma( \Sigma^{-1}(V_1)  \otimes \Sigma^{-1} (V_2 \otimes V_3) ) \xra{\Sigma( \id \otimes \psi^L_{V_2,V_3})}  \\
		&\Sigma( \Sigma^{-1}(V_1) \otimes \Sigma ( \Sigma^{-1}(V_2) \otimes \Sigma^{-1}(V_3)) ) \xra{\Sigma(\phi^{-1})} \Sigma^{2}( \Sigma^{-1}(V_1) \otimes \Sigma^{-1}(V_2) \otimes \Sigma^{-1}(V_3))  
	\end{split}
\end{equation*}
then we have $\psi^L_{((V_1,V_2),V_3)} = - \psi^L_{(V_1,(V_2,V_3))}$.
\elm

\bpf
We first give an informal argument.
Label the functors $\Sigma$ and $\Sigma^{-1}$ at the target of these two maps. For example, write the target as $\Sigma_a \Sigma_b ( \Sigma_1^{-1}(V_1) \otimes \Sigma_2^{-1}(V_2) \otimes \Sigma_3^{-1}(V_3))$. Read the maps  $\psi^L_{((V_1,V_2),V_3)}$ and $\psi^L_{(V_1,(V_2,V_3))}$ backwards, then $\psi^L_{((V_1,V_2),V_3)}$ instructs that $\Sigma_b$ annihilate with $\Sigma_1$, and then $\Sigma_a$ annihilate with $\Sigma_2$; while $\psi^L_{(V_1,(V_2,V_3))}$ instructs that $\Sigma_b$ annihilate with $\Sigma_2$, and then $\Sigma_a$ annihilate with $\Sigma_1$. These two differ by a sign. 

It is clear how to formalize this argument in our main example of interest $\cC = \GrMod(\bK)$ and $\cC = \Ch(\bK)$. For the general case, one can put $\Sigma^2$ to the left, and regard both $\psi^L_{((V_1,V_2),V_3)}$ and $\psi^L_{(V_1,(V_2,V_3))}$ as maps $\Sigma^{-3}(V_1 \otimes V_2 \otimes V_3) \ra \Sigma^{-1}(V_1) \otimes \Sigma^{-1}(V_2) \otimes \Sigma^{-1}(V_3)$, and apply Remark \ref{coherence_for_odd}.
\epf

Once we have fixed this convention $\psi^L$, we define $d : \Omega( \cQ )|_{\frf^{\scO}} \ra \Sigma( \Omega( \cQ )|_{\frf^{\scO}} )$ to be the unique derivation (see Lemma \ref{der_free_operad}) whose restriction to $\Sigma^{-1} (\cQ|_{\frf^{\scO}})$ is the map
\begin{equation}  \label{d_on_Q}
 d \, : \, \Sigma^{-1} (\cQ|_{\frf^{\scO}}) \raq \Sigma (\Omega^{(2)}(\cQ)|_{\frf^{\scO}} )
\end{equation}
defined by an analogue of \eqref{m_12_lim_diag}. Namely, for each connected component of $\pi_0^{-1}(C_0)$, choose a representative $(C_1,C_2,C_0,\varphi)$, and consider the map
\begin{equation*}  
d_{(C_1,C_2,C_0,\varphi)} \, : \, \Sigma^{-1}(\cQ(C_0)) \xraq{ \Sigma^{-1}(\cQ(\varphi)) } \Sigma^{-1}(\cQ(C_1) \otimes \cQ(C_2)) \xraq{\psi^L} \Sigma( \Sigma^{-1}(\cQ(C_1)) \otimes \Sigma^{-1}(\cQ(C_2)) )
\end{equation*}
We will regard this map as landing in $\Sigma( \Omega^{(2)}(\cQ)(C_0) )$ by putting it to the component $(C_1 \otimes \ldots \otimes C_r, \varphi)$ in \eqref{cobar_f_Mod_2}. As such, $d_{(C_1,C_2,C_0,\varphi)} : \Sigma^{-1}(\cQ(C_0)) \ra \Sigma( \Omega^{(2)}(\cQ)(C_0) )$ is independent of the choice of the representative $(C_1,C_2,C_0,\varphi)$ in each connected component of $\pi_0^{-1}(C_0)$. 
Define \eqref{d_on_Q} to be the sum of these maps, one for each connected component of $\pi_0^{-1}(C_0)$ (there are finitely many of them).

\bthm  \label{cobar_d_square_zero}
The derivation $d : \Omega( \cQ )|_{\frf^{\scO}} \ra \Sigma( \Omega( \cQ )|_{\frf^{\scO}} )$ satisfies $\Sigma(d) \circ d = 0$. 
\ethm

\bpf
$\Sigma(d) \circ d$ is a derivation by Lemma \ref{square_zero_deriv_lemma}. In order to show that it is zero, it suffices to consider its restriction to $\Sigma^{-1}(\cQ(C_0))$. 
Thus, we compute the composition
\begin{equation}  \label{d_square_cobar_Q}
	 \Sigma^{-1} (\cQ|_{\frf^{\scO}}) \xraq{d} \Sigma (\Omega^{(2)}(\cQ)|_{\frf^{\scO}} ) \xraq{\Sigma(d)} (\Omega^{(3)}(\cQ)|_{\frf^{\scO}} ) 
\end{equation}
By the description \eqref{cobar_f_Mod_2} of $\Omega^{(2)}(\cQ)|_{\frf^{\scO}}$, we see that $\Sigma(d)$ can be applied to either the first or second tensor product. Hence we may write $\Sigma(d) = \Sigma(d \otimes \id + \id \otimes d)$.
The composition \eqref{d_square_cobar_Q} can be described in terms of the structure maps in $\cQ$ and of the maps in Lemma \ref{psi_V123_lemma} (we will write $V_i = \cQ(C_i)$ for $i = 1,2,3$ below). Recall the categories $\cM^{\scO}_{2 \ra 1 \la 3}$, $\cM^{\scO}_{1 \la 3 \ra 2}$ and $\cM^{\scO}_{1 \la 2 \la 3}$ appearing in Section \ref{Lie_alg_sec}.
First, the map $\Sigma(d \otimes \id) \circ d$ is a sum of maps of two types:
\begin{enumerate}
	\item[(1)] For each connected component in $\pi_0^{-1}(C_0) \subset \cM^{\scO}_{1 \la 2 \la 3}$, choose a representative $(C_1,C_2,C_3,\varphi)$. Then there is a term
	\begin{equation}  \label{d_square_cobar_Q_term_1}
	\begin{split}
		&\Sigma^{-1}(\cQ(C_0)) \xraq{ \Sigma^{-1}(\cQ(\varphi)) } \Sigma^{-1}(\cQ(C_1) \otimes \cQ(C_2) \otimes \cQ(C_3)) \xraq{\psi^L_{(V_1,V_2),V_3}} \\
		& \Sigma( \Sigma^{-1}(\cQ(C_1)) \otimes \Sigma^{-1}(\cQ(C_2)) \otimes \Sigma^{-1}(\cQ(C_3)) )
	\end{split}
	\end{equation}
\item[(2)] For each connected component in $\pi_0^{-1}(C_0) \subset \cM^{\scO}_{2 \ra 1 \la 3}$, choose a representative $(C_1,C_2,C_3,\varphi)$. Then there is a term \eqref{d_square_cobar_Q_term_1}.
\end{enumerate}

Similarly, the map $\Sigma(\id \otimes d) \circ d$ is a sum of maps of two types:
\begin{enumerate}
	\item[(3)] For each connected component in $\pi_0^{-1}(C_0) \subset \cM^{\scO}_{1 \la 2 \la 3}$, choose a representative $(C_1,C_2,C_3,\varphi)$. Then there is a term \eqref{d_square_cobar_Q_term_1} with $\psi^L_{(V_1,V_2),V_3}$ replaced by $\psi^L_{V_1,(V_2,V_3)}$.
	\item[(4)] For each connected component in $\pi_0^{-1}(C_0) \subset \cM^{\scO}_{1 \la 3 \ra 2}$, choose a representative $(C_1,C_2,C_3,\varphi)$. Then there is a term \eqref{d_square_cobar_Q_term_1} with $\psi^L_{(V_1,V_2),V_3}$ replaced by $\psi^L_{V_1,(V_2,V_3)}$.
\end{enumerate}

By Lemma \ref{psi_V123_lemma}, the terms (1) and terms (3) cancel each other. 
Notice that terms in (2) come in pairs, where one switch the role of $C_2$ and $C_3$. By \eqref{psi_L_sigma} and Lemma \ref{psi_V123_lemma}, these two terms cancel each other. Similarly, terms in (4) come in pairs that cancel each other.
\epf

By Lemma \ref{square_zero_deriv_lemma}, this derivation $d$ then induces a $\frG^{\scO}$-operad in $\COM(\cC,\Sigma)$. 
We will mostly apply this to the case $\cC = \COM(\cC_0,\Sigma)$ (think of $\cC_0 = \GrMod(\bK)$ and $\cC = \Ch(\bK)$). 
Recall from Example \ref{chain_complex_example_2} that an object in $\COM(\cC,\Sigma)$ is an object $V \in \cC_0$ together with a pair $(d_1,d_2)$ of anti-commuting square-zero differentials.
Thus, if we start with $(\cQ,d_1) \in \COP(\frG^{\scO}, \cC )$, then Theorem \ref{cobar_d_square_zero} gives $(\Omega(\cQ),d_1,d_2) \in \COP(\frG^{\scO}, \COM(\cC,\Sigma) )$. We will take $d = d_1 + d_2$, so that $(\Omega(\cQ),d) \in \COP(\frG^{\scO}, \cC )$.

\bdf  \label{cobar_def}
Suppose that $\cC = \COM(\cC_0,\Sigma)$. For an $\scO$-colored ribbon co-dioperad $(\cQ,d_1) \in \COP(\frG^{\scO}, \cC )$, the $\scO$-colored ribbon dioperad $(\Omega(\cQ),d) \in \COP(\frG^{\scO}, \cC )$, $d = d_1 + d_2$, is called the \emph{cobar ribbon dioperad} of $\cQ$.
\edf

\bdf
Suppose that $\cC = \COM(\cC_0,\Sigma)$. For a Lie algebra object $L$ in $\cC$,
and for a map $x : \mathbf{1} \ra \Sigma(L)$ in $\cC_0$, denote by $d(x)$ the map $\mathbf{1} \xra{x} \Sigma(L) \xra{\Sigma(d)} \Sigma^2(L)$, and by 
$[x,x]$ the map $\mathbf{1} =  \mathbf{1} \otimes \mathbf{1} \xra{x \otimes x} \Sigma(L) \otimes \Sigma(L) \cong \Sigma^2(L \otimes L) \xra{\Sigma^2([-,-])} \Sigma^2(L)$. A \emph{Maurer-Cartan element} is such a map $x$ that satisfies $2d(x) + [x,x] = 0$.
Denote by $\MC(L)$ the set of Maurer-Cartan elements of $L$. 
\edf

\bdf  \label{twisting_mor_def}
Suppose that $\cC = \COM(\cC_0,\Sigma)$. Given $\cQ \in \COP(\frG^{\scO}, \cC )$ and $\cP \in \OP(\frG^{\scO}, \cC )$, then a \emph{twisting morphism} from $\cQ$ to $\cP$ is a Maurer-Cartan element in the convulution Lie algebra $[\cQ,\cP]_{\frf}$ (see \eqref{conv_Lie}). Denote the set of twisting morphisms from $\cQ$ to $\cP$ by
\begin{equation*}
	\Tw(\cQ,\cP) \, := \, \MC( [\cQ,\cP]_{\frf} )
\end{equation*}
\edf

\bthm  \label{twisting_cobar_adj}
If $2$ in invertible in $\bK$, then there is a canonical bijection
\begin{equation*}
	\Tw(\cQ,\cP) \, \cong \, \Hom_{\OP(\frG^{\scO},\cC)}(\Omega(\cQ),\cP)
\end{equation*}
\ethm

\bpf
First, notice that we have
\begin{equation*}
	\Hom_{\cC_0}( \mathbf{1},  \Sigma([\cQ,\cP]_{\frf^{\scO}})) \, \cong \, 
	\Hom_{\frf^{\scO}\text{-}\Mod_{\cC_0}}( \Sigma^{-1}(\cQ|_{\frf^{\scO}}) , \cP|_{\frf^{\scO}} ) \, \cong \, \Hom_{\OP(\frG^{\scO},\cC_0)}(\Omega(\cQ),\cP)
\end{equation*}
Indeed, the right hand side follows from the definition of $\Omega(\cQ)$ as a free $\frG^{\scO}$-operad in $\cC_0$, while the left hand side follows from the fact that $\cQ|_{\frf^{\scO}}$ is defined the left $\frf^{\scO}$-module obtained from the right $\frf^{\scO}$-module $\cQ|_{(\frf^{\scO})^{\op}}$ by using the inverse involution on the groupoid $\frf^{\scO}$.
%

Thus, it suffices to check that the Maurer-Cartan equation on the left hand side corresponds to the compatibility with differentials on the right hand side (it suffices to check compatibility with differentials on the generator $\Sigma^{-1}(\cQ|_{\frf^{\scO}})$ of $\Omega(\cQ)$). 
As in Definition \ref{cobar_def}, write the differential on $\Omega(\cQ)$ as $d = d_1 + d_2$. If $F \in \Hom_{\cC_0}( \mathbf{1},  \Sigma([\cQ,\cP]_{\frf^{\scO}}))$ corresponds to $\widetilde{F} \in \Hom_{\OP(\frG^{\scO},\cC_0)}(\Omega(\cQ),\cP)$, then $d(F)$ is given by $(d_{\cP} \circ \widetilde{F} - \widetilde{F} \circ d_1)|_{\frf^{\scO}}$, while $[F,F] = 2 m_{1 \la 2}(F,F)$, where $m_{1 \la 2}(F,F)$ is given by 
$(\widetilde{F} \circ d_2)|_{\frf^{\scO}}$.
\epf

\section{Modular ribbon properads}  \label{mrp_sec}

\bdf
A \emph{modular ribbon graph} $\widetilde{\Gamma}$ is a ribbon graph $\Gamma = (V,F,p,\sigma,\mu)$ together with a finite set $W$, a map of sets $c : V \ra W$ and a function $g : W \ra \bN$. 

In a modular ribbon graph, the set $V$ is called the set of \emph{simple vertices}; while the set $W$ is called the set of \emph{compound vertices}.

Denote by $q : F \ra W$ the composition $q = c \circ p$. The \emph{underlying graph} of a modular ribbon graph is the graph $\Gamma_W = (W,F,q, \sigma)$. Notice that this is not the underlying graph of $\Gamma$, and in general has no induced ribbon structure. 

An \emph{orientation} of $\widetilde{\Gamma}$ is an orientation of $\Gamma$, or equivalently, of $\Gamma_W$.
\edf

Each modular ribbon graph has an associated surface given by the connected sum
\begin{equation}  \label{assoc_surface_modular}
	\Sigma(\widetilde{\Gamma}) \, := \, \Sigma(\Gamma) \underset{c}{\#} \coprod_{w \in W} \Sigma_{g(w)}
\end{equation}
Here, $\Sigma_{g(w)}$ is the genus $g(w)$ closed oriented surface. We form a connected sum between the surface $\Sigma(\Gamma)$ and $\coprod_{w \in W} \Sigma_{g(w)}$ as follows: each $v \in V$ gives a point on the surface $\Sigma(\Gamma)$, which we will still denote as $v$. Choose a point $\phi(v) \in \Sigma_{g(w)}$ on the surface corresponding to $w = c(v)$. Choose them so that $\phi(v) \neq \phi(v')$ for any distinct $v,v' \in c^{-1}(w)$. Then form a connected sum between $v \in \Sigma(\Gamma)$ and $\phi(v) \in \Sigma_{g(w)}$, one for each $v \in V$ ({\it i.e.}, remove small $\epsilon$-balls around $v$ and $\phi(v)$, and then glue along the obtained boundary).

The following combinatorially defined notions are motivated by the construction \eqref{assoc_surface_modular}:

\bdf
The \emph{set of connected components} of $\widetilde{\Gamma}$ is the set of connected components of the underlying graph $\Gamma_W$. Clearly, each connected component determines a modular ribbon subgraph $\widetilde{\Gamma}' \subset \widetilde{\Gamma}$, so that $\widetilde{\Gamma}$ is the disjoint union of its connected components.

The \emph{set of boundary components} and the \emph{set of compactified boundary components} of $\widetilde{\Gamma}$ is defined to be that of $\Gamma$, as defined in Definition \ref{graph_def_1}. {\it i.e.}, $\pi_0(\partial \Sigma(\widetilde{\Gamma}) ) := \pi_0(\partial \Sigma(\Gamma) )$ and $\pi_0(\partial \overline{\Sigma}(\widetilde{\Gamma}) ) := \pi_0(\partial \overline{\Sigma}(\Gamma) )$.

If $\widetilde{\Gamma}$ is connected, then its \emph{Euler characteristics} is given by 
\begin{equation*}
	\chi( \Sigma(\widetilde{\Gamma}) ) \, := \, - |V(\Gamma)| - |E(\Gamma)| + \sum_{w \in W} ( 2 - 2g(w))
\end{equation*}
Its \emph{genus} is the non-negative integer $g(\Sigma(\widetilde{\Gamma})) \in \bN$ determined by
\begin{equation*}
	\chi( \Sigma(\widetilde{\Gamma}) ) \, = \, 2 - 2g(\Sigma(\widetilde{\Gamma})) - | \pi_0(\partial \overline{\Sigma}(\widetilde{\Gamma}) ) | 
\end{equation*}
\edf

Now we define an edge contraction operation for modular ribbon graphs.
Namely, given a modular ribbon graph $\widetilde{\Gamma}$, and given a subset $E_0 \subset E(\Gamma)$, we now define a contracted modular ribbon graph $\widetilde{\Gamma}' = (V',F',p',\sigma',\mu',W',c',g')$.
Below, we will denote the subset $\bigcup_{e \in E_0} e$ of $F$ also as $E_0$.

Let $\Gamma_{E_0}$ be the ribbon graph obtained from $\Gamma = (V,F,p,\sigma,\mu)$ by cutting open every edge not in $E_0$. {\it i.e.}, $\Gamma_{E_0} = (V,F,p,\sigma_{E_0},\mu)$, where $\sigma_{E_0}(f) = \sigma(f)$ if $f \in E_0$, and $\sigma_{E_0}(f) = f$ if $f \notin E_0$.

We will now define $\widetilde{\Gamma}'$. We will first define the ribbon multi-corolla $(V',F',p',\mu')$, then the modular ribbon multi-corolla $(V',F',p',\mu',W',c',g')$, and finally the map $\sigma'$. In fact, $\sigma'$ is easy to define in the onset, but we postpone its definition for expository reasons, because of its increasing complexity

In the definition of $\widetilde{\Gamma}'$ that we proceed to below, it is helpful to keep in mind the dependence on various data:
\begin{enumerate}
	\item The ribbon multi-corolla $(V',F',p',\mu')$ depends only on $\Gamma_{E_0} = (V,F,p,\sigma_{E_0},\mu)$.
	\item The ribbon graph $\Gamma' = (V',F',p',\sigma',\mu')$ depends only on $\Gamma = (V,F,p,\sigma,\mu)$. But $\Gamma'$ is not the contraction of $\Gamma$ by $E_0$. Instead, we will call $\Gamma'$ the \emph{dual contraction} of $\Gamma$ by $E_0$, and denote it by  $\Gamma' = \Gamma /\!/ E_0$.
	\item The modular ribbon multi-corolla $(V',F',p',\mu',W',c',g')$ depends only on $\Gamma_{E_0}$ and $(W,c,g)$.
\end{enumerate}

Let $V' := \pi_0(\partial \overline{\Sigma}( \Gamma_{E_0} ))$. Let $F' := L(\Gamma_{E_0}) = F \setminus E_0$. 

For any ribbon graph $\Lambda$, there is a well-defined map $F(\Lambda) \ra \pi_0(\partial \overline{\Sigma}( \Lambda ))$ defined by passing to the $\sigma \circ \mu$ orbit. In particular, restricting to the set of legs, we have $L(\Lambda) \ra \pi_0(\partial \overline{\Sigma}( \Lambda ))$.
Apply this to $\Lambda = \Gamma_{E_0}$, we obtain a map $p' : F' \ra V'$.

Let $\mu' : F' \ra F'$ be defined by $\mu'(f) = (\mu \circ \sigma)^m(\mu(f))$ where $m \geq 0$ is the smallest integer such that $(\mu \circ \sigma)^m(\mu(f)) \in F'$.

This defines the ribbon multi-corolla $(V',F',p',\mu')$ (except that we haven't proved that it is indeed a ribbon multi-corolla, see Proposition \ref{dual_contr_is_ribbon} below). Notice that, although our definition of $\mu'$ involves $\sigma$, we obtain the same permutation if $\sigma$ is replaced by $\sigma_{E_0}$. Thus, what we have defined so far depends only on $\Gamma_{E_0}$.
To define the ribbon graph $\Gamma'$, simply take $\sigma' := \sigma|_{F'}$.

Let $W' := W/_{\sim_{E_0}}$, where $\sim_{E_0}$ is the equivalence relation on $W$ generated by $q(f) \sim_{E_0} q(\sigma(f))$ if $f \in E_0$.
For $w' \in W$, we denote by  $W(w') \subset W$ the equivalence class.

The ribbon graph $\Gamma_{E_0}$ can be written as a disjoint union $\Gamma_{E_0} = \coprod_{w' \in W'} \Gamma_{E_0}(w')$, where $\Gamma_{E_0}(w')$ is the union of connected components consisting of vertices $v \in V$ such that $c(v) \in W(w')$.
This gives a decomposition $V' = \coprod_{w' \in W'} \pi_0(\partial \overline{\Sigma}( \Gamma_{E_0}(w') ))$, which therefore gives a map $c' : V' \ra W'$. 

The ribbon graph $\Gamma_{E_0}(w')$ can be extended to a modular ribbon graph  $\widetilde{\Gamma_{E_0}(w')}$ by taking $W(w')$ and the same maps $c$ and $g$.
Define $g' : W' \ra \bN$ by
\begin{equation*}
	g'(w') \, := \, g(\Sigma( \widetilde{\Gamma_{E_0}(w')} ) )
\end{equation*}

\blm  \label{boundary_bijection_lem}
There is a canonical bijection%
\footnote{We haven't proved that $\Gamma'$ is a ribbon graph yet, so that strictly speaking we should not write $\pi_0(\partial \overline{\Sigma}( \widetilde{\Gamma}' ))$ yet. Here, we simply define it to be $(V' \setminus p'(F')) \amalg (F'/(\sigma' \circ \mu'))$.}
\begin{equation*}
 \pi_0(\partial \overline{\Sigma}( \widetilde{\Gamma} )) := (V \setminus p(F)) \amalg (F/(\sigma \circ \mu))  \xraq[\cong]{\Phi}  (V' \setminus p'(F')) \amalg (F'/(\sigma' \circ \mu')) =: \pi_0(\partial \overline{\Sigma}( \widetilde{\Gamma}' ))
\end{equation*}
Moreover, if $b \in \pi_0(\partial \overline{\Sigma}( \widetilde{\Gamma} ))$, denote by $|b| \subset F$ the corresponding $\sigma \circ \mu$-orbit (we define $|b| = \emptyset$ if $b \in V \setminus p(F)$). Then we have $|\Phi(b)| = |b| \setminus E_0$, with the induced cyclic order as a subset $|\Phi(b)| \subset |b|$.
\elm

\bpf
On $V \setminus p(F)$, it is the obvious injection $\Phi : V \setminus p(F) \rinto V' \setminus p'(F')$. On $F/(\sigma \circ \mu)$, if $\vec{f} = (f_1,\ldots,f_r)$ is an orbit under $\sigma \circ \mu$, then remove all the elements in $E_0$, then notice that $\mu'$ is defined in such a way that the remaining set $\vec{f}' = (f_{i_1},\ldots,f_{i_s})$, with the induced cyclic order, is an orbit under $\sigma' \circ \mu'$ whenever it is non-empty. If $\vec{f}'$ is non-empty, then $\Phi$ sends $\vec{f}$ to $\vec{f}'$. If $\vec{f}'$ is empty, then by definition $\vec{f}$ represents an element of $V'$ with no flags attached to it, and we define $\Phi(\vec{f}) \in V'\setminus p'(F')$ to be this element.
\epf

\bpp  \label{dual_contr_is_ribbon}
The data $\widetilde{\Gamma}' = (V',F',p',\sigma',\mu',W',c',g')$ is a modular ribbon graph. Moreover, the associated surfaces \eqref{assoc_surface_modular} of $\widetilde{\Gamma}$ and $\widetilde{\Gamma}'$ are diffeomorphic.
\epp

\bpf
The only non-trivial content of the first statement is that $\mu$ is a ribbon structure; {\it i.e.}, for each $v' \in V'$, the subset $p^{\prime -1}(v') \subset F'$ is either empty of consists of a single $\mu'$-orbit. 
Since the statement is about $(V',F',p',\mu')$, it only depends on $\Lambda = \Gamma_{E_0}$.
Applying Lemma \ref{boundary_bijection_lem} to the dual contraction of the ribbon graph $\Lambda = \Gamma_{E_0}$ by $E_0$, then it says that $V' = (V' \setminus p'(F')) \amalg (F'/\mu')$, which is precisely the statement we want to prove.

For the second statement, notice that Lemma \ref{boundary_bijection_lem} shows that $\overline{\Sigma}(\widetilde{\Gamma}')$ and $\overline{\Sigma}(\widetilde{\Gamma})$ have bijective boundary components, and the bijection preserves the set of legs $\overline{\Sigma} \setminus \Sigma$ on each (even preserving the cyclic orders on them). Thus it suffices to show that each connected component have the same Euler characteristics. Our definition of $g'(w')$ is designed so that this holds.
\epf

\bdf
The modular ribbon graph $\widetilde{\Gamma}' = (V',F',p',\sigma',\mu',W',c',g')$ defined above is called the \emph{dual contraction of $\widetilde{\Gamma}$ by $E_0$}. It is denoted as $\widetilde{\Gamma}' = \widetilde{\Gamma} / \! / E_0$.
If $\widetilde{\Gamma}$ is oriented, then $\widetilde{\Gamma}' = \widetilde{\Gamma} / \! / E_0$ is understood to inherit the induced orientation.
\edf

\blm  \label{two_step_dual_contr}
Suppose that $E_0 = E_1 \amalg E_2$, then $E_2$ induces a set of edges in $\widetilde{\Gamma} /\!/ E_1$, which we will continue to denote as $E_2$. There is a canonical isomorphism $\widetilde{\Gamma} /\!/ E_0 \cong (\widetilde{\Gamma} /\!/ E_1) /\!/ E_2$.
\elm

\bpf
Write $\widetilde{\Gamma}'_0 = \widetilde{\Gamma} /\!/ E_0$ and $\widetilde{\Gamma}'_{12} = (\widetilde{\Gamma} /\!/ E_1) /\!/ E_2$. 
First, we show that the underlying ribbon multi-corolla $(V'_0,F'_0,p'_0,\mu'_0)$ and $(V'_{12},F'_{12},p'_{12},\mu'_{12})$ are isomorphic. Notice that both of these depend only on $\Lambda = \Gamma_{E_0}$. 
It is clear that $F'_0 = F'_{12} = F \setminus E_0$.
Apply Lemma \ref{boundary_bijection_lem} to the dual contraction of $\Gamma_{E_0}$ by $E_1$ then gives an isomorphism $\Phi : V'_0 \xra{\cong} V'_{12}$, which moreover intertwines with $p'_0$ and $p'_{12}$. 
Moreover, Lemma \ref{boundary_bijection_lem}, when applied to the dual contraction $\Gamma_{E_0} /\!/ E_0$ and the two-step dual contraction $\Gamma_{E_0} /\!/ E_0$, gives a description of orbits of both $\mu'_{0}$ and $\mu'_{12}$ as orbits of $\mu$ with elements in $E_0$ removed. This shows that $\mu'_{0} = \mu'_{12}$.
The identifications $W'_0 = W'_{12}$ and $\Phi \circ c'_0 = c'_{12} \circ \Phi$ and $g'_0 = g'_{12}$ are straightforward. 
\epf

\bdf
Given modular ribbon graphs $\widetilde{\Gamma}$ and $\widetilde{\Gamma}'$,
an \emph{edge contraction map} from $\widetilde{\Gamma}$ to $\widetilde{\Gamma}'$ consists of a subset $E_0 \subset E(\Gamma)$, together with an isomorphism $\widetilde{\Gamma} /\!/ E_0 \cong \widetilde{\Gamma}'$.
If $\widetilde{\Gamma}$ and $\widetilde{\Gamma}'$ are oriented, then we require the isomorphism $\widetilde{\Gamma} /\!/ E_0 \cong \widetilde{\Gamma}'$ to preserve the orientations.

A \emph{modular ribbon multi-corolla} $\widetilde{C}$ is defined analogous to a modular ribbon graph, except that we do not specify $\sigma$. In other words, $\widetilde{C}$ consists of a ribbon multi-corolla $C = (V,F,p,\mu)$ together with a finite set $W$, a map of sets $c : V \ra W$ and a function $g : W \ra \bN$. 
A \emph{modular ribbon corolla} is a modular ribbon multi-corolla where $W = \{*\}$.

Given oriented modular ribbon multi-corolla $\widetilde{C} = (V,F,p,\theta,\mu,W,c,g)$ and $\widetilde{C'} = (V',F',p',\theta',\mu',W',c',g')$, then a \emph{merging operation} from $\widetilde{C}$ to $\widetilde{C'}$ is an extension of $\widetilde{C}$ into an oriented modular ribbon graph $\widetilde{\Gamma} = (V,F,p,\sigma,\theta,\mu,W,c,g)$, together with an edge contraction map from $\widetilde{\Gamma}$ to $\widetilde{\Gamma}'$, where $\widetilde{\Gamma}'$ is defined by the data $\widetilde{C}'$ with $\sigma = \id$. Notice that this contraction map must necessarily contract all edges of $\widetilde{\Gamma}$. 

\begin{equation}  \label{merge_diag_4}
	\begin{tikzcd}[column sep = 0]
		& \widetilde{\Gamma} = (V,F,p,\sigma,\mu,\theta,W,c,g) \ar[ld, dashed, "\text{cut all edges}"'] \ar[rd, "\text{dual contract all edges}"] \\
	\widetilde{C} = (V,F,p,\mu,\theta,W,c,g) & & \widetilde{C}' = (V',F',p',\mu',\theta',W',c',g')
	\end{tikzcd}
\end{equation}
\edf

Any $\scO$-coloring of $\widetilde{\Gamma}$ (this simply means an $\scO$-coloring of the ribbon graph $\Gamma$) induces an $\scO$-coloring of both $\widetilde{C}$ and $\widetilde{C}'$. If both $\widetilde{C}$ and $\widetilde{C}'$ are $\scO$-colored, then we require the merging operation \eqref{merge_diag_4} to respect that, in the sense that $\widetilde{\Gamma}$ is $\scO$-colored in a way that induces those given $\scO$-colorings on $\widetilde{C}$ and $\widetilde{C}'$.

By Lemma \ref{two_step_dual_contr}, it is clear that one can compose merging operations. We will say that a merging operation \eqref{merge_diag_4} has \emph{no oriented cycles} if the underlying graph $\Gamma_W$ of $\widetilde{\Gamma}$ has no oriented cycles. Since a dual edge contraction induces an ordinary contraction edge on the underlying graph $\Gamma_W$, it is clear that this condition is closed under successive dual edge contractions. Thus, merging operations with no oriented cycles are closed under compositions.

\bdf
Let $\frG^{\scO}_{{\rm mrp}}$ be the category whose objects are $\scO$-colored oriented modular ribbon multi-corolla, and whose morphisms are merging operations (respecting the $\scO$-coloring and the orientation) that have no oriented cycles.

Let $\frF^{\scO}_{{\rm mrp}} \subset \frG^{\scO}_{{\rm mrp}}$ be the subcategory with all objects and with the invertible morphisms. Let $\frf^{\scO}_{{\rm mrp}} \subset \frF^{\scO}_{{\rm mrp}}$ be the full subcategory consisting of $\scO$-colored oriented modular ribbon corolla.
\edf
 
\bpp
The data $\frf^{\scO}_{{\rm mrp}} \subset \frF^{\scO}_{{\rm mrp}} \ra \frG^{\scO}_{{\rm mrp}}$ is a regular pattern.
\epp

\bpf
The map \eqref{Day_restr_lax_set} is clearly a bijection.
\epf

\bdf
An operad over $\frG^{\scO}_{{\rm mrp}}$ is called an \emph{$\scO$-colored modular ribbon properad} (or \emph{$\scO$-colored modular fat properad}).
\edf

One can describe $\frG^{\scO}_{{\rm mrp}}$ by generators and relations in a fashion similar to \eqref{gen_rel_diop}. We will say that a merging operation between oriented modular ribbon graphs is \emph{elementary} if the contraction graph has the shape $[\bullet \Leftarrow \bullet]$. Here, the dots refer to composite vertices ({\it i.e.,} elements in $W$) and the arrow ``$\Leftarrow$'' means there is at least one edge between the two composite vertices (all oriented in the indicated direction). Use the same notation for other shapes.

Then, one can show that $\frG^{\scO}_{{\rm mrp}}$ is $\otimes$-generated over $\frF^{\scO}_{{\rm mrp}}$ by the elementary merging operations. Moreover, for any merging operation whose contraction graph is of the form
\begin{equation*}
	[
	\begin{tikzcd} [column sep = 14]
		\bullet & \bullet \ar[l, Rightarrow] & \bullet \ar[l, Rightarrow] \ar[ll, bend right, Rightarrow]
	\end{tikzcd}
    ]
    \quad \text{or} \quad
    	[
    \begin{tikzcd} [column sep = 14]
    	\bullet & \bullet \ar[l, Rightarrow] & \bullet \ar[l, Rightarrow]
    \end{tikzcd}
    ]
    \quad \text{or} \quad
    [
    \begin{tikzcd} [column sep = 14]
    	\bullet  \ar[r, Rightarrow] &  \bullet & \bullet \ar[l, Rightarrow]
    \end{tikzcd}
    ]
    \quad \text{or} \quad
       [
    \begin{tikzcd} [column sep = 14]
    	\bullet   &  \bullet \ar[l, Rightarrow] \ar[r, Rightarrow] & \bullet 
    \end{tikzcd}
    ]
\end{equation*}
there is a relation that says that the two ways of writing it as a composition of two elementary merging operations are the same.
We leave it to the reader to formulate a precise statement, similar to the formulation of \eqref{gen_rel_diop} in terms of Proposition \ref{frG_gen_rel}.

\beg  \label{end_mrp}
Recall from Example \ref{end_diop} the regular pattern $\frG_4^{\scO^2}$ and the symmetric monoidal functor $\Endpr(\cA) : \frG_4^{\scO^2} \ra \cC_1$. There is an obvious symmetric monoidal functor $\frG^{\scO}_{{\rm mrp}} \ra \frG_4^{\scO^2}$ obtained by sending a modular ribbon graph $\widetilde{\Gamma}$ to its underlying graph $\Gamma_W$. Compose these two functors, and we have an $\scO$-colored modular ribbon properad $\Endmrp(\cA) : \frG^{\scO}_{{\rm mrp}} \ra \cC_1$.
\eeg

\bdf
Let $\cP : \frG^{\scO}_{{\rm mrp}} \ra \cC_1$ be an $\scO$-colored modular ribbon properad in a symmetric monoidal category $(\cC_1,\otimes)$. Suppose that $(\cC_2,\otimes)$ is a symmetric monoidal $\cC_1$-category, then a $\cP$-algebra in $\cC_2$ is a collection  $\cA$ that associates an object $\cA(x,y) \in \cC_2$ to each pair $(x,y) \in \scO^2$, together with a map $\cP \ra \Endmrp(\cA)$ of modular ribbon properads in $\cC_1$, where $\Endmrp(\cA)$ is as defined in Example \ref{end_mrp}.
\edf

We have the analogues of Definition \ref{unital_def_1} and Lemma \ref{unital_ext_1}.

\bdf
Let $C_{x,y}$ be as in Definition \ref{unital_def_1}. Extend it into an $\scO$-colored modular ribbon corolla $\widetilde{C}_{x,y}$ by taking $W = \{*\}$ and $g(*) = 0$. 

An $\scO$-colored modular ribbon properad $\cP : \frG_{{\rm mrp}}^{\scO} \ra \cC$ is said to be \emph{unital} if for any $x,y \in \scO$, there exists a map $\id_{x,y} : \mathbf{1} \ra \cP(\widetilde{C}_{x,y})$ in $\cC$ that acts as an identity element under $1$-edge contractions.

Maps between unital $\scO$-colored modular ribbon properads are required to preserve the maps $\id_{x,y}$. The category of unital $\scO$-colored modular ribbon properads is denoted as $\OP^{u}(\frG_{{\rm mrp}}^{\scO},\cC)$.
\edf

\blm  \label{unitalization_mrp}
Let $(\cC,\otimes)$ be symmetric monoidal category with finite coproducts such that $\otimes$ preserves finite coproducts in each variable. 
The forgetful functor $\OP^{u}(\frG_{{\rm mrp}}^{\scO},\cC) \ra \OP(\frG_{{\rm mrp}}^{\scO},\cC)$ has a left adjoint $\cP \mapsto \cP^+$. Moreover, the underlying $\frf_{{\rm mrp}}^{\scO}$-module of $\cP^+$ is given by
\begin{equation*}  
	\cP^+(\widetilde{C}) \, = \, \begin{cases*}
		\cP(\widetilde{C}) \amalg \mathbf{1} &  if  $\widetilde{C} \cong \widetilde{C}_{x,y}$ for some $x,y \in \scO$\\
		\cP(\widetilde{C}) & otherwise
	\end{cases*}
\end{equation*}
for any $\scO$-colored modular ribbon corolla $\widetilde{C}$.
\elm

To every unital $\scO$-colored modular ribbon properad $\cP$ in $\cC$, we now associate an $\scO^2$-colored PROP in $\cC$, which will be denoted as $\Cob(\cP)$.
We assume that $\cC$ is cocomplete, and $\otimes : \cC \times \cC \ra \cC$ is cocontinuous in each variable.

The object set of $\Cob(\cP)$ is 
\begin{equation*}
	\Ob(\Cob(\cP)) \,:= \, \coprod_{n \geq 0} (\scO^2)^n
\end{equation*} 
For tuples $(x,y) = ((x_1,y_1),\ldots,(x_n,y_n))$ and $(x',y') = ((x'_1,y'_1),\ldots,(x'_m,y'_m))$, define $\widetilde{\frF}_{{\rm mrp}}((x,y),(x',y'))$ to be the groupoid whose objects are pairs $(\widetilde{C},\varphi)$, where $\widetilde{C} \in \frF_{{\rm mrp}}^{\scO}$ is an $\scO$-colored modular ribbon multi-corolla, and $\varphi = (\varphi^+, \varphi^-)$ is a pair of bijections $\varphi^+ : \{1\ldots,n\} \xra{\cong} F^+$ and  $\varphi^- : \{1\ldots,m\} \xra{\cong} F^-$, where $F^{\pm} = \theta^{-1}(\pm)$, such that $\nu( \varphi^+(i) ) = (x_i,y_i)$ and $\nu( \varphi^-(i) ) = (x'_i,y'_i)$, where $\nu : F \ra \scO^2$ is as defined in Example \ref{end_diop} ({\it i.e.,} $\nu(f) = (o(f''),o(f'))$ if $f \in F^+$ and $\nu(f) = (o(f'),o(f''))$ if $f \in F^-$). A morphism in $\widetilde{\frF}_{{\rm mrp}}((x,y),(x',y'))$ is a morphism in $\frF_{{\rm mrp}}^{\scO}$ intertwining with $\varphi$.
Define
\begin{equation*}
	\Cob(\cP)((x,y),(x',y')) \, := \, \colim_{(\widetilde{C},\varphi) \in \widetilde{\frF}_{{\rm mrp}}((x,y),(x',y'))} \, \cP(\widetilde{C})
\end{equation*}

Compositions are defined in the obvious way. Since $\cP$ is assumed to be unital, this composition has units, and hence defines a $\cC$-enriched category $\Cob(\cP)$. In fact, it is a symmetric monoidal $\cC$-category. On object set, the monoidal product is defined by concatenating two tuples. On Hom-objects, it is induced by taking the disjoint union of labelled $\scO$-colored modular ribbon multi-corolla.

\bdf  \label{cob_cat_def}
The symmetric monoidal $\cC$-category $\Cob(\cP)$ is called the \emph{cobordism category} of $\cP$.
\edf

\bpp  \label{cob_cat_prop}
Let $\Endmrp(\cA) : \frG^{\scO}_{{\rm mrp}} \ra \cC_1$ be as in Example \ref{end_mrp}. Notice that it is unital.
Let $\cP \in \OP^{u}(\frG_{{\rm mrp}}^{\scO},\cC_1)$ be a unital $\scO$-colored modular ribbon properad in $\cC_1$. Then a map $\cP \ra \Endmrp(\cA)$ in $\OP^{u}(\frG_{{\rm mrp}}^{\scO},\cC_1)$ is equivalent to a symmetric monoidal $\cC_1$-functor $\Cob(\cP) \ra \cC_2$ that sends the object $(x,y)$ to $\cA(x,y)$.
\epp

\bpf
It will be convenient to use an equivalent, but larger, model for $\Cob(\cP)$. Namely, define a symmetric monoidal $\cC_1$-category $\Cob'(\cP)$ whose objects are pairs $(S,\nu)$, where $S$ is a finite set, and $\nu: S \ra \scO^2$ is a map of sets. Hom-objects are defined in a similar way, and the monoidal structure is defined by disjoint union. There is a symmetric monoidal equivalence $\Cob(\cP) \xra{\simeq} \Cob'(\cP)$.

A symmetric monoidal $\cC_1$-functor $F : \Cob'(\cP) \ra \cC_2$ is equivalent to the data 
\begin{enumerate}
	\item An object $\cA(S,\nu) \in \cC_2$ for any $(S,\nu) \in \Ob(\Cob'(\cP))$.
	\item A map $\cP(\widetilde{C}) \ra [\cA(S,\nu) , \cA(S',\nu')]$ for any $\widetilde{C} \in \frF_{{\rm mrp}}^{\scO}$ and $\varphi = (\varphi^+,\varphi^-)$, where $\varphi^+ : S \xra{\cong} F^+$ and $\varphi^+ : S' \xra{\cong} F^-$ are bijections intertwining with $\nu$ and $\nu'$ with the maps $\nu : F^{\pm} \ra \scO^2$ defined as above.
	\item An isomorphism $\cA((S,\nu)\amalg (S',\nu')) \cong \cA(S,\nu) \otimes  \cA(S',\nu')$, satisfying the associativity constraints.
\end{enumerate}
The data (2) is required to be
\begin{enumerate}
	\item[(a)] functorial under isomorphism $\widetilde{C} \cong \widetilde{C}'$ in $\frF_{{\rm mrp}}^{\scO}$;
	\item[(b)] functorial under merging operations;
	\item[(c)] compatible with the monoidal structure in (3).
\end{enumerate}
Since $\cP$ is unital, the functoriality requirement (b), applied to the ``identity cobordism'' implies the functoriality of the data (2) under change of labelling $S \cong T$ and $S' \cong T'$. Thus, in (2), it suffices to consider $S = F^+$, $S = F^-$ and $\varphi^{\pm} = \id$. Also, by (c), it suffices to specify the data (2) for $\scO$-colored modular ribbon corolla $\widetilde{C}$. The relevant data is then precisely a map of the underlying $\frf_{{\rm mrp}}^{\scO}$-modules $\cP|_{\frf_{{\rm mrp}}^{\scO}} \ra \Endmrp(\cA)|_{\frf_{{\rm mrp}}^{\scO}}$. The condition (b) then translates to this map being functorial with respect to the merging operations.  
\epf

We now prove an analogue of Theorem \ref{Lie_colim_lim_diop}. We will only work with the colimit version. 

Let $(\frG_{{\rm mrp}}^{\scO})_{a,b} \subset \frG_{{\rm mrp}}^{\scO}$ and $(\frF_{{\rm mrp}}^{\scO})_{a,b} \subset \frF_{{\rm mrp}}^{\scO}$ be the full subcategories consisting of $\scO$-colored modular ribbon multi-corolla such that
\begin{equation*}
	|L^+(\Gamma)| - \chi(\Sigma(\widetilde{\Gamma})) = a
	\qquad \text{and} \qquad 
	|L^-(\Gamma)| - \chi(\Sigma(\widetilde{\Gamma})) = b
\end{equation*}
where $L^{\pm}(\Gamma) = \theta^{-1}(\pm) \cap L(\Gamma)$, and let $(\frf_{{\rm mrp}}^{\scO})_{a,b} = (\frF_{{\rm mrp}}^{\scO})_{a,b} \cap \frf_{{\rm mrp}}^{\scO}$. We have the following obvious

\blm  \label{bigrading_mrp_lemma}
\begin{enumerate}
	\item $(\frF_{{\rm mrp}}^{\scO})_{a_1,b_1} \otimes (\frF_{{\rm mrp}}^{\scO})_{a_2,b_2} \subset (\frF_{{\rm mrp}}^{\scO})_{a_1+a_2,b_1+b_2}$.
	\item $\frG_{{\rm mrp}}^{\scO} = \coprod_{a,b} \, (\frG_{{\rm mrp}}^{\scO})_{a,b}$. {\it i.e.,} there are no morphisms between objects with different bigradings.
\end{enumerate}
\elm

Assume that $(\cC,\otimes)$ is a cocomplete symmetric monoidal additive category such that $\otimes$ is cocontinuous. Let $\cP$ be an $\scO$-colored modular ribbon properad in $\cC$. Define 
\begin{equation*}
	\begin{split}
		{\rm colim}_{\frf_{{\rm mrp}}^{\scO}}^{\gr}(\cP)_{a,b} \, &:= \, \colim \, [ \, \cP  :  (\frf_{{\rm mrp}}^{\scO})_{a,b} \ra \cC \, ] \\
			{\rm colim}_{\frF_{{\rm mrp}}^{\scO}}^{\gr}(\cP)_{a,b} \, &:= \, \colim \, [ \, \cP  :  (\frF_{{\rm mrp}}^{\scO})_{a,b} \ra \cC \, ]
	\end{split}
\end{equation*}
From Lemma \ref{bigrading_mrp_lemma}(1), we see that
\begin{equation}  \label{F_mrp_colim_Sym}
	{\rm colim}_{\frF_{{\rm mrp}}^{\scO}}^{\gr}(\cP) \, \cong \, \Sym \, ( \, {\rm colim}_{\frf_{{\rm mrp}}^{\scO}}^{\gr}(\cP) \, )
\end{equation}
as a bigraded object in $\cC$.

For each $n \geq 0$, define the map
\begin{equation*}
	m_n  \, : \, {\rm colim}_{\frF_{{\rm mrp}}^{\scO}}^{\gr}(\cP) \, \otimes \, {\rm colim}_{\frF_{{\rm mrp}}^{\scO}}^{\gr}(\cP) \raq {\rm colim}_{\frF_{{\rm mrp}}^{\scO}}^{\gr}(\cP)
\end{equation*}
by the formula
\begin{equation*}
	m_n|_{\cP(C_1) \otimes \cP(C_2)} \, = \, \sum_{\alpha = [C_1 \stackrel{(n)}{\Longleftarrow} C_2]} \cP(\alpha) 
\end{equation*}
where the sum is over all the possible ways of pairing $p$ elements in $L^+(\Gamma_1)$ with $p$ elements in $L^-(\Gamma_2)$ in a way that is compatible with the $\scO$-coloring. If we call such a pairing $\alpha$, then it determines a merging operation from $C_1 \otimes C_2$ to $C_0$, which is well-defined up to isomorphism. Applying the functoriality of $\cP$ then determines a map $\cP(C_1) \otimes \cP(C_2) \ra \cP(C_0)$. Its post-composition with $\cP(C_0) \ra {\rm colim}_{\frF_{{\rm mrp}}^{\scO}}^{\gr}(\cP)$ is independent of any choice, and is denoted as $\cP(\alpha)$.
One can obtain a precise formulation of this along the lines of \eqref{m_12_colim_diag}, the details of which will be omitted.

\bpp
\begin{enumerate}
	\item The map $m_0$ is the symmetric product in \eqref{F_mrp_colim_Sym}.
	\item The formal expression $m = \sum_{n \geq 0} m_n \hbar^n$ is associative. {\it i.e.,} for each $n \geq 0$, we have
	\begin{equation} 
		\sum_{p = 0}^n \, m_{n-p}(m_{p}(x,y),z) = \sum_{q = 0}^n \, m_{n-q}(x,m_{q}(y,z))
	\end{equation}
where we use the notation as discussed in the paragraph preceding Definition \ref{pseudo_pre_Lie_def}.
\end{enumerate}
\epp

\bpf
(1) is obvious. For (2), notice that both sides, when restricted to $\cP(C_1) \otimes \cP(C_2) \otimes \cP(C_3)$, are sum of the maps
\begin{equation*}
	\cP \Biggl( 	\Bigl[
	\begin{tikzcd}
		C_1 & C_2 \ar[l, Rightarrow, "(p)"] & C_3 \ar[l, Rightarrow, "(q)"] \ar[ll, bend right, Rightarrow, "(r)"']
	\end{tikzcd}
	\Bigr] \Biggr)
\end{equation*}
where $p \geq 0$, $q \geq 0$, $r \geq 0$, and $p + q + r = n$.
\epf

In the case $\cC = \Ch(\bK)$, the expression $m = \sum_{n \geq 0} m_n \hbar^n$ is well-defined on ${\rm colim}_{\frF_{{\rm mrp}}^{\scO}}^{\gr}(\cP)[\hbar]$, since $m_n(x,y) = 0$ for $n \gg 0$ for any given pair $(x,y)$. This gives a $\bK[\hbar]$-linear dg algebra structure on ${\rm colim}_{\frF_{{\rm mrp}}^{\scO}}^{\gr}(\cP)[\hbar]$.

\section{Pre-Calabi-Yau categories}  \label{PCY_sec}

We work with $\cC_0 = \GrMod(\bK)$ and $\cC = \Ch(\bK)$, and we consider Example \ref{end_diop}, where we regard $\cC$ as being enriched over itself via the internal Hom. Thus, $\cA$ is a collection that associates a chain complex $\cA(x,y) \in \Ch(\bK)$ for every pair $x,y \in \scO$, and we form the $\scO$-colored dg ribbon dioperad $\Endrd(\cA)$.

Take the $(1,m)$-twist $\Endrd(\cA) \otimes \cS^{1,m}$ as in Example \ref{ab_twist}, and
apply Theorem \ref{Lie_colim_lim_diop} to obtain $(\bZ_{\geq -1})^2$-graded dg Lie algebra structures on ${\rm lim}_{\frf^{\scO}}^{\gr}(\Endrd(\cA) \otimes \cS^{1,m} )$ and ${\rm colim}_{\frf^{\scO}}^{\gr}(\Endrd(\cA) \otimes \cS^{1,m} )$, as well as a norm map between them:
\begin{equation}  \label{norm_map_Endrd}
	\Nm \, : \, {\rm colim}_{\frf^{\scO}}^{\gr}(\Endrd(\cA) \otimes \cS^{1,m} ) \raq {\rm lim}_{\frf^{\scO}}^{\gr}(\Endrd(\cA) \otimes \cS^{1,m} )
\end{equation}

Notice that each $\frf^{\scO}_{a,b}$ is a groupoid, where the automorphism groups of objects are finite cyclic groups. Moreover, if $\scO$ is finite, then the groupoid has finitely many connected components. Clearly, if $\bK \supset \bQ$, then the norm map from cyclic group coinvariants to invariants is an isomorphism. Hence we have the following

\blm  \label{char_zero_Nm_lemma}
If $\bK \supset \bQ$ and $\scO$ is finite, then the norm map \eqref{norm_map_Endrd} is an isomorphism of bigraded dg Lie algebras.
\elm

\brm  \label{char_zero_Nm_rmk}
A similar statement holds if $\scO$ is infinite, provided that we modify ${\rm colim}_{\frf^{\scO}}^{\gr}(\cP )$ by replacing a direct sum over $\scO$-indexings into a direct product. We will denote this modification $\widehat{{\rm colim}}_{\frf^{\scO}}^{\gr}(\cP )$.
\erm

\bdf  \label{PCY_def}
Let $\Lambda^{{\rm PCY}} \subset (\bZ_{\geq -1})^2$ be the sub-semigroup given by $\Lambda^{{\rm PCY}} = \{  (a,b) \in (\bZ_{\geq -1})^2 \, | \, b \geq 0 \text{ and } (a,b) \neq (-1,0),(0,0)\}$.

Assume that $\bK$ is a field of characteristic $0$. Let $m = 2-n$, then an \emph{$n$-pre-Calabi-Yau category} over the object set $\scO$ consists of a chain complex $\cA(x,y) \in \Ch(\bK)$ for every pair $x,y \in \scO$, together with a Maurer-Cartan element in the dg Lie algebra
\begin{equation}  \label{PCY_DGLA_1}
	\prod_{(a,b) \in \Lambda^{{\rm PCY}}}  \, {\rm lim}_{\frf^{\scO}}^{\gr}(\Endrd(\cA) \otimes \cS^{1,m} )_{a,b}
\end{equation}
\edf


We now unravel Definition \ref{PCY_def} into a more elementary form. For notational simplicity, we consider $\scO = \{*\}$ (in which case a pre-Calabi-Yau category is called a pre-Calabi-Yau algebra), and write $A = \cA(*,*)$.
We may describe the underlying $\frf$-module of $\Endrd(A)$ as in Remark \ref{comb_descr_rmk}.
Then for any $b \geq 0$, let $p = b+1$, and we have
\begin{equation}  \label{Hoch_limt_End}
	\prod_{a \in \bZ_{\geq -1}}  \, {\rm lim}_{\frf^{\scO}}^{\gr}(\Endrd(\cA) \otimes \cS^{1,m} )_{a,b} \, = \, \left( \prod_{(n_1,\ldots,n_p) \in \bN^p} \Homcom_{\bK}(A[1]^{\otimes n_1} \otimes \ldots \otimes A[1]^{\otimes n_p} , A[-m]^{\otimes p})[m+1] \right)^{C_p}
\end{equation}
Thus the data of a pre-Calabi-Yau algebra consists of an infinite sum $\pi = \pi_1 + \pi_2 + \ldots$, where $\pi_p$ is a degree $-1$ element in \eqref{Hoch_limt_End} for $p = b+1$, satisfying the Maurer-Cartan equation. 

Since we are always restricting to the $b \geq 0$ part here, the $b = 0$ part forms a Lie subalgebra, while the $b > 0$ part forms a Lie ideal. Thus, $\pi_1$ itself satisfies the Maurer-Cartan equation. Since we are working over the range $\Lambda^{{\rm PCY}}$, we may write $\pi_1 = \mu_2 + \mu_3 + \ldots$, where $\mu_n$ is a degree $-2$ element in $\Homcom_{\bK}(A[1]^{\otimes n} , A)$. Unravelling the definition, we see that $\pi_1 = \mu_2 + \mu_3 + \ldots$ is an $A_{\infty}$-algebra structure on the chain complex $A$. 
The part $\pi_{\geq 2} = \pi_2 + \pi_3 + \ldots$ is then an extra structure, called a \emph{pre-Calabi-Yau structure} on the $A_{\infty}$-algebra $(A,\pi_1)$.

By the same consideration, we see that Theorem \ref{Lie_alg_intro} is a consequence of Theorems \ref{Lie_colim_lim_diop} and \ref{norm_map_Lie}:

\bpf[Proof of Theorem \ref{Lie_alg_intro}]
For any $A_{\infty}$-algebra $(A,\pi_1)$, use the Maurer-Cartan element $\pi_1$ to twist the differential of the dg Lie algebra \eqref{PCY_DGLA_1}. This gives a twisted differential on \eqref{Hoch_limt_End} for each fixed $p$. Notice that the right hand side of \eqref{Hoch_limt_End} is isomorphic to (an $m$-shifted version of) $C_H^{(p)}(A)^{C_p}$, as defined in \eqref{higher_HH_intro_2}. One can verify that the twisted differential on \eqref{Hoch_limt_End} coincides with the Hochschild differential in \eqref{higher_HH_intro_2}. Similarly, $C_H^{(p)}(A)_{C_p}$ can be identified with the colimit total object with twisted differential.
\epf

We now give some reformulations of the notion of a pre-Calabi-Yau category. Let $\bK_{{\rm const}}^c \in \COP(\frG^{\scO},\cC)$ be the $\scO$-colored ribbon co-dioperad given by sending every multi-corolla to $\bK$, and every merging operation to the identity. Let $\bK^c_{\Lambda^{{\rm PCY}}}$ be the quotient of $\bK^c_{{\rm const}}$ obtained by setting $\bK^c_{{\rm const}}(C)$ to zero if the multi-corolla $C$ contains any corolla not in $\Lambda^{{\rm PCY}}$. This is a co-ideal since  $\Lambda^{{\rm PCY}} \subset (\bZ_{\geq -1})^2$ is a sub-semigroup, so that we have $\bK^c_{\Lambda^{{\rm PCY}}} \in \COP(\frG^{\scO},\cC)$. 
Thus, a pre-Calabi-Yau structure is precisely a twisting morphism from $\bK^c_{\Lambda^{{\rm PCY}}}$ to $\cP = \Endrd(\cA) \otimes \cS^{1,m}$ (see Definition \ref{twisting_mor_def}), or equivalently, a twisting morphism from $\bK^c_{\Lambda^{{\rm PCY}}}  \otimes \cS^{-1,-m} $ to $\cP = \Endrd(\cA)$.

Let $\cP_{{\rm PCY}} := \Omega( \bK^c_{\Lambda^{{\rm PCY}}}  \otimes \cS^{-1,-m} )$. 
Then applying Theorem \ref{twisting_cobar_adj}, we see that an $n$-pre-Calabi-Yau category is precisely a map of dg ribbon dioperad $\cP_{{\rm PCY}} \ra \Endrd(\cA)$. In other words, it is an algebra over $\cP_{{\rm PCY}}$.

One can also reformulate it in terms of an $\scO$-colored PROP. Consider the map $\tau : \frG^{\scO} \ra \frG_{{\rm mrp}}^{\scO}$ of regular patterns. Denote by 
$\Indrm : \OP( \frG^{\scO} , \cC ) \ra  \OP( \frG_{{\rm mrp}}^{\scO} , \cC ) $ the extension functor.
Then an $n$-pre-Calabi-Yau category can also be formulated as an algebra over $\Indrm(\cP_{{\rm PCY}})$.

Apply the unitalization construction in Lemma \ref{unitalization_mrp} to $\Indrm(\cP_{{\rm PCY}})$, and consider its cobordism category $\Cob( \Indrm(\cP_{{\rm PCY}})^+ )$ (see Definition \ref{cob_cat_def}). Then by Proposition \ref{cob_cat_prop}, we see that an $n$-pre-Calabi-Yau category is the same as a symmetric monoidal functor $\Cob( \Indrm(\cP_{{\rm PCY}})^+ ) \ra \Ch(\bK)$.

\brm
We expect the $\scO^2$-colored PROP $\Cob( \Indrm(\cP_{{\rm PCY}})^+ )$ (or a suitable modification of it) to play the role of a pre-Calabi-Yau analogue of the combinatorial model considered in \cite{Cos07} for the symmetric monoidal dg category that governs open topological conformal field theory (TCFT). This question is also studied in \cite{KTV}.
\erm

In Definition \ref{PCY_def}, we assumed that ${\rm char}(\bK)=0$. In this case, we can either apply ${\rm lim}_{\frf^{\scO}}^{\gr}$ or $\widehat{{\rm colim}}_{\frf^{\scO}}^{\gr}$ to $\cP = \Endrd(\cA) \otimes \cS^{1,m} $ in order to define pre-Calabi-Yau categories, since these two are isomorphic in characteristic zero (see Lemma \ref{char_zero_Nm_lemma} and Remark \ref{char_zero_Nm_rmk}).
Concretely, taking ${\rm lim}_{\frf^{\scO}}^{\gr}$ gives the $C_p$-invariants in \eqref{Hoch_limt_End} while taking $\widehat{{\rm colim}}_{\frf^{\scO}}^{\gr}$ gives $C_p$-coinvariants. In positive characteristics, these two are different from each other. Moreover, there are two other natural choices, corresponding to taking homotopy $C_p$-invariants and homotopy $C_p$-coinvariants. 
In order to define the two latter versions of pre-Calabi-Yau categories, one have to choose specific representatives of these homotopy (co)invariants, so that one still has a dg Lie algebra structure.

Let $\bK_{{\rm const}}^c \in \COP(\frG^{\scO},\cC)$ be the constant $\scO$-colored ribbon co-dioperad as above. Similarly, let $\bK_{{\rm const}} \in \OP(\frG^{\scO},\cC)$ be the $\scO$-colored ribbon dioperad with constant value $\bK$. Then for any $\cP \in \OP(\frG^{\scO},\cC)$, we may write
\begin{equation*}
	{\rm lim}_{\frf^{\scO}}^{\gr}(\cP) = [\bK_{{\rm const}}^c , \cP]^{\gr}_{\frf^{\scO}} 
	\qquad \text{and} \qquad 
	{\rm colim}_{\frf^{\scO}}^{\gr}(\cP) = {\rm colim}_{\frf^{\scO}}^{\gr}(\bK_{{\rm const}} \otimes \cP)
\end{equation*}
where the latter is the Hadamard product.

Thus, to define the homotopy invariant version, it suffices to resolve $\bK_{{\rm const}}^c \in \COP(\frG^{\scO},\cC)$ so that the underlying $\frf^{\scO}$-module is projective. Similarly, to define the homotopy coinvariant version, it suffices to resolve $\bK_{{\rm const}} \in \OP(\frG^{\scO},\cC)$ so that the underlying $\frf^{\scO}$-module is projective.

We now give a solution to a set-theoretic version of the latter problem. Namely, consider  $\ast_{{\rm const}} \in \OP(\frG^{\scO},\Set)$ to be the $\scO$-colored ribbon dioperad with constant value $\{*\} \in \Set$. We now resolve it within $\OP(\frG^{\scO},\Setdel)$. 
Consider the functor
\begin{equation*}
\cQ_0^{{\rm pre}} \, : \, \frG^{\scO} \raq \Setdel \, , \qquad x \mapsto N((\frG^{\scO})_{/x})
\end{equation*}
where $(\frG^{\scO})_{/x}$ is the under-category, and $N(-)$ is the nerve. Since the under-category has a final object, its nerve is contractible. Hence, the unique map $\cQ_0^{{\rm pre}} \ra \ast_{{\rm const}}$ is a pointwise homotopy equivalence. This part of the argument holds if $\frG^{\scO}$ is replaced by any other category.

Now we use the symmetric monoidal structure on $\frG^{\scO}$ to define the functor $ (\frG^{\scO})_{/x} \times (\frG^{\scO})_{/y} \ra (\frG^{\scO})_{/(x \otimes y)}$. Applying the nerve, we see that $\cQ_0^{{\rm pre}}$ has a lax symmetric monoidal structure.
Applying Proposition \ref{lax_assoc_operad}, we may take the associated operad $\cQ_0 := (\cQ_0^{{\rm pre}})^a$. Since the underlying $\frf^{\scO}$-module is unchanged, the unique map $\cQ_0 \ra \ast_{{\rm const}}$ is still a pointwise homotopy equivalence. Thus, this gives a resolution of $\ast_{{\rm const}}$ in $\OP(\frG^{\scO},\Setdel)$.

One can then use this to resolve $\bK_{{\rm const}} \in \OP(\frG^{\scO},\cC)$ for $\cC = \Ch(\bK)$. Namely, let $N_*(\bK[-]) : \Setdel \ra \Ch_{\geq 0}(\bK)$ be the Dold-Kan normalization functor, which has a lax symmetric monoidal structure via the Eilenberg-Zilber shuffle map. 
In particular, $N_*(\bK[\cQ_0^{{\rm pre}}]) : \frG^{\scO} \ra \Ch_{\geq 0}(\bK)$ is lax symmetric monoidal.
Applying Proposition \ref{lax_assoc_operad}, we then obtain a resolution of $\bK_{{\rm const}} \in \OP(\frG^{\scO},\cC)$.
Using this, one may define a homotopy cyclic coinvariant version of ${\rm colim}_{\frf^{\scO}}^{\gr}(\cP)$ for $\cP = \Endrd(\cA) \otimes \cS^{1,m}$, and a corresponding notion of pre-Calabi-Yau categories.

\brm
We do not know of an explicit resolution of $\bK_{{\rm const}}^c \in \COP(\frG^{\scO},\cC)$. However, we make some observations here. The first observation is that, in order to find a set-theoretic resolution, it suffices to find an object in $\COP(\frG^{\scO},\Set)$ whose underlying $\frf^{\scO}$-module have non-empty and free automorphism group action. Indeed, if this is done, then, as in the Barratt-Eccles model of the $E_{\infty}$-operad (see, {\it e.g.}, \cite{BF04}), one can post-compose it with the symmetric monoidal functor $E : \Set \ra \Setdel$, given by $E(X)_n = X^{n+1} = \Hom_{\Set}([n],X)$. It is a standard fact that $E(X)$ is contractible whenever $X$ is non-empty (because any choice of $x \in X$ gives an extra degeneracy of $E(X)$).

The second observation is that, once a set-theoretical resolution $\cQ_0^c \xra{\sim} \ast^c_{{\rm const}}$ in $\COP(\frG^{\scO},\Setdel)$ is found, one may apply $[\bK(\cQ_0^c),\cP]^{\gr}_{\frf^{\scO}}$ to each simplicial degree, and obtain a cosimplicial dg Lie algebra, to which we may take its homotopy limit.
\erm

\brm  \label{Fuk_remark}
The problem of forming the homotopy cyclic (co)invariant versions of pre-Calabi-Yau categories is related to the problem of constructing a pre-Calabi-Yau structure on the Fukaya category. On the one hand, while the counting of holomorphic disks should ``morally'' give elements in the right hand side of \eqref{Hoch_limt_End}, since these counting depends on perturbation data, technically they are not $C_p$-invariant, but should only be homotopy invariant in some sense. On the other hand, the Fukaya category is supposed to be defined over $\bZ$, so the appropriate notion of pre-Calabi-Yau structure that it carries should be applicable in all characteristics. The author thanks S. Ganatra for explaining him this point.
\erm

\end{document}